\documentclass[review,onefignum,onetabnum]{siamart171218} 

\pdfoutput=1

\usepackage{soul}
\usepackage{cite}

\usepackage{graphicx,subcaption, mwe,amsmath,amsfonts,latexsym,amssymb,amsmath,array,fancyhdr,lscape}

\usepackage{hyperref}
\hypersetup{
    colorlinks=false, 
    linktoc=all,     
    linkcolor=black,  
}
\usepackage{bm}
\usepackage{array, caption, tabularx, makecell, booktabs}
\usepackage{mathtools}
\usepackage{multirow}

\usepackage{calc}

\newcommand{\citeCount}[1]{}
\newcommand{\dt}{\Delta t}
\newcommand{\pd}[2]{\frac{\partial #1}{\partial #2}}
\newcommand{\pdn}[3]{\frac{\partial^#3 #1}{\partial #2^#3}}

\newcommand{\dd}[2]{\frac{d #1}{d #2}}
\newcommand{\ddn}[3]{\frac{d^#3 #1}{d #2^#3}}

\newcommand{\pdh}[2]{\frac{\partial_h #1}{\partial_h #2}}
\newcommand{\pdhn}[3]{\frac{\partial_h^#3 #1}{\partial_h #2^#3}}

\newcommand{\iv}{\mathbf{ i}}

\newcommand{\nv}{\mathbf{ n}}

\newcommand{\tv}{\mathbf{ t}}

\newcommand{\xv}{\mathbf{ x}}

\newcommand{\Ac}{{\mathcal A}}
\newcommand{\Bc}{{\mathcal B}}

\newcommand{\Kc}{{\mathcal K}}

\newcommand{\Gc}{{\mathcal G}}

\usepackage{tikz}

\newlength{\tfwidth}
\newlength{\tfheight}
\newlength{\tfxa}
\newlength{\tfxb}
\newlength{\tfya}
\newlength{\tfyb}
\newcommand{\trimFigWithBox}[6]{%
\setlength\fboxsep{0pt}%
\setlength\fboxrule{1.0pt}
\fbox{\includegraphics[width=#2, clip, trim=#3 #4 #5 #6]{#1}}%
}
\newcommand{\trimFigNoBox}[6]{%
\setlength\fboxsep{1pt}
\setlength\fboxrule{0.0pt}
\fbox{\includegraphics[width=#2, clip, trim=#3 #4 #5 #6]{#1}}%
}
\newcommand{\trimFigHeightWithBox}[6]{%
\setlength\fboxsep{0pt}%
\setlength\fboxrule{1.0pt}
\fbox{\includegraphics[height=#2, clip, trim=#3 #4 #5 #6]{#1}}%
}
\newcommand{\trimFigHeightNoBox}[6]{%
\setlength\fboxsep{1pt}
\setlength\fboxrule{0.0pt}
\fbox{\includegraphics[height=#2, clip, trim=#3 #4 #5 #6]{#1}}%
}

\newsavebox\figBox

\newcommand{\trimw}[6]{%
\sbox\figBox{\includegraphics{#1}}
\setlength{\tfwidth}{\the\wd\figBox}
\setlength{\tfheight}{\the\ht\figBox}
\setlength{\tfxa}{\tfwidth*\real{#3}}%
\setlength{\tfxb}{\tfwidth*\real{#4}}%
\setlength{\tfya}{\tfheight*\real{#5}}%
\setlength{\tfyb}{\tfheight*\real{#6}}%
\trimFigNoBox{#1}{#2}{\tfxa}{\tfya}{\tfxb}{\tfyb}%
}

\newcommand{\trimwb}[6]{%

\sbox\figBox{\includegraphics{#1}}
\setlength{\tfwidth}{\the\wd\figBox}
\setlength{\tfheight}{\the\ht\figBox}
\setlength{\tfxa}{\tfwidth*\real{#3}}%
\setlength{\tfxb}{\tfwidth*\real{#4}}%
\setlength{\tfya}{\tfheight*\real{#5}}%
\setlength{\tfyb}{\tfheight*\real{#6}}%
\trimFigWithBox{#1}{#2}{\tfxa}{\tfya}{\tfxb}{\tfyb}%
}

\newcommand{\trimh}[6]{%
\sbox\figBox{\includegraphics{#1}}
\setlength{\tfwidth}{\the\wd\figBox}
\setlength{\tfheight}{\the\ht\figBox}
\setlength{\tfxa}{\tfwidth*\real{#3}}%
\setlength{\tfxb}{\tfwidth*\real{#4}}%
\setlength{\tfya}{\tfheight*\real{#5}}%
\setlength{\tfyb}{\tfheight*\real{#6}}%
\trimFigHeightNoBox{#1}{#2}{\tfxa}{\tfya}{\tfxb}{\tfyb}%
}

\newcommand{\trimhb}[6]{%

\sbox\figBox{\includegraphics{#1}}
\setlength{\tfwidth}{\the\wd\figBox}
\setlength{\tfheight}{\the\ht\figBox}
\setlength{\tfxa}{\tfwidth*\real{#3}}%
\setlength{\tfxb}{\tfwidth*\real{#4}}%
\setlength{\tfya}{\tfheight*\real{#5}}%
\setlength{\tfyb}{\tfheight*\real{#6}}%
\trimFigHeightWithBox{#1}{#2}{\tfxa}{\tfya}{\tfxb}{\tfyb}%
}
\usepackage{algorithm}
\usepackage{pdfpages} 

\begin{document}

\title{
  {Stable and accurate  numerical methods for generalized Kirchhoff-Love plates}
  \thanks{Submitted to the editors DATE.
    \funding{This research was supported by   the Louisiana Board of Regents Support Fund under contract No. LEQSF(2018-21)-RD-A-23.}}}

\author{Duong~T.~A.~Nguyen\thanks{Department of Mathematics, University of Louisiana at Lafayette, Lafayette, LA 70504, USA.
(\email{duong.nguyen1@louisiana.edu}).
}
\and Longfei~Li\thanks{Corresponding Author. Department of Mathematics, University of Louisiana at Lafayette, Lafayette, LA 70504, USA.
(\email{longfei.li@louisiana.edu})}
\and Hangjie~Ji\thanks{Department of Mathematics, University of California Los Angeles, Los Angeles, CA 90095, USA. (\email{hangjie@math.ucla.edu}).}
}
\maketitle

\begin{abstract}

Efficient and accurate numerical algorithms are developed to solve a generalized Kirchhoff-Love plate model subject to three common physical boundary conditions: (i) clamped; (ii) simply supported; and (iii) free. We solve the model  equation by   discretizing the spatial derivatives using  second-order finite-difference schemes, and then advancing the semi-discrete problem in time with either an explicit predictor-corrector  or an implicit Newmark-Beta time-stepping algorithm. Stability analysis is conducted for the schemes and the results are used to determine stable time steps in practice.
 A series of carefully  chosen  test problems  are solved  to demonstrate the properties and applications  of our numerical approaches.
The numerical results confirm the stability and  2nd-order accuracy of the algorithms, and  are also comparable  with  experiments for similar  thin plates. As an application, we  illustrate a strategy  to  identify the natural frequencies of  a plate using our numerical methods in conjunction with a fast Fourier transformation (FFT) power spectrum analysis  of the  computed data. Then we  take advantage of one of the  computed  natural frequencies to simulate the interesting physical phenomena known as resonance and  beat  for a generalized  Kirchhoff-Love  plate. 


\end{abstract}

\begin{keywords}
 thin plates,  Kirchhoff-Love theory, finite difference method, predictor-corrector method, Newmark-Beta scheme, eigenvalues and eigenmodes, resonance.
\end{keywords}

\begin{AMS}
  65M06 , 65M12, 74S20
\end{AMS}





\tableofcontents

\section{Introduction}

Thin-walled elastic solids, often referred to as plates or shells,  are  ubiquitous  in  engineerings and applied sciences.  Examples of plates or shells  can be found in many common mechanical and biological structures such as  dome-shaped stadium  rooftops,  airplane fuselages, vessel walls and  aortic valves, etc.  Adequately  understanding  the intrinsic properties  of plates (shells) is crucial for the various applications involving these structures. Therefore,  the investigation of mathematically modeling plate-like structures    and the subsequent  development of numerical approximations for their solutions have long been  active areas of research. Noting that the difference between a plate and a shell lies in its precast  stress-free shape, which is flat for a plate and curved for a shell.

To   study  plate structures analytically,  numerous   theories  have been developed over the years aiming at predicting the various key physical characteristics; see  \cite{Reissner76,TimoshenkoWoinowsky59,Leissa69,Love1888,KoiterSimmonds73} and the references therein. The classical Kirchhoff-Love  plate theory, which  was developed way back in  1888 under the assumptions  that the thickness of the plate remains fixed and any straight lines normal to the reference surface remain straight and normal to the reference surface after deformation,  captures the bending dynamics of a plate in response to a transverse load  and  determines the propagation of waves in the plate \cite{Love1888}. As an extension to the Kirchhoff-Love model,  the Mindlin-Reissner plate theory  takes a first-order shear deformation  into account and no longer assumes that straight lines normal to the reference surface remain normal during a deformation \cite{Mindlin51}. As is reviewed in \cite{KoiterSimmonds73},  there are also many other plate theories that are able to  describe more sophisticated nonlinear physical phenomena, which make them viable choices for modeling complicated engineering applications. For example, the Koiter shell theory \cite{Koiter60} and its recent variant that  incorporates viscoelasticity
\cite{CanicEtal06} are often used in biomedical engineering to model    artery walls. 

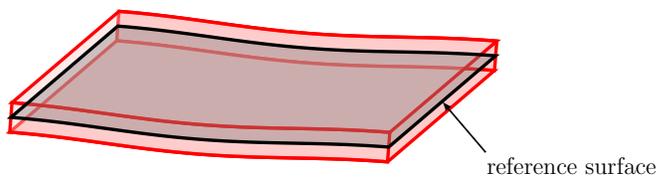
\begin{figure}[h]
  \centering
  \resizebox{9cm}{!}{ 
    %

%

{
\def\xa{-.5}
\def\xb{12.5}
\def\xs{.3}
\def\ys{.3}
\def\lb{8}
\def\a{.075}%
\def\b{45}
\def\hsh{.5}
\def\extra{1.}
\def\yTop{2.5}
\def\yBot{-2.75}
%
\newcommand{\textFont}{\normalss}
\begin{tikzpicture}[scale=1]
\useasboundingbox (\xa,.25) rectangle (12.5,5.25);  
\begin{scope}[xshift=0cm,yshift=\yTop cm]
%
  \begin{scope}[xshift=2.5cm,rotate=-4.5]
    %
    %
    \fill[fill=red!20,draw=red,line width=2pt,xshift=.3cm,yshift=1.3cm] 
           plot[samples=100, domain=0.:\lb] (\x, {(0.25*\hsh+\a*cos(\b*\x))}) --
           plot[samples=100, domain=\lb:0] (\x, {(-\hsh+\a*cos(\b*\x))}) -- cycle;
           %
    %
    \fill[fill=red!20,draw=red,line width=2pt,xshift=.3cm,yshift=1.3cm] 
           plot[samples=100, domain=0.:\lb] (\x, {(-\hsh+\a*cos(\b*\x))}) --
           plot[samples=100, domain=\lb:0] (\x-7*\xs, {(-7*\ys-\hsh+\a*cos(\b*\x))}) -- cycle;           
    \fill[fill=red!20,draw=red,line width=2pt,xshift=.3cm,yshift=1.3cm,fill opacity=0.4] 
           plot[samples=100, domain=0.:\lb] (\x-7*\xs, {(-7*\ys+0.25*\hsh+\a*cos(\b*\x))}) --
           plot[samples=100, domain=\lb:0] (\x-7*\xs, {(-7*\ys-\hsh+\a*cos(\b*\x))}) -- cycle;
    \fill[fill=red!20,draw=red,line width=2pt,xshift=.3cm,yshift=1.3cm,fill opacity=0.4] 
           plot[samples=100, domain=0.:\lb] (\x-7*\xs, {(-7*\ys+0.25*\hsh+\a*cos(\b*\x))}) --
           plot[samples=100, domain=\lb:0]  (\x, {(0.25*\hsh+\a*cos(\b*\x))}) -- cycle;
     \fill[fill=black!50,draw=black,line width=2pt,xshift=.3cm,yshift=1.3cm,fill opacity=0.4] 
           plot[samples=100, domain=0.:\lb] (\x-7*\xs, {(-7*\ys-0.5*0.75*\hsh+\a*cos(\b*\x))}) --
           plot[samples=100, domain=\lb:0]  (\x, {(-0.5*0.75*\hsh+\a*cos(\b*\x))}) -- cycle;

     \pgfmathparse{(\lb-7*\xs+\lb)/2}
     \let\xx\pgfmathresult
     
     \pgfmathparse{((-7*\ys-0.5*0.75*\hsh+\a*cos(\b*\lb))+(-0.5*0.75*\hsh+\a*cos(\b*\lb)))/2}
     \let\yy\pgfmathresult

    \draw[black] ({\xx+1},{\yy-1}) node[anchor=north west,xshift=.3cm,yshift=1.3cm] {{\Large reference surface}};
    \draw[-latex,black,line width=1.pt,xshift=.3cm,yshift=1.3cm] ({\xx+1},{\yy-1}) -- (\xx,\yy);


     
    

  \end{scope}
%
\end{scope}
%
%
\end{tikzpicture}
}
    }
  \caption{Cartoon illustration of a deformed  thin plate and its reference surface.}\label{fig:shellCartoon}
\end{figure}

These plate theories are in general  derived by utilizing the disparity in the  length scales of the  thin structures, and significantly  reduce the complexity of the    three-dimensional (3D)  continuum mechanics  problem to a  two-dimensional one (2D). The governing equations of a plate theory  typically deal with variables defined only on a reference surface  that resides on a 2D domain; for  an  isotropic and homogeneous plate, its middle (or center) surface is  used as  reference.  See Figure~\ref{fig:shellCartoon}  for a schematic illustration of a 3D thin plate and its 2D reference surface.
Physical assumptions of the underlining plate theories also provide means of calculating the load-carrying and deflection characteristics of the original thin-walled structures; and therefore, the complete deformation and stress fields of a 3D thin   structure can be inferred from the  solution of its  reference surface. It is generally expected that the thinner  the  structure, the more accurate the  plate theory.

From both  the analytical and numerical points of view, 2D plate theories are  immensely more tractable than  3D  solid mechanical models. The  plate theories are especially appealing to researchers exploring multi-physical problems, such as  fluid-structure interaction (FSI) problems   involving thin-walled elastic structures \cite{fis2014,LiHenshaw2016,CanicEtal06},  whereby multiple physical subproblems  are dealt with simultaneously.  Although greatly simplified from   the  full 3D continuum mechanics problem,   governing equations for  plates  are still too complicated to be solved analytically, except for a limited number of cases with simple specifications \cite{Rudolph-2004}.  Efficient and accurate numerical approximations for the solutions are therefore of greater interest in practice. However, due to  numerical challenges posed by the high-order spatial  derivatives that are associated with the  bending effect of plates, the development of stable and accurate numerical methods for solving plate equations is non-trivial.   Many numerical approaches have been developed for solving various plate models  based on common discretization methods such as finite difference \cite{bilbao2008family,JiLi2019}, finite element  (FEM) \cite{Jean-Loui-1982,Adna-1993,Eugeni-2000,daVeiga-2007,Bischoff-2004,Motasoares-2006,Perotti-2013,Huang-2011,BecacheEtal2004}, and  boundary element (BEM) \cite{AFrangi-1999}, to name just  a few. More recently,  new computational methods developed  from the  isogeometric analysis have also emerged; see  \cite{BensonEtal2010} for an example of solving  the Reissner-Mindlin shell using the  isogeometric analysis.

In this paper,
we  present  accurate and efficient numerical approximations of   a generalized Kirchhoff-Love model that  incorporates  additional important physics, such as  linear restoration,   tension and visco-elasticity, for  isotropic and homogeneous thin plates.  The  generalized model equation is spatially discretized with a standard second-order accurate finite difference method and integrated in time using either an explicit predictor-corrector or an implicit  time-stepping scheme.   Stability analysis is performed and the  results are utilized to determine  stable time steps for the proposed numerical schemes. Stable and accurate numerical  boundary conditions are also investigated for the most common physical boundary conditions (i.e., clamped, simply supported  and free).   Carefully designed   test problems are solved  using  all the proposed schemes for numerical validations. Interesting applications using the numerical methods are also discussed.

 The remainder of the paper is organized as follows. In Section \ref{sec:governingEquations}, we present  the governing equation and its boundary conditions for a  generalized Kirchhoff-Love model. The numerical algorithms  for solving  the model equation  are discussed in Section \ref{sec:numericalMethods}. We analyze  the stability of the numerical schemes and lay out a strategy for determining stable time steps for the  algorithms in Section~\ref{sec:analysis}. Numerical results that provide verification of the stability and accuracy of the schemes,  as well as cross-validation with  experiments,  are presented  in Section~\ref{sec:results}.  Finally, concluding remarks  are made   in Section~\ref{sec:conclusions}.

\section{Governing equations}\label{sec:governingEquations}
The classical Kirchhoff-Love plate model  concerns the small  deflection of thin plates that is used as a simplified theory  of solid mechanics  to  determine the stresses and deformations in the plates subject to external forcings.     The governing  equation   for  an isotropic and homogeneous plate is a  single  time-dependent biharmonic partial differential  equation (PDE)  for the transverse displacement  of the plate's middle surface. It  is   derived by   balancing the external loads with  the internal  bending force  that tends to restore the plate to its stress-free state.
In this paper, we consider a 
plate model that is generalized from the  classical Kirchhoff-Love equation by including  additional terms to account for more physical effects, such as   linear restoration, tension and visco-elasticity.

To be specific, this   work concerns  developing numerical algorithms for solving the following generalized Kirchhoff-Love  model for  an isotropic and homogeneous plate with constant thickness $h$, 
\begin{equation} \label{eq:generalizedKLPlate}
{\rho}{h}\pdn{w}{t}{2}=-{K}_0w+{T}\nabla^2w-{D} \nabla^4w-{K}_1\pd{w}{t}+{T}_1\nabla^2\pd{w}{t}+F(\xv,t),
\end{equation}
where $w(\xv,t)$ with $\xv\in\Omega\subset\mathbb{R}^2$ is   the transverse displacement of the middle surface subject to some given body force $F$. Here,  ${\rho}$ denotes density, ${K}_0$ is the linear stiffness coefficient that acts as a linear restoring force, ${T}$ is the tension coefficient, and ${D}={E}{h}^3/(12(1-\nu))$  represents the flexural rigidity  with $\nu$  and $E $ being the Poisson's ratio and  Young's modulus, respectively. The term with coefficient ${K}_1$ is a linear damping term, while the term with coefficient ${T}_1$ is a visco-elastic damping term that tends to smooth high-frequency oscillations in space. Noting that the visco-elastic damping is often added to model vascular structures in haemodynamics \cite{CanicEtal06}.

On the boundary, one of   the following   physical boundary conditions is imposed; namely, 
for $\xv\in\partial\Omega$,  we have
{\small
\begin{align}
  &\text{clamped:} & w=0, && \pd{ w}{\nv}=0; \label{eq:clampedBC}\\
  &\text{supported:}  &  w=0, && \pdn{ w}{\nv}{2}+\nu\pdn{ w}{\tv}{2}=0; \label{eq:supportedBC}\\
  &\text{free:}  &   \pdn{ w}{\nv}{2}+\nu\pdn{ w}{\tv}{2}=0, &&   \pd{}{\nv}\left[\pdn{ w}{\nv}{2}+\left(\nu-2\right)\pdn{ w}{\tv}{2}\right]=0, \label{eq:freeBC}
\end{align}
}
where  $\partial/\partial{\nv}$ and  $\partial/\partial{\tv}$ are the normal and tangential derivatives defined on  the boundary of the domain. It is important to point out that, for a rectangular plate,  the free boundary conditions \eqref{eq:freeBC} must be complemented by a corner condition that imposes zero forcing \cite{bilbao2008family}; in other words,  we also impose  $\partial^2w/\partial x\partial y=0$ at the corners of a rectangular plate. 
Note that  for notational brevity   the functional dependence on $(\xv,t)$ has been suppressed  in the statement  of the boundary conditions.

Appropriate initial conditions need to be specified to complete the statement of the governing equations. Specifically, we assign
\begin{equation}\label{eq:IC}
w(\xv,0)=w_0(\xv) ~\text{and}~\pd{w}{t}(\xv,0)=v_0(\xv)
\end{equation}
as the initial conditions with $w_0(\xv)$ and $v_0(\xv)$ representing two given functions that prescribe  the plate's  initial displacement and velocity.


\section{Numerical methods} \label{sec:numericalMethods} 
In this section, we present the numerical approaches to solve the governing equation \eqref{eq:generalizedKLPlate} subject to the boundary conditions \eqref{eq:clampedBC}--\eqref{eq:freeBC}. Standard centered  finite difference methods of second-order accuracy  are used for the discretization of all the spatial derivatives, and then the resulted  semi-discrete equations are integrated in time using an appropriate  time-stepping scheme.

Let $\Omega_h$ denote a mesh covering the domain $\Omega$, and let $\xv_\iv\in \Omega_h$ denote  the coordinates of a grid point with multi-index  $\iv=(i_1,i_2)$. The time-dependent  grid function that approximates  the displacement on the mesh is given by $w_\iv(t)\approx w(\xv_{\iv},t)$. Similarly, $F_\iv$ is used to denote the given forcing function evaluated at $\xv_\iv$; namely, $F_\iv(t)=F(\xv_\iv,t)$.   We spatially discretize the governing equation and its boundary conditions by replacing the differential operators with the corresponding finite-difference operators (distinguished with a subscript $h$) to derive the  semi-discrete equations,
\begin{equation}\label{eq:discreteGeneralizedKLPlate}
{\rho}{h}\ddn{w_{\iv}}{t}{2}=-{K}_0w_{\iv}+{T}\nabla_h^2w_{\iv}-{D} \nabla_h^4w_{\iv}-{K}_1\dd{w_{\iv}}{t}+{T}_1\nabla_h^2\dd{w_{\iv}}{t}+F_{\iv}, ~\forall \xv_\iv\in\Omega_h,
\end{equation}
as well as the following discrete boundary conditions.    For $\forall \xv_{\iv_b}\in\partial\Omega_h$, the numerical boundary conditions are given by 
{\small
\begin{align}
  & \text{clamped:}  & w_{\iv_b}=0, && \pdh{ w_{\iv_b}}{\nv}=0;  \label{eq:discreteClampedBC}\\
  &\text{supported:}  &  w_{\iv_b}=0, && \pdhn{ w_{\iv_b}}{\nv}{2}+\nu\pdhn{ w_{\iv_b}}{\tv}{2}=0; \label{eq:discreteSupportedBC}\\
  &\text{free:}  &  \pdhn{ w_{\iv_b}}{\nv}{2}+\nu\pdhn{ w_{\iv_b}}{\tv}{2}=0, &&   \pdh{}{\nv}\left[\pdhn{ w_{\iv_b}}{\nv}{2}+\left(\nu-2\right)\pdhn{ w_{\iv_b}}{\tv}{2}\right]=0. \label{eq:discreteFreeBC}
\end{align}
}

For numerical purposes, we  rewrite \eqref{eq:discreteGeneralizedKLPlate}  into a system of first-order ODEs.  
If we  denote $v_{\iv}$ and $a_{\iv}$  the   numerical approximations of the velocity and acceleration at grid point $\xv_\iv$, equation \eqref{eq:discreteGeneralizedKLPlate} can thus be conveniently written as
\begin{align}
    &\dd{w_\iv}{t}(t)= v_{\iv}(t),\label{eq:ODEw}\\
  &\dd{v_\iv}{t}(t)= a_{\iv}(t),\label{eq:ODEv}\\
  &{\rho}{h}a_\iv(t)=-{\Kc_h} w_{\iv}(t)-{\Bc_h}v_{\iv}(t)+F_{\iv}(t),\label{eq:ODEa}
\end{align}
where the operators for the internal forces $\Kc_h$ and the damping forces $\Bc_h$ are introduced below    to simplify the notations;
\begin{equation}\label{eq:KBOperators}
\Kc_h = {K}_0-{T}\nabla_h^2+{D} \nabla_h^4\quad\text{and}\quad
\Bc_h={K}_1- {T}_1\nabla_h^2.
\end{equation}

In this paper,  two time-stepping methods  are considered to advance the ODE system  in time. In particular, one of the methods is  an explicit predictor-corrector scheme that consists of a second-order Adams-Bashforth (AB2) predictor and a second-order Adams-Moulton (AM2) corrector, while  the other one  is an implicit   Newmark-Beta scheme of second-order  accuracy \cite{Newmark59}. We refer to the former scheme as  PC22 scheme and the latter one as    NB2 scheme for short.

To simplify the discussion, the algorithms are developed for a fixed time-step $\dt$ so that $t_n =n\dt$.
Let the numerical solutions of  \eqref{eq:ODEw}--\eqref{eq:ODEa}   at time $t_n$ be $w_\iv^n\approx w_\iv(t_n)$,  $v_\iv^n\approx v_\iv(t_n)$, and   $a_\iv^n\approx a_\iv(t_n)$, and denote $F_{\iv}^{n}=F(\xv_\iv,t_n)$. The goal of a  time-stepping algorithm is to determine the solutions at a new time given solutions at   previous time levels.

First, we  describe the PC22 scheme  in Algorithm~\ref{alg:pc22}. 
\begin{algorithm}[H]
 {\bf Input:} {solutions at two previous time levels; i.e.,   $(w_\iv^n,v_\iv^n,a_\iv^n)$ and $(w_\iv^{n-1},v_\iv^{n-1},a_\iv^{n-1})$ }\\
 {\bf Output:} {solutions at the new time level; i.e., $(w_\iv^{n+1},v_\iv^{n+1},a_\iv^{n+1})$  }\\
  {\bf Procedures:} \\
 {\em Stage I:  predict solutions using a second-order Adams-Bashforth (AB2) predictor}
 \begin{equation*}
   \forall \xv_\iv\in\Omega_h:\quad
   \begin{cases}
     w_\iv^p= w_\iv^n+\dt\left(\frac{3}{2}v_\iv^{n}-\frac{1}{2}v_\iv^{n-1}\right)\\
     v_\iv^p= v_\iv^n+\dt\left(\frac{3}{2}a_\iv^{n}-\frac{1}{2}a_\iv^{n-1}\right) \\
     a^{p}_\iv=\frac{1}{{\rho}{h}}\left(-{\Kc_h} w^{p}_{\iv}-{\Bc_h}v^{p}_{\iv}+F^{n+1}_{\iv}\right)
   \end{cases}
 \end{equation*}
 {\em Stage II:   correct solutions using a  second-order Adams-Molton (AM2) corrector}
  \begin{equation*}
   \forall \xv_\iv\in\Omega_h:\quad
   \begin{cases}
     w_\iv^{n+1}= w_\iv^n+\dt\left(\frac{1}{2}v_\iv^{n}+\frac{1}{2}v_\iv^{p}\right)\\
     v_\iv^{n+1}= v_\iv^n+\dt\left(\frac{1}{2}a_\iv^{n}+\frac{1}{2}a_\iv^{p}\right) \\
     a^{n+1}_\iv=\frac{1}{{\rho}{h}}\left(-{\Kc_h} w^{n+1}_{\iv}-{\Bc_h}v^{n+1}_{\iv}+F^{n+1}_{\iv}\right)
   \end{cases}
  \end{equation*}
  
      {\bf Remark}:{\em Boundary conditions are applied after both the predictor and corrector stages to fill in the solutions at ghost and/or boundary grid points. Note that numerical  boundary conditions for $v$ and $a$  are derived from those for  $w$  by taking the appropriate time derivatives of \eqref{eq:discreteClampedBC}--\eqref{eq:discreteFreeBC}.}
 \caption{PC22 time-stepping scheme}\label{alg:pc22}
\end{algorithm}

Second, we consider the Newmark-Beta scheme for solving our problem \eqref{eq:ODEw}--\eqref{eq:ODEa}. The so-called  Newmark-Beta scheme is a general procedure proposed by Newmark for the solution of  problems in structural dynamics \cite{Newmark59}. Given acceleration, the scheme updates the velocity and displacement by solving
\begin{equation} \label{eq:nbWV}
\begin{cases}
  w_{\iv}^{n+1}= w_{\iv}^{n}+\dt v_{\iv}^n +\frac{\dt^2}{2}\left[ (1-2\beta)a_{\iv}^{n}+2\beta a_{\iv}^{n+1} \right],  \\
  v_{\iv}^{n+1} = v_{\iv}^{n}+\dt \left[ (1-\gamma)a_{\iv}^{n}+\gamma a_{\iv}^{n+1} \right],
\end{cases}
\end{equation}
where the acceleration in our case is given by
\begin{equation} \label{eq:nbA}
{\rho}{h}a^{n+1}_\iv=-{\Kc_h} w_{\iv}^{n+1}-{\Bc_h}v_{\iv}^{n+1}+F_{\iv}^{n+1}.
\end{equation}
We note that the scheme  is unconditionally stable if $1/2\leq \gamma\leq 2\beta$, whereas it is conditionally stable if $\gamma>\max\{1/2,2\beta\}$.

Instead of solving the above implicit system for $w_\iv^{n+1}$, $v_\iv^{n+1}$ and $a_\iv^{n+1}$ all  at the same time, we use \eqref{eq:nbWV} to eliminate $w_\iv^{n+1}$ and  $v_\iv^{n+1}$ in \eqref{eq:nbA}, and then solve a smaller system for $a_\iv^{n+1}$ only.  The complete algorithm for this scheme is summarized in Algorithm~\ref{alg:nb2}.
\begin{algorithm}[H]
 {\bf Input:} {solutions at the previous time level; i.e.,   $(w_\iv^n,v_\iv^n,a_\iv^n)$  }\\
 {\bf Output:} {solutions at the new time level; i.e., $(w_\iv^{n+1},v_\iv^{n+1},a_\iv^{n+1})$  }\\
 {\bf Procedures:} \\
 {\em Stage I.  compute a first-order prediction for displacement and velocity}
 \begin{equation*} 
   \forall \xv_\iv\in\Omega_h:\quad
\begin{cases}
  w_{\iv}^{p}= w_{\iv}^{n}+\dt v_{\iv}^n +\frac{\dt^2}{2} (1-2\beta)a_{\iv}^{n}  \\
  v_{\iv}^{p} = v_{\iv}^{n}+\dt  (1-\gamma)a_{\iv}^{n}
\end{cases}
 \end{equation*}
{\em Stage II.   solve a system of equations for acceleration at $t_{n+1}$}
\begin{equation*}
 \forall \xv_\iv\in\Omega_h:\quad  \left( {\rho}{h}+\beta \dt^2\Kc_h+\gamma\dt\Bc\right)a^{n+1}_\iv=-{\Kc_h} w_{\iv}^{p}-{\Bc_h}v_{\iv}^{p}+F_{\iv}^{n+1}
\end{equation*}

{\em Stage III.  solve for displacement and velocity  at $t_{n+1}$} \\
$$
\text{
  compute $w_\iv^{n+1}$ and $v_\iv^{n+1}$ explicitly  from \eqref{eq:nbWV}
}
$$

{\bf Remark}: {\em In this paper, we set $\beta = 1/4$ and $\gamma = 1/2$. With this choice of parameters,  the scheme  is  second-order accurate and unconditionally stable.  We also note that boundary conditions are applied after stages I and III to fill in the solutions of $w$ and $v$ at ghost and/or boundary grid points. For stage II,  equations for acceleration at ghost and boundary nodes are replaced with boundary  conditions.} 
 \caption{NB2 time-stepping scheme}\label{alg:nb2}
\end{algorithm}

\section{Stability analysis and time step determination}\label{sec:analysis}
We study the  stability of the schemes and use the analytical results to determine  stable time steps in practical computations. As is already pointed out in  \cite{Newmark59} that the implicit NB2  time-stepping scheme is   unconditionally stable, the focus of the stability analysis here is on the explicit PC22 scheme.

\subsection{Stability of the PC22 scheme}

\begin{figure}[h] 
    \centering
    \includegraphics[width=2.in,height=2.in]{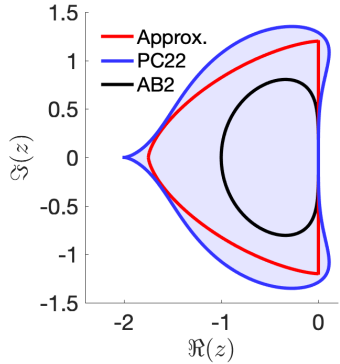} 
  \caption{Regions of absolute stability for the PC22 time-step scheme and  the scheme with an  AB2 predictor only. The approximated region of stability using a  half super-ellipse is also depicted in the plot. Here, $z=\lambda\dt$ with $\lambda$ being the time-stepping eigenvalue and $\dt$ representing the time step. $\Re(z)$ and $\Im(z)$ represent the real and imaginary parts of the complex number $z$.}
    \label{fig:AbsoluteStabilityRegion}
\end{figure}

Applying the PC22 scheme to the  Dahlquist test equation $\eta'=\lambda \eta$ leads  to the characteristic polynomial for a complex-valued amplification factor $\zeta$. Letting $z=\lambda\dt$, the roots of the characteristic equation are found to be
\begin{equation}\label{eq:RootsOfCharacteristicEqnPC22}
\zeta(z)=\dfrac{1}{2}\left(1+z+\frac{3}{4}z^2\pm\sqrt{\left(1+z+\frac{3}{4}z^2\right)^2-z^2}\right).
\end{equation}
The region of absolute stability for the PC22 time-stepping scheme is the set of complex values $z$ for which the roots of the characteristic polynomial satisfy   $|\zeta(z)| \leq 1$.

The stability region can be used to find the  time step restriction  for a typical problem. However, it is not straightforward to obtain a stable time step by  solving the inequality $|\zeta(z)| \leq 1$  directly from \eqref{eq:RootsOfCharacteristicEqnPC22}, so a half super-ellipse is introduced as an approximation of the stability region. To be specific, we  define the half super-ellipse by 
\begin{equation}\label{eq:approxStabilityRegion}
  \left|\frac{\Re(z)}{a} \right|^n+\left|\dfrac{\Im(z)}{b} \right|^n \leq 1 \quad\text{and}\quad \Re(z) \leq 0,
  \end{equation}
where $\Re(z)$ and $\Im(z)$ denote the real and imaginary parts of $z$, respectively.  We want the  half super-ellipse  to be completely enclosed by the actual region of stability and to be as large as possible. Given a time-stepping eigenvalue $\lambda$, it is much easier to find a sufficient condition for stability  by  requiring  $\lambda\dt$ to be  inside the approximated region defined by \eqref{eq:approxStabilityRegion}.
It is found that the above  half super-ellipse  makes  a good approximation for the stability region of  the PC22  time-stepping scheme by setting  $a=1.75,b=1.2$ and $n=1.5$.

The region of absolute stability for the PC22  time-stepping scheme, together with the approximated region, is shown in Figure~\ref{fig:AbsoluteStabilityRegion}. For comparison purposes, we also plot the stability region for the scheme  that only uses an AB2 predictor in the same  figure.  We can see that  by including  a corrector step the PC22 scheme has a much  larger stability region than the predictor alone, and the stability region includes the imaginary axis so that the scheme can be used for problems with no dissipations. From the plot, we can also see that the half super-ellipse that is chosen to be an approximation  fits perfectly inside the original stability region for the PC22 scheme. 
\subsection{Time-step determination}
A strategy for the determination of  stable time steps to be used in   Algorithms~\ref{alg:pc22} \& \ref{alg:nb2} is outlined here. We first transform   the semi-discrete problem \eqref{eq:ODEw}--\eqref{eq:ODEa} into Fourier Space, and then derive a stable time step   by imposing the condition that  the product of the time step  and any eigenvalue of the Fourier transformation of the difference operators lies  inside the stability region of a particular time-stepping method for all wave numbers. 

For simplicity of presentation, we assume the plate resides on a unit square domain ($\Omega=[0,1]\times[0,1]$), and the solution is 1-periodic in both $x$ and $y$ directions. Results on a more general domain can be readily obtained by mapping the general domain to a unit square.  After Fourier transforming the homogeneous version of   equations  \eqref{eq:ODEw}--\eqref{eq:ODEa} (i.e., assuming $F_\iv(t)\equiv0$), an  ODE system  for the transformed variables ($\hat{w}$, $\hat{v}$) is derived with  $\omega=(\omega_x,\omega_y)$ denoting any  wavenumber pair, 
\begin{equation}\label{m1}
\begin{bmatrix} 
\hat{w}(\omega,t)\\
\hat{v}(\omega,t)
\end{bmatrix}_t =
{\hat{Q}(\omega)}
\begin{bmatrix} 
\hat{w}(\omega,t)\\
\hat{v}(\omega,t)
\end{bmatrix},
~~\text{where}~~
\hat{Q}(\omega)= \begin{bmatrix} 
0 & 1 \\
-\hat{\Kc}(\omega) & -\hat{\Bc}(\omega) 
\end{bmatrix}.
\end{equation}
Here $\hat{\Kc}$ and $\hat{\Bc}$ are the Fourier transformations of the difference operators  $\Kc_h/(\rho h)$ and $\Bc_h/(\rho h)$, respectively. Let $k_{\omega_x}={2\sin({\omega_xh_x}/{2})}/{h_x}$ and $k_{\omega_y}={2\sin({\omega_yh_y}/{2}})/{h_y}$, where $h_x$ and $h_y$ are the grid spacings in the corresponding directions; then we have
\begin{alignat*}{2}
\hat{\Kc}(\omega) &=\frac{1}{\rho{h}}\left[K_0+T\left(k^2_{\omega_x}+k^2_{\omega_y}\right)+D\left(k^4_{\omega_x}+2 k^2_{\omega_x}k^2_{\omega_y}+ k^4_{\omega_y}\right)\right],\\
\hat{\Bc}(\omega) &=\frac{1}{\rho{h}}\left[K_1+T_1\left(k^2_{\omega_x}+k^2_{\omega_y}\right)\right].
\end{alignat*}
Noting that both $\hat{\Kc}$ and $\hat{\Bc}$ are non-negative for any $\omega$. In the  analysis to follow,  their   maximum values denoted by $\hat{\Kc}_M$ and $\hat{\Bc}_M$ are of interest, which  are attained when $\omega_xh_x=n\pi$ and $\omega_yh_y=m\pi$ ($n,m\in \mathbb{Z}$);  namely,
\begin{alignat}{2}
\hat{\Kc}_M&=\frac{1}{\rho{h}}\left[K_0+4{T}\left(\frac{1}{h_x^2}+\frac{1}{h_y^2}\right)+16D\left(\frac{1}{h_x^2}+\frac{1}{h_y^2}\right)^2\right] ,\label{eq:KM}\\
\hat{\Bc}_M&=\frac{1}{\rho{h}}\left[K_1+4T_1\left(\frac{1}{h_x^2}+\frac{1}{h_y^2}\right)\right] .\label{eq:BM}
\end{alignat}

A numerical method is stable provided all  the eigenvalues of $\hat{Q}(w) \dt$ lie within the stability region of the time-stepping method. The eigenvalues of the coefficient matrix $\hat{Q}(w) $ for the problem  \eqref{m1} are 
\begin{equation}\label{eq:eigenValueStability}
\hat{\lambda}(\omega)=-\frac{\hat{\Bc}(\omega)}{2}\pm  \sqrt{\left(\frac{\hat{\Bc}(\omega)}{2}\right)^2-\hat{\Kc}(\omega)}.
\end{equation}
For stability analysis, it suffices to consider the eigenvalue with the largest possible magnitude denoted by $\hat{\lambda}_M$. To find  $\hat{\lambda}_M$, we consider the following situations.

\newcommand{\determinant}{\left({\hat{\Bc}(\omega)}/{2}\right)^2-\hat{\Kc}(\omega)}
 \textbf{Under-damped case.} If $\determinant<0$, we obtain complex eigenvalues,
  $$
  \hat{\lambda}(\omega)=-\frac{\hat{\Bc}(\omega)}{2}\pm  i\sqrt{\hat{\Kc}(\omega)-\left(\frac{\hat{\Bc}(\omega)}{2}\right)^2}.
  $$ 
  In this case, we may define
\begin{equation}\label{eq:lambdaMUnderDamped}
  \hat{\lambda}_M=-\frac{\hat{\Bc}_M}{2}\pm  i\sqrt{\hat{\Kc}_M-\left(\frac{\hat{\Bc}_M}{2}\right)^2},
\end{equation}
since 
$|\hat{\lambda}(\omega)|=\sqrt{\hat{\Kc}(\omega)} \leq \sqrt{\hat{\Kc}_M} =|\hat{\lambda}_M|$. Here $\hat{\Kc}_M$ and $\hat{\Bc}_M$ are the maximum values of  $\hat{\Kc}(\omega)$ and $\hat{\Bc}(\omega)$ that are given by \eqref{eq:KM} and \eqref{eq:BM}.

\textbf{Over-damped case.} If $\determinant>0$, the eigenvalues are real and are of the same form as \eqref{eq:eigenValueStability}.
In this case, we may define
\begin{equation}\label{eq:lambdaMOverDamped}
  \hat{\lambda}_M=-\hat{\Bc}_M.
\end{equation}
This is because
$$
|\hat{\lambda}(\omega)|\leq \frac{\hat{\Bc}(\omega)}{2}+  \sqrt{\left(\frac{\hat{\Bc}(\omega)}{2}\right)^2-\hat{\Kc}(\omega)}\leq \hat{\Bc}(\omega) \leq \hat{\Bc}_M=|\hat{\lambda}_M|.
$$

We note that $\hat{\lambda}_M$ introduced in \eqref{eq:lambdaMUnderDamped} and \eqref{eq:lambdaMOverDamped} represent the eigenvalues of the worst-case scenario for the under-damped  and over-damped cases, respectively. A sufficient condition that ensures stability for the PC22 scheme is found by letting  $z=\hat{\lambda}_M\Delta t$ lie  in  the approximated stability  region that is defined by  the half super-ellipse in \eqref{eq:approxStabilityRegion}.  Since the approximated stability  region is a subset of the actual one,  a time step  that is sufficient to guarantee the stability of Algorithm~\ref{alg:pc22} can be chosen as following,
\begin{equation} \label{eq:timeStep}
\Delta t ={C_{\text{sf}}}{\left(\left|\frac{\Re(\hat{\lambda}_M)}{a}\right|^n+\left|\frac{\Im(\hat{\lambda}_M)}{b}\right|^n\right)^{-1/n}},
\end{equation}
where $C_\text{sf}\in (0,1] $ is a stability  factor (sf) that multiplies an estimate of the largest stable time step based on the above analysis. Unless otherwise noted, we choose  $C_\text{sf} =0.9$ for the PC22 scheme throughout this paper.

In terms of the NB2 scheme, we know that it  is implicit in time and stable for any time step. However, for accuracy reasons, we choose its time step based on the  condition for the explicit PC22 time-stepping scheme \eqref{eq:timeStep}, but  with  a much  larger  stability  factor. Typically, we choose   $C_\text{sf}=90$ for  Algorithm~\ref{alg:nb2}.

\section{Numerical results}\label{sec:results}
We now present the results for a series of test problems  to demonstrate the properties and applications  of our numerical approaches.
Mesh refinement studies using  problems with known exact solutions are first considered to verify the stability and  accuracy of the schemes.    Free and forced vibrations of thin plates with various geometrical and physical configurations  are then  solved to  further demonstrate  the numerical properties of our schemes and to compare with existing results. In particular,  the simulation of one  test problem is   cross-validated with reported  experimental results.  As an application, we illustrate a strategy using our numerical methods, together with  fast Fourier transformation (FFT), to identify the natural frequencies of a plate, and then numerically investigate  the interesting physical  phenomena known as   resonance and beat.

\subsection{Method of manufactured solutions}
As a first test, we verify the accuracy and stability of the algorithms using the method of manufactured solutions by adding  forcing functions to  the PDE \eqref{eq:generalizedKLPlate}  and the boundary conditions \eqref{eq:clampedBC}--\eqref{eq:freeBC} so that a chosen function becomes an exact solution.  The   exact solution is  chosen to be
\begin{equation}\label{eq:manufacturedExact} 
w_e(x,y,t)=\sin^4\left( \pi (x+1)\right) \sin^4\left( \pi (y+1)\right) \cos(2\pi t).\\
\end{equation}
In order to validate the  algorithms on both Cartesian and curvilinear grids, we consider  a square plate ($\Omega_S=[-1,1]\times[-1,1]$) and an annular plate ($\Omega_A=\{\xv: 0.5 \leq |\xv| \leq 1\} $). Physical parameters of the governing equation are specified as   $\rho h=1,K_0=2,T=1,D=0.01,K_1=5,T_1=0.1$ and $\nu=0.1$ for this test.

\begin{figure}[h] 
\begin{center}
  \begin{subfigure}[b]{0.45\linewidth}
    \centering
    \includegraphics[width=1.\linewidth]{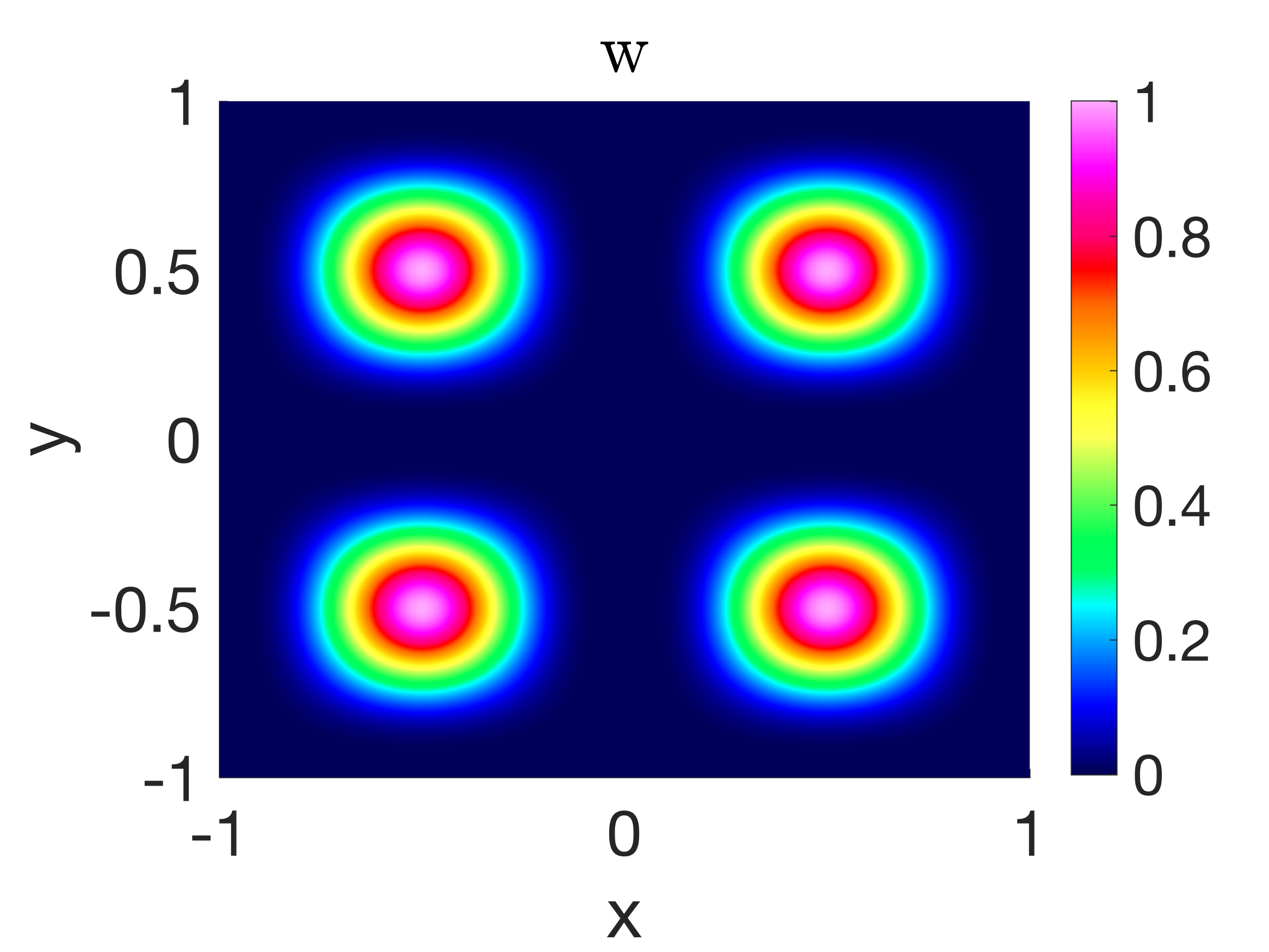}
      \end{subfigure} 
      \begin{subfigure}[b]{0.45\linewidth}
     \includegraphics[width=1.\linewidth]{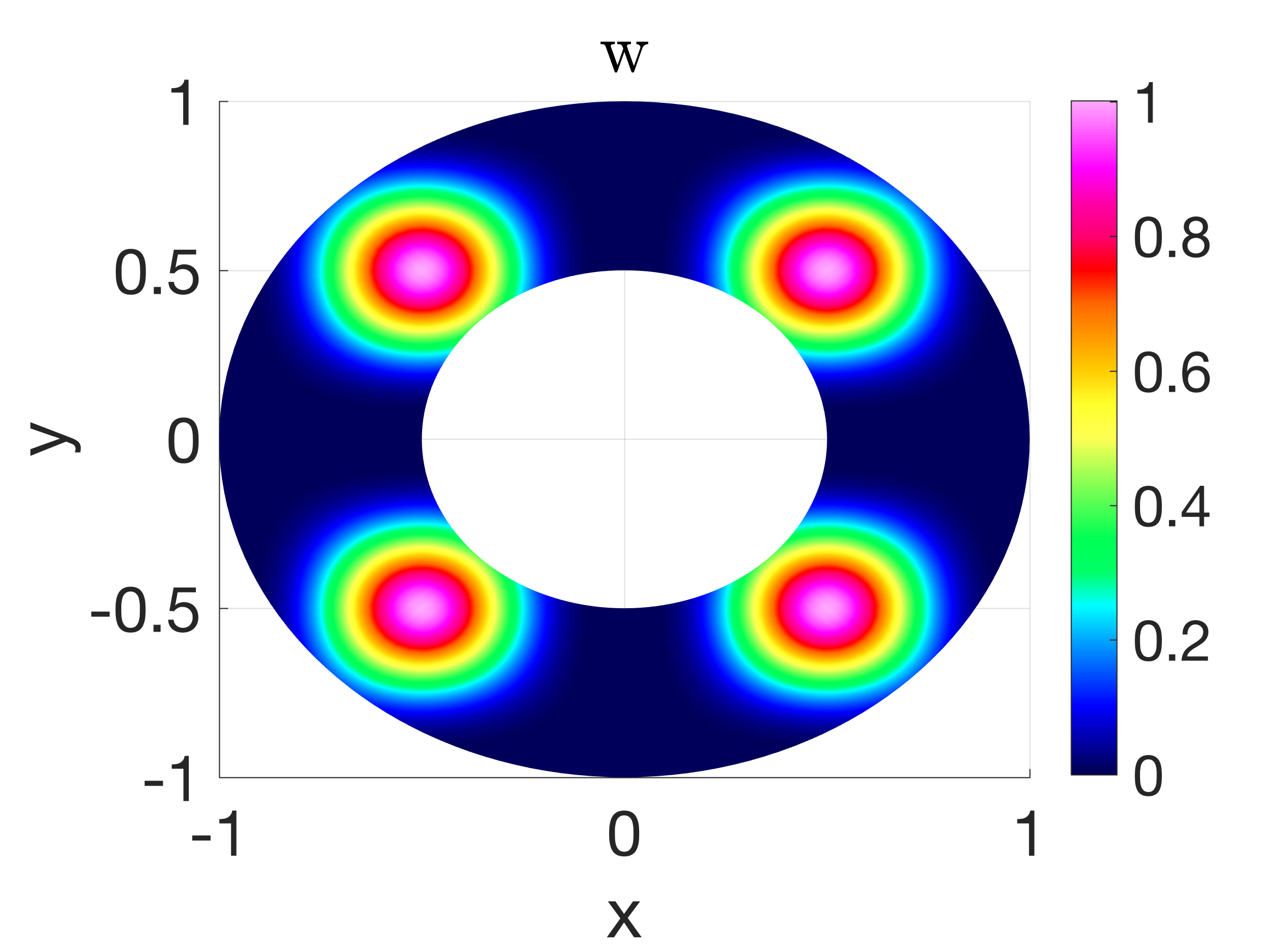}
  \end{subfigure} 
    \end{center}
  \caption{Contour plots of the computed displacement $w$ on grid $\Gc_{160}$ when $\text{t}=1$. The solution presented for the square plate are solved using the PC22 scheme subject to the free boundary conditions, while the plot for the annular plate is generated from the  NB2 scheme with simply supported boundary conditions. Results obtained using either algorithm are similar regardless of the  boundary conditions.
}
  \label{fig:ManufactureSol}
\end{figure}

\begin{figure}[h!] 
\centering
  \begin{subfigure}[b]{0.32\linewidth}
    \centering
    Clamped\\\vspace{0.1in}
    \includegraphics[width=1\linewidth]{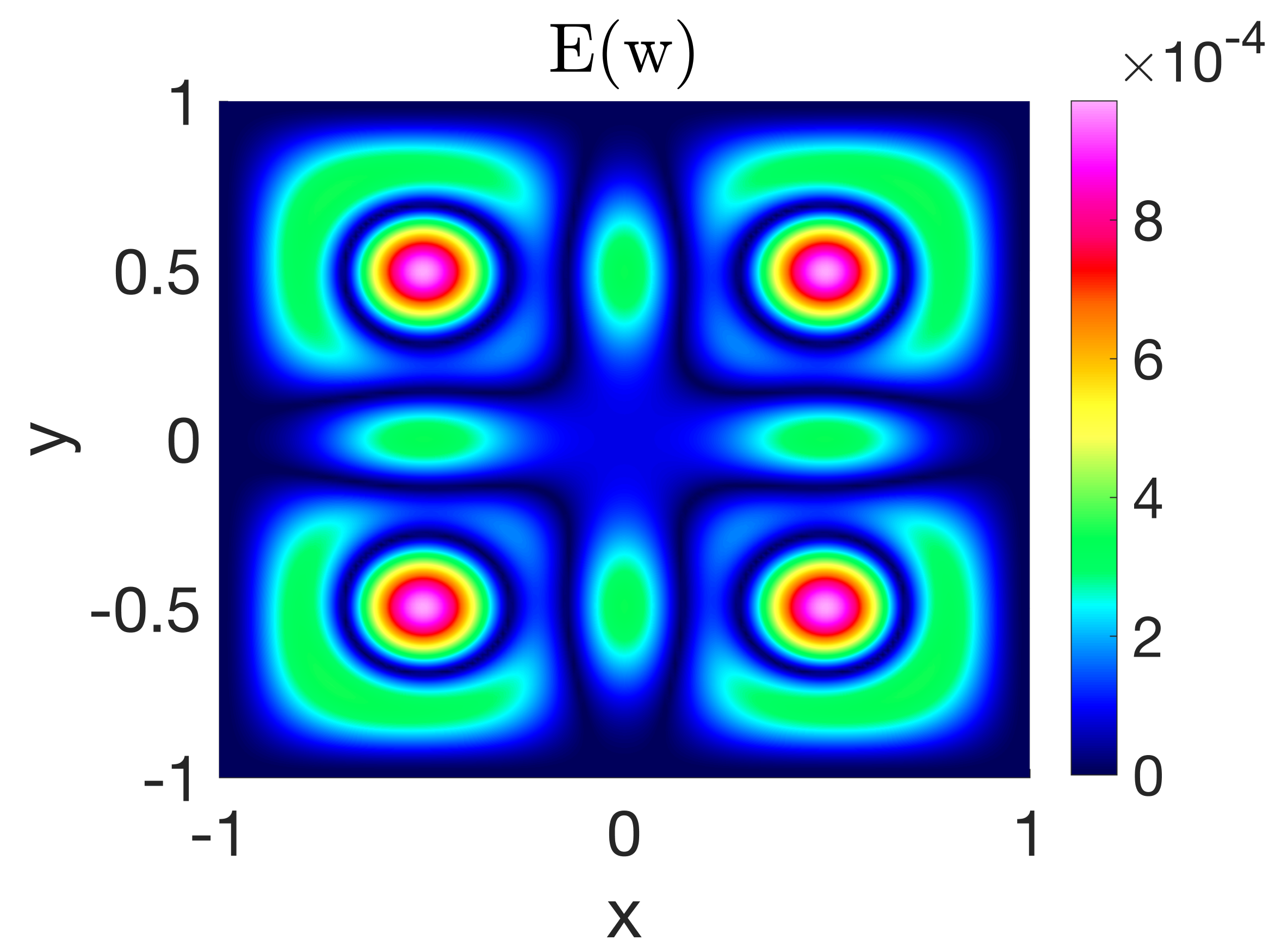}
  \end{subfigure} 
  \begin{subfigure}[b]{0.32\linewidth}
    \centering
    Simply Supported\vspace{0.1in}
    \includegraphics[width=1\linewidth]{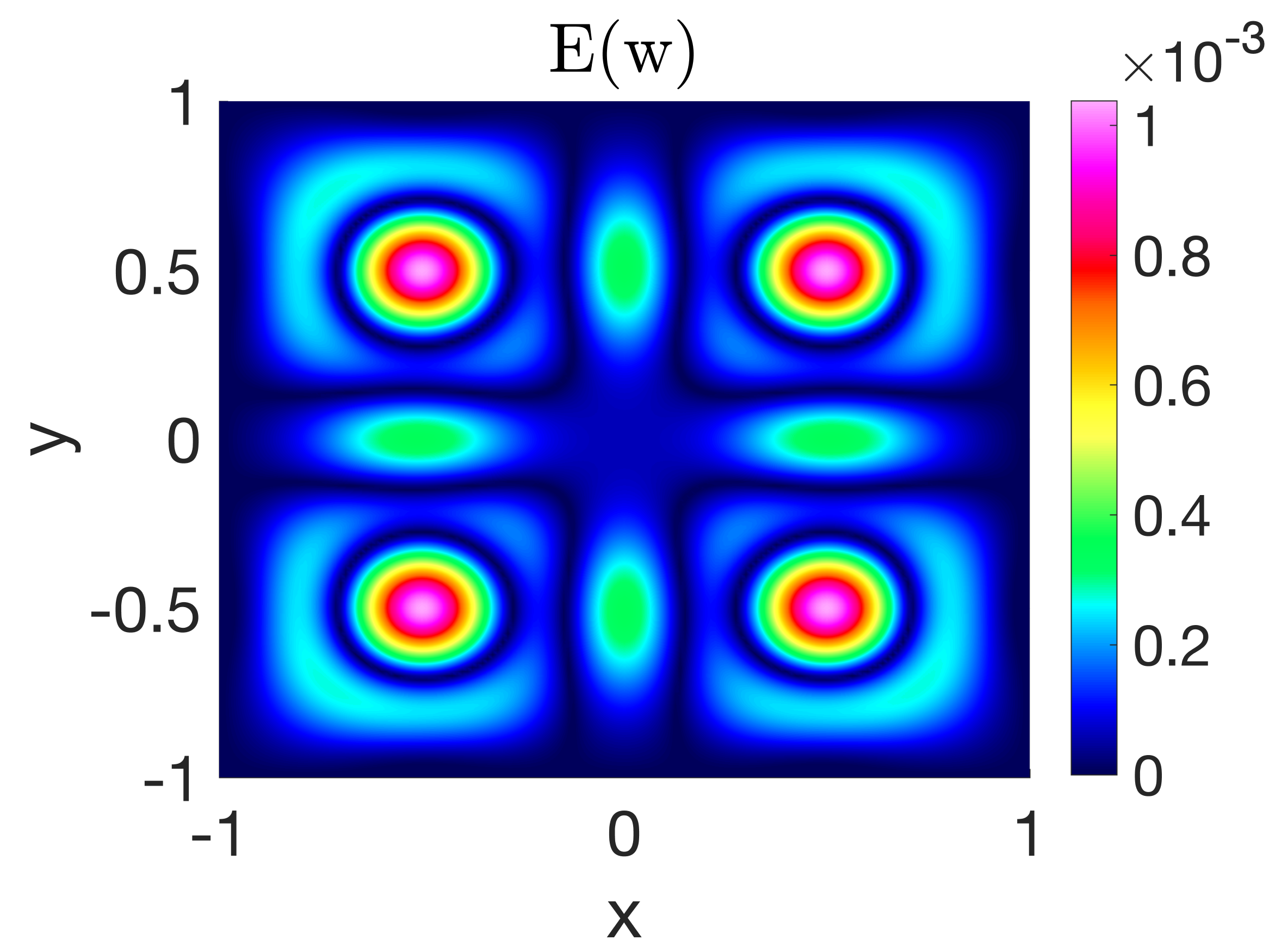}
  \end{subfigure}
  \begin{subfigure}[b]{0.32\linewidth}
    \centering
    Free\vspace{0.1in}
    \includegraphics[width=1\linewidth]{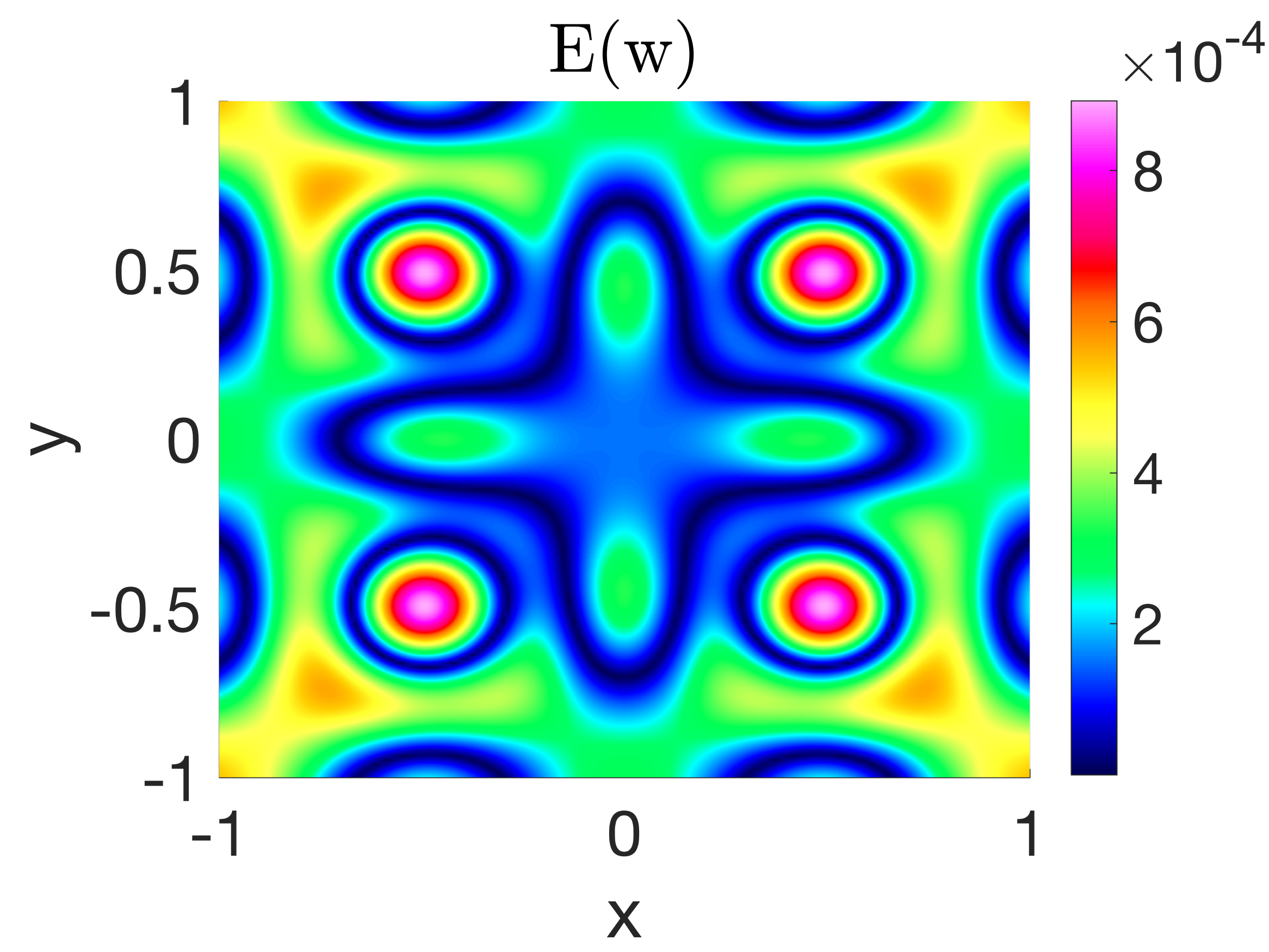}
  \end{subfigure}
      \begin{subfigure}[b]{0.32\linewidth}
        \centering
        \bigskip
    \includegraphics[width=1\linewidth]{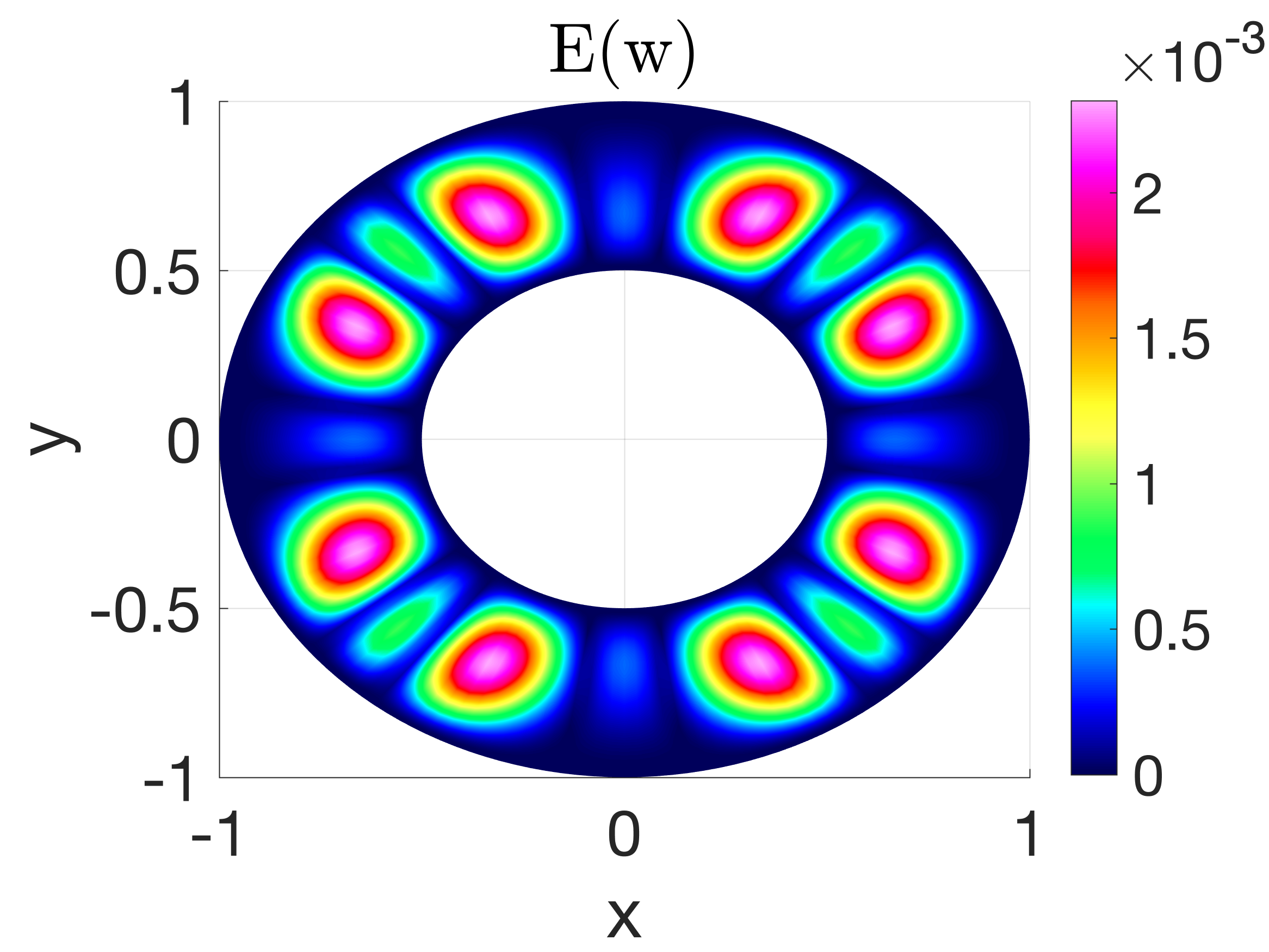}
      \end{subfigure}
  \begin{subfigure}[b]{0.32\linewidth}
    \centering
    \includegraphics[width=1\linewidth]{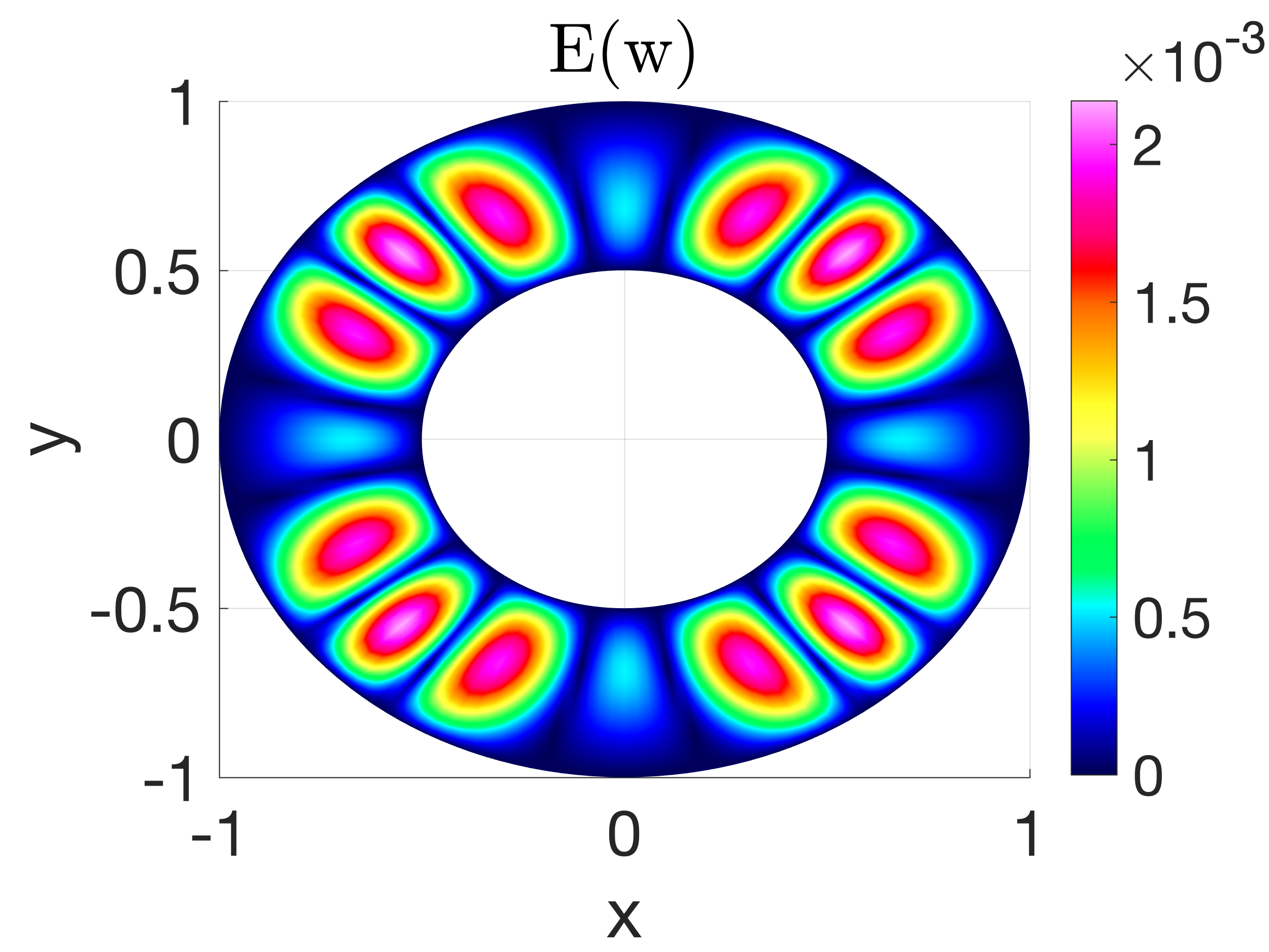}
  \end{subfigure}
  \begin{subfigure}[b]{0.32\linewidth}
   \centering
   \includegraphics[width=1\linewidth]{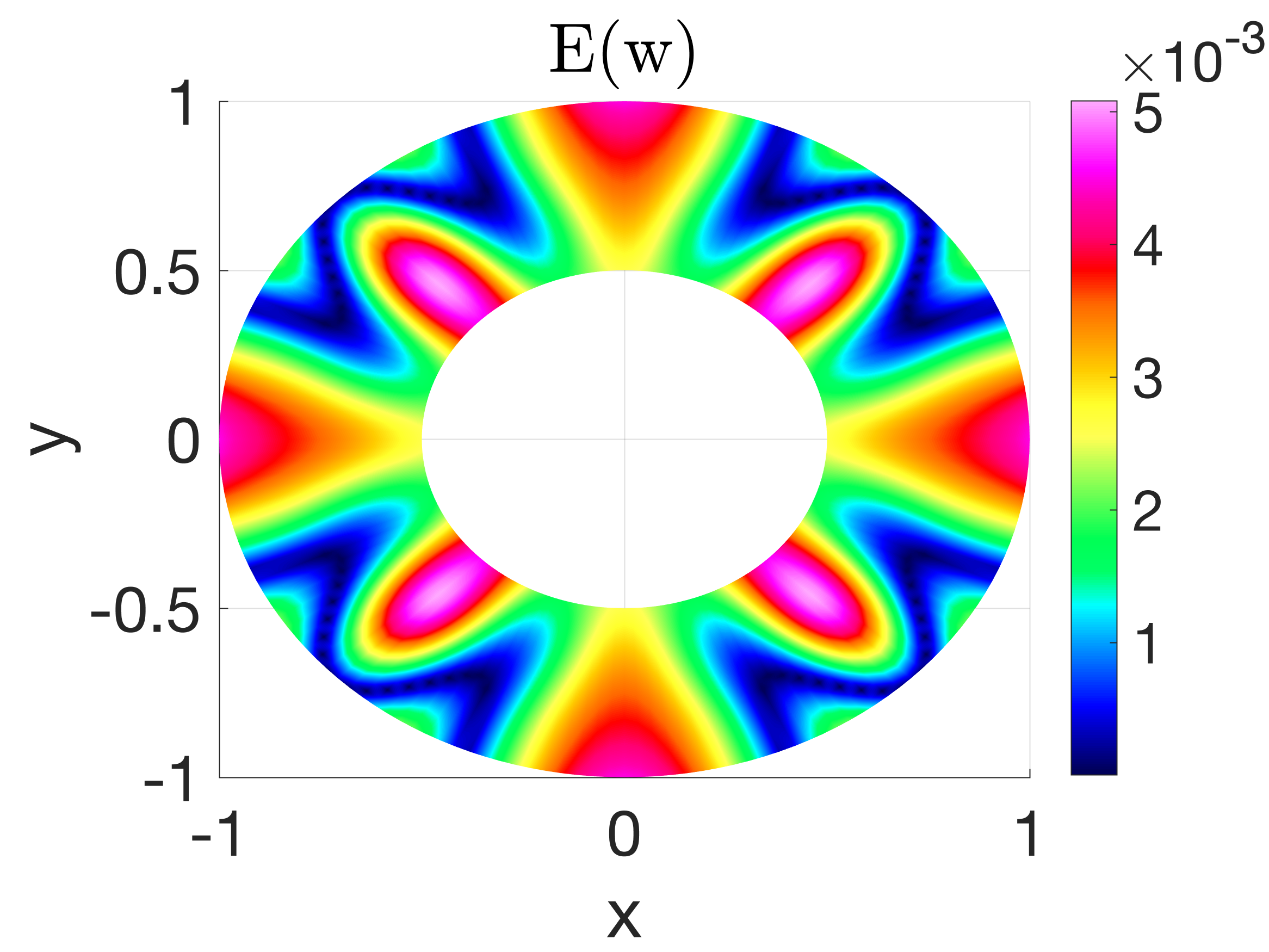}\\
     \end{subfigure}
  \caption{ Contour plots showing the errors of the numerical solutions for the displacement $w$ with various boundary conditions when $\text{t}=1$. Results shown here are obtained   on grid $\Gc_{160}$ using the PC22  method for the square plate and the NB2 method for the annular plate. }
  \label{fig:ManufactureErr}
  \end{figure}  

Given initial conditions from the exact solution at $t=0$, we solve   the test problem  to  $t=1$ on a sequence of refined grids $\Gc_N$, where $N=10\times 2^j$ represents  the number of grid points in each axial direction with $j$ ranging from $0$ through $4$.  All the  boundary conditions listed  in \eqref{eq:clampedBC}--\eqref{eq:freeBC} as well as   both of the proposed schemes (i.e., PC22 and NB2) are considered, but  we selectively present two of the  numerical solutions obtained on the finest considered grid (i.e., $\Gc_{160}$) in Figure~\ref{fig:ManufactureSol}.  In particular, the solution presented for the square plate is  solved using the PC22 scheme subject to the free boundary conditions, while the plot for the annular plate  is generated from the  NB2 scheme with simply supported boundary conditions. We note that 
numerical solutions  for the other cases are  similar, since the test problem is designed to have the same exact solution \eqref{eq:manufacturedExact}.

Let $E(w)=|w_{\iv}(t)-w_e(\xv_\iv,t)|$ denote the error function of  a numerical solution $w_{\iv}(t)$, and we show in Figure~\ref{fig:ManufactureErr}  the  contour plots of $E(w)$  on $\Gc_{160}$   to demonstrate the accuracy of our schemes for all the boundary conditions.  Results shown for the square plate  are obtained using the PC22  scheme, and  those for the annular plate are  solved with the NB2 scheme.  We observe that the numerical solutions subject to all the boundary conditions are accurate in the sense that the errors are small and smooth throughout the domain including the boundaries. Errors for all the other cases behave similarly, so their plots are omitted here to save space.

{
\newcommand{\figWidth}{6cm}
\def\xa{13.}
\def\ya{10.5}
\newcommand{\trimfig}[2]{\trimw{#1}{#2}{0.08}{0.05}{0.08}{0.0}}
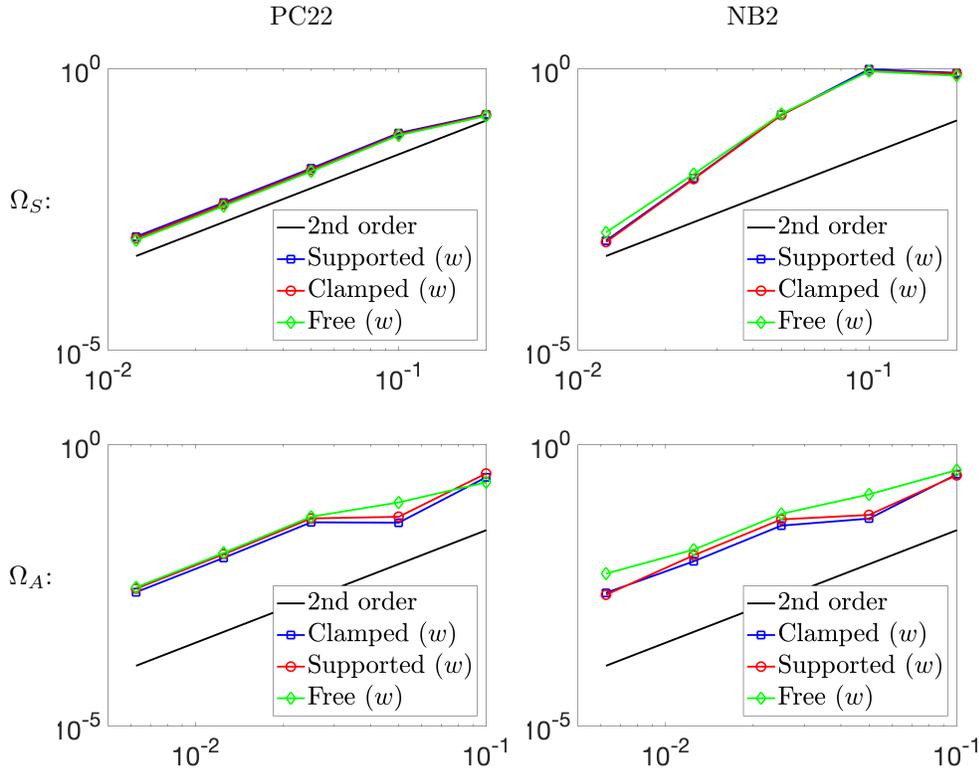
\begin{figure}[h]
\begin{center}
\begin{tikzpicture}[scale=1]
  \useasboundingbox (0.0,0.0) rectangle (\xa,\ya);  

\draw(0.25,5.) node[anchor=south west,xshift=0pt,yshift=0pt] {\trimfig{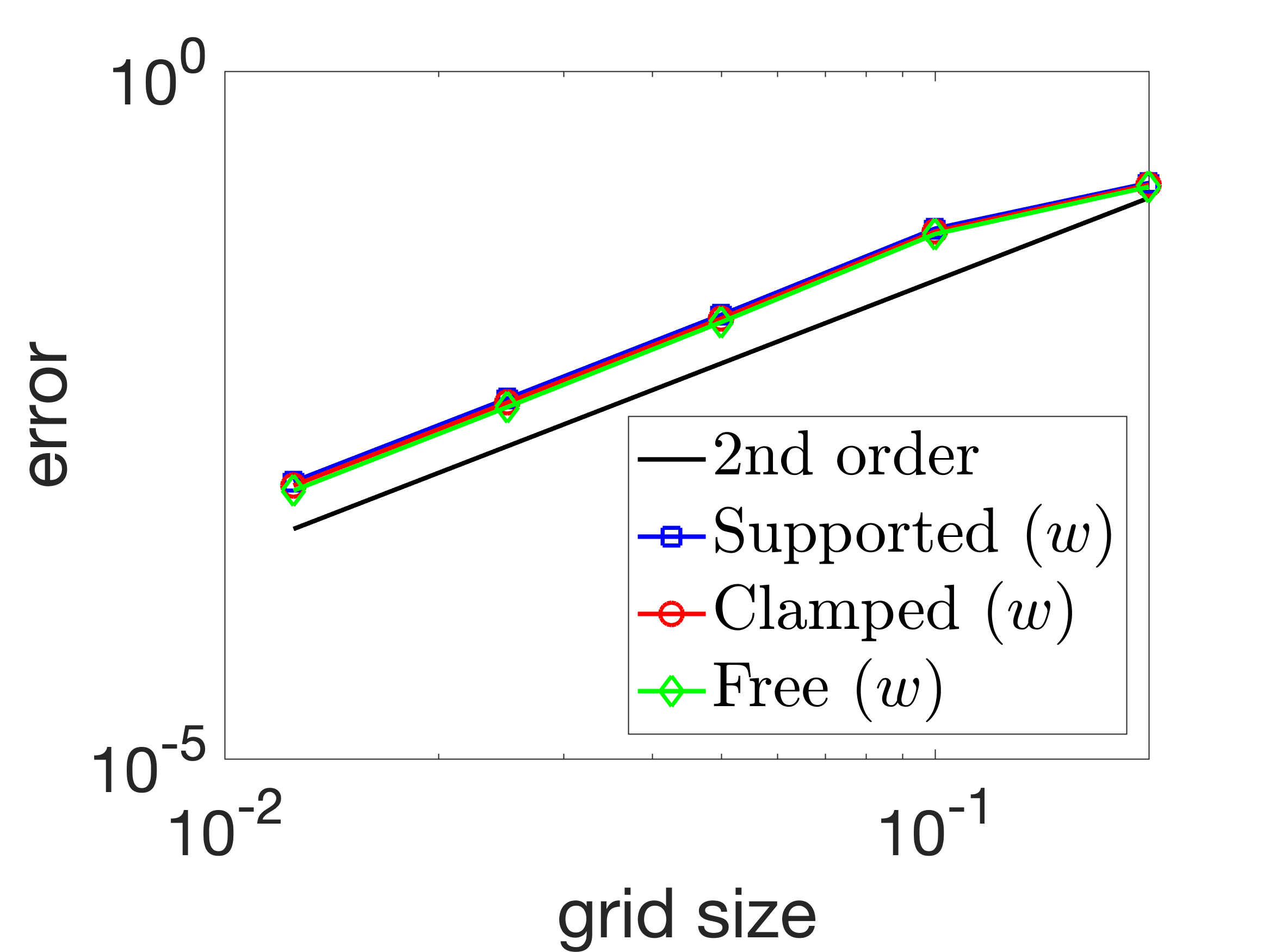}{\figWidth}};
\draw(6.5,5.) node[anchor=south west,xshift=0pt,yshift=0pt] {\trimfig{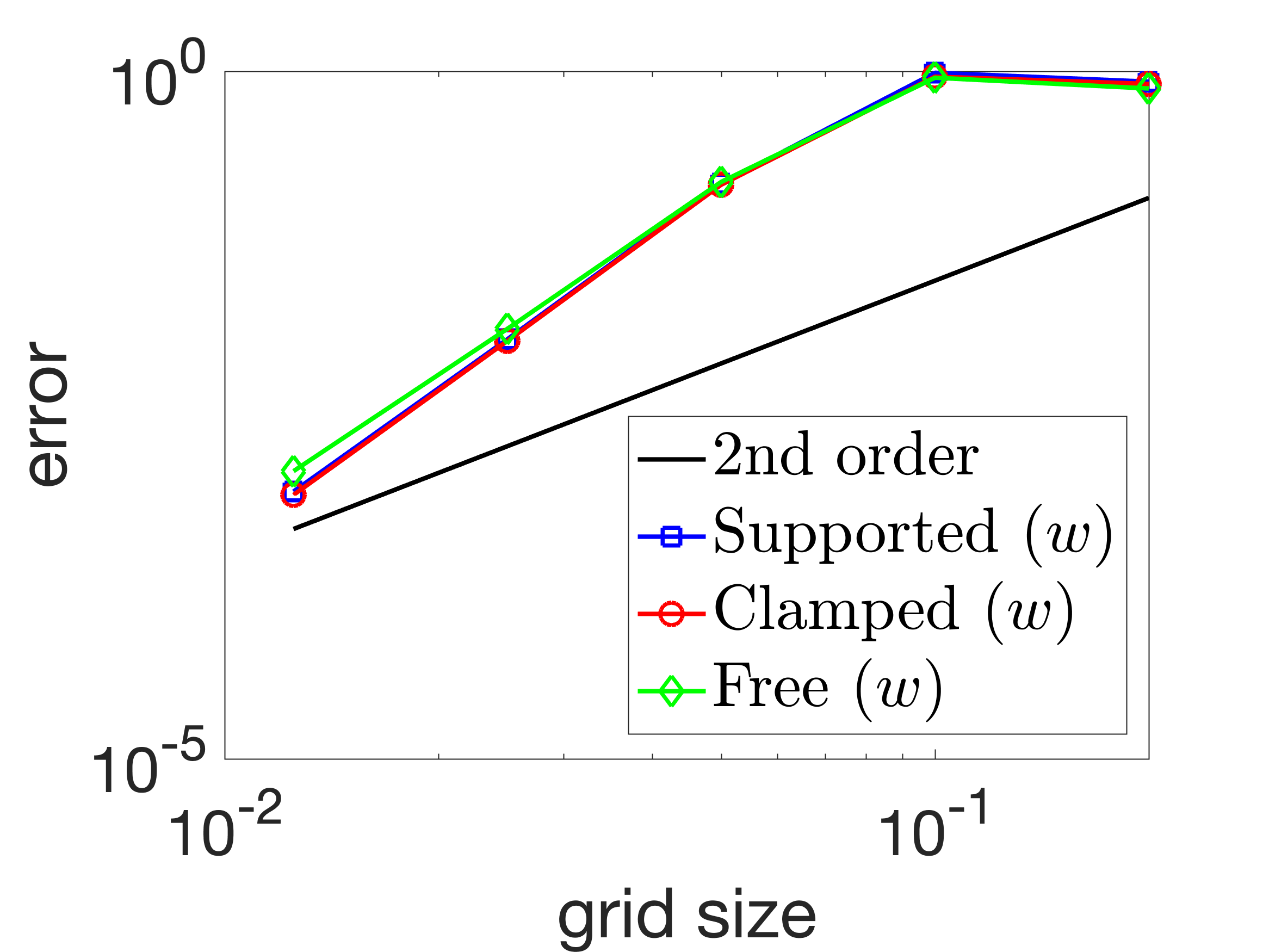}{\figWidth}};

\draw(0.25,0.0) node[anchor=south west,xshift=0pt,yshift=0pt] {\trimfig{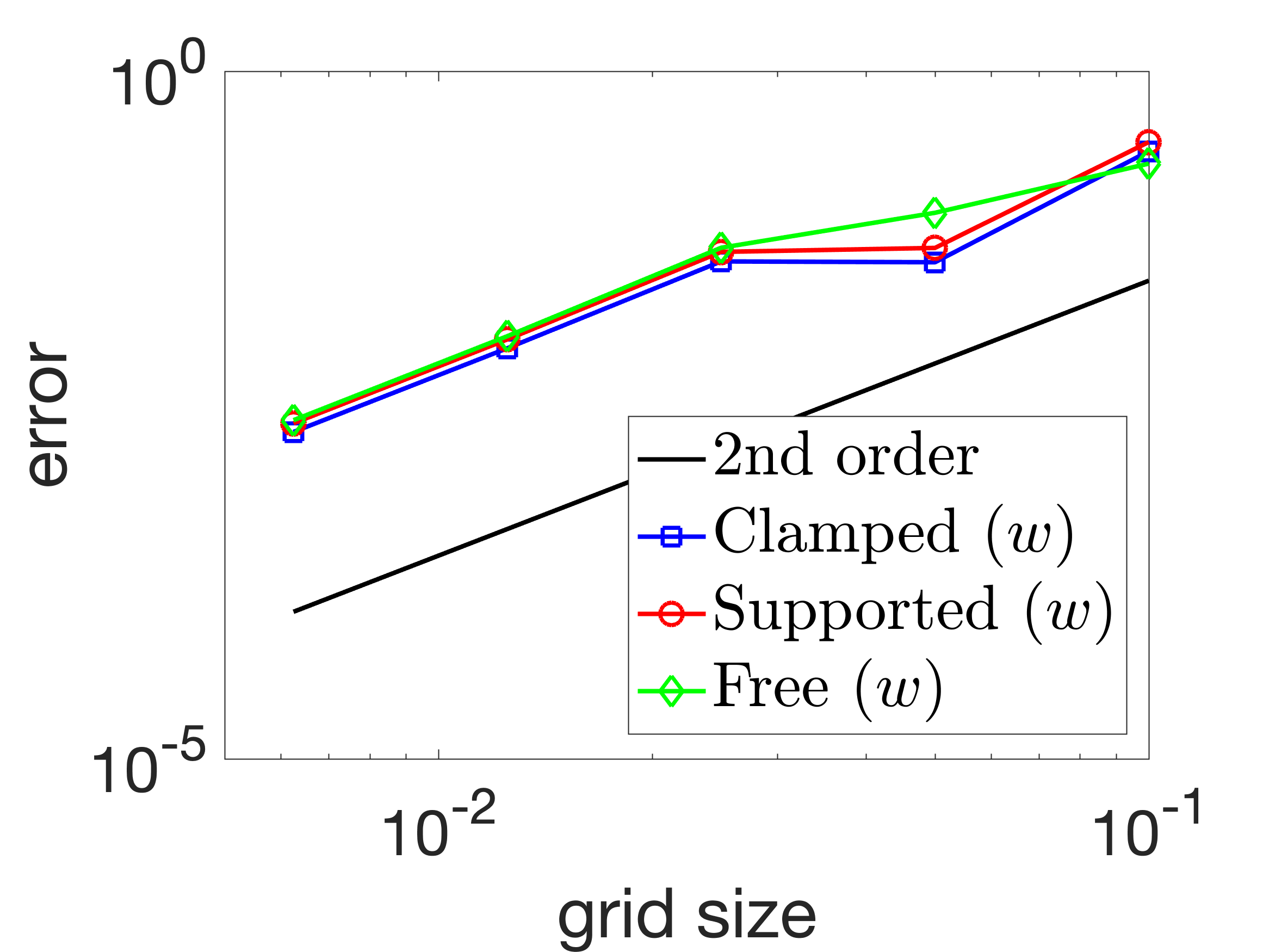}{\figWidth}};
\draw(6.5,0.0) node[anchor=south west,xshift=0pt,yshift=0pt] {\trimfig{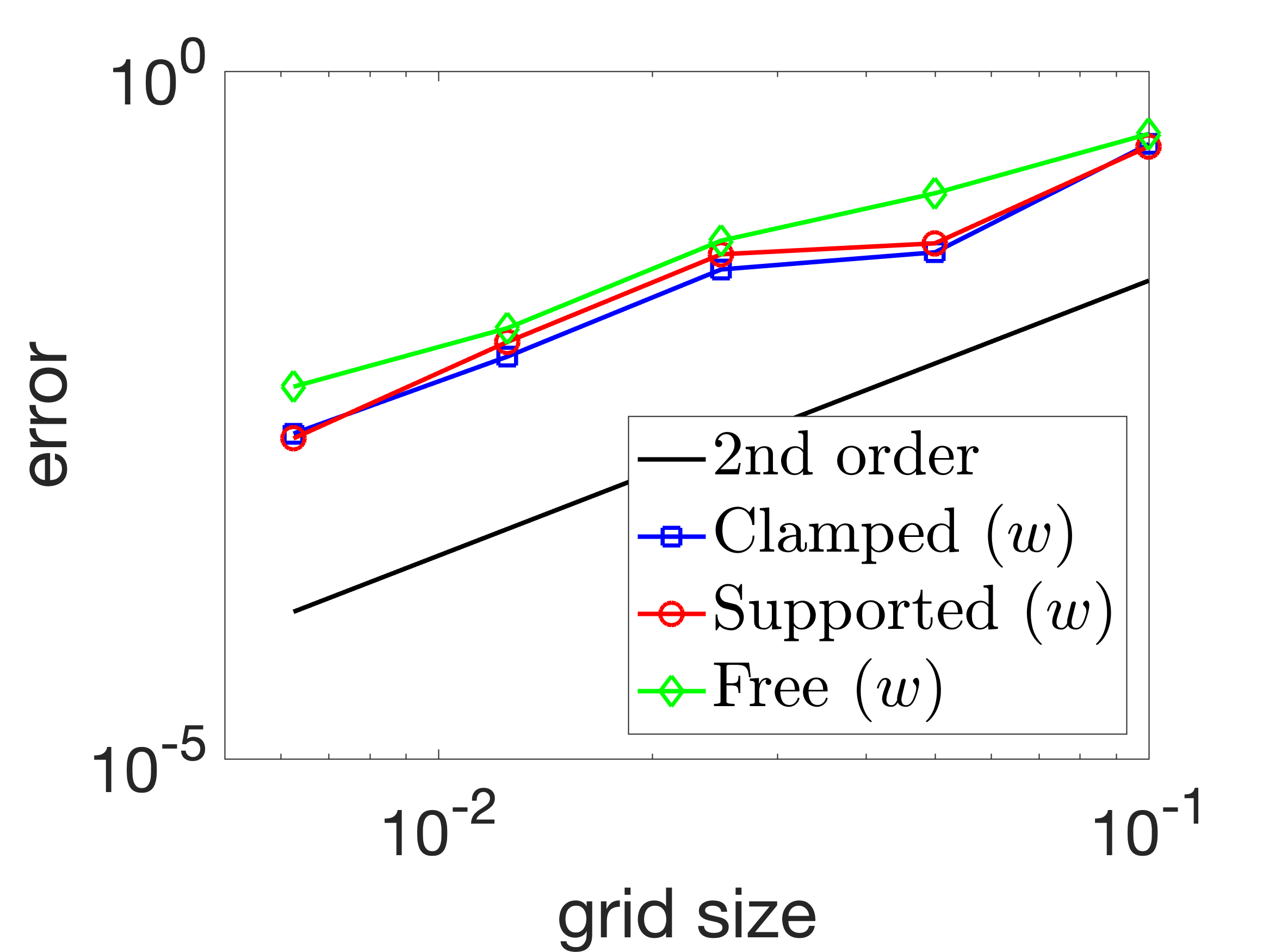}{\figWidth}};

\draw(0,2.5)  node[anchor=south ,xshift=0.4cm,yshift=0pt] {$\Omega_A$:};
\draw(0,7.5)  node[anchor=south,xshift=0.4cm,yshift=0pt] {$\Omega_S$:};
\draw(4,10)  node[anchor=south,xshift=0pt,yshift=0pt] {PC22};
\draw(10.,10)  node[anchor=south,xshift=0pt,yshift=0pt] {NB2};

%
\end{tikzpicture}

\end{center}
\caption{Convergence rates of  $w$  for both the square  ($\Omega_S$) and the annular ($\Omega_A$) plates subject to all the boundary conditions  are  presented.    The upper panel of plots  show  the results of the square  plate, and the lower panel illustrates  the results of the annular plate.  The  errors for both numerical methods (PC22 and NB2)  are computed at $t=1$ and measured in the maximum norm.}\label{fig:ManufactureSquareConv}
\end{figure}
}

  Convergence studies for both the square and annular plates subject to all the boundary conditions  are  performed for   the numerical results obtained on the sequence of  grids $\Gc_N$'s using both numerical methods. We plot the maximum norm errors $||E(w)||_\infty$ against the grid size  together with a second-order reference curve in log-log scale to reveal the order of accuracy. In the top row of Figure~\ref{fig:ManufactureSquareConv}, we show the results for the square plate; and in the  bottom row of the same figure, we show the results for the annular plate.  For all the boundary conditions and both of the numerical schemes, we  observe   the expected second-order  accuracy regardless of the plate shape.

  \subsection{Vibration of plates}\label{sec:vibrationOfPlates}
  Mechanical vibrations are problems of great interest in  engineering and material sciences. For  a thin plate-like  structure, a 2D plate theory is capable of   giving an excellent approximation to the actual 3D motion. The vibration of a plate  can be caused    either by  displacing the plate from its stress-free state or by exerting an  external  forcing, where  the former is  referred to as   free vibration and the latter is called forced vibration.
  Among the numerous plate models,  the Kirchhoff-Love theory  is  most commonly used. To further validate  the numerical properties of  the proposed schemes, we consider  the vibration problems of the generalized  Kirchhoff-Love plate \eqref{eq:generalizedKLPlate}, which  include  the study of    natural frequencies and mode shapes of vibration, and   propagating or  standing waves in the plate.

 \subsubsection{Vibration with known analytical  solutions}
 The classical Kirchhoff-Love plate with some simple specifications  can be solved analytically. 
 We consider a thin plate on a rectangular domain (i.e., $\Omega=[0,L]\times [0,H]$).  Analytical solutions  to the classical Kirchhoff-Love plate  equation subject to  simply supported boundary conditions \eqref{eq:supportedBC} are available  for both the free and forced vibration cases. We solve   each case numerically and compare  our approximations  with the analytical solutions to reveal the stability and accuracy of our schemes.  
 
 {\bf $\bullet$ Free vibration.} Consider the  free vibration case; i.e., the forcing function in \eqref{eq:generalizedKLPlate} is zero ($F(\xv,t)\equiv0$). In this case, the governing equation   can be  analytically solved using separation of variables or Fourier transformation.    Let  $A_{mn}$ and $B_{mn}$  denote  the coefficients to be determined by  the initial conditions and the orthogonality of Fourier components, then the general solution to this simple plate  can be expressed as the following infinite series,
 \begin{equation}\label{eq:generalSolution}
   w(\xv,t)=\sum_{m=1}^{\infty}\sum_{n=1}^{\infty}\sin\frac{m\pi x}{L}\sin\frac{n\pi y}{H}(A_{mn}\cos\omega_{mn}t+B_{mn}\sin\omega_{mn}t),
 \end{equation}
 where the natural frequencies of vibration for this plate is found to be
 \begin{equation}\label{eq:NaturalFrequencySupported}
\omega_{mn}={\pi^2}\left(\frac{m^2}{L^2}+\frac{n^2}{H^2}\right)\sqrt{\frac{D}{ \rho h}}.
\end{equation}

 Standing wave test problems can be constructed by specifying modes of vibration as the given functions for the initial conditions \eqref{eq:IC}; namely, 
 \begin{equation}\label{eq:freeVibrationIC}
 w_0(\xv)=\sin\frac{m\pi x}{L}\sin\frac{n\pi y}{H},~\text{and}~v_0(\xv)=0.
 \end{equation}
 Enforcing the above  initial conditions on the general solution \eqref{eq:generalSolution}, we deduce the exact standing wave solution for each 2-tuple $(m,n)$, 
 \begin{equation}\label{eq:freeVibrationExact}
 w_e(\xv,t) =\sin\frac{m\pi x}{L}\sin\frac{n\pi y}{H}\cos\omega_{mn}t.
 \end{equation}

\begin{figure}[h] 
  \begin{subfigure}[b]{0.32\linewidth}
    \centering
    \includegraphics[width=1\linewidth]{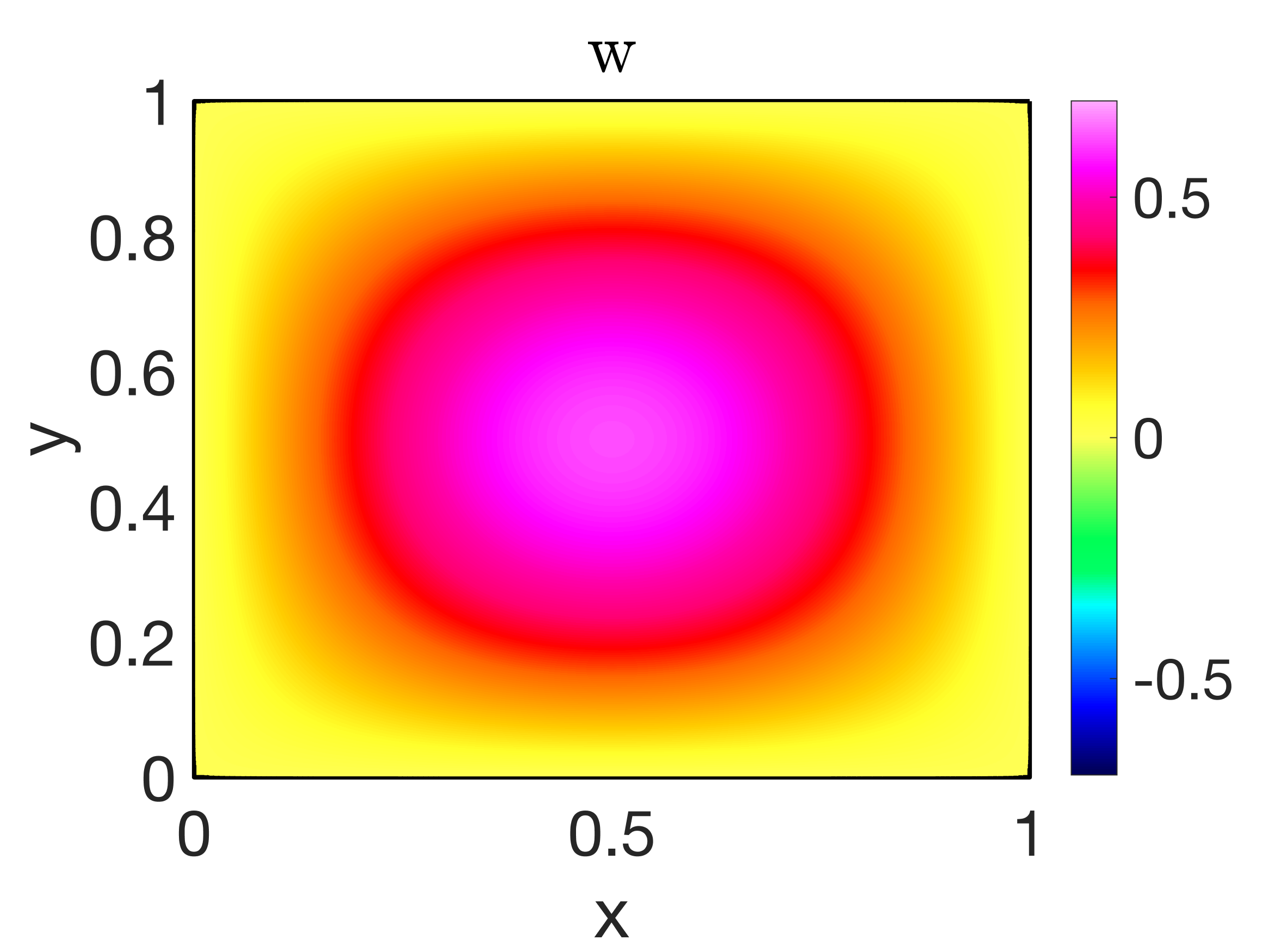}
        \caption{\footnotesize $m=1,n=1$}
  \end{subfigure}
  \begin{subfigure}[b]{0.32\linewidth}
    \centering
    \includegraphics[width=1\linewidth]{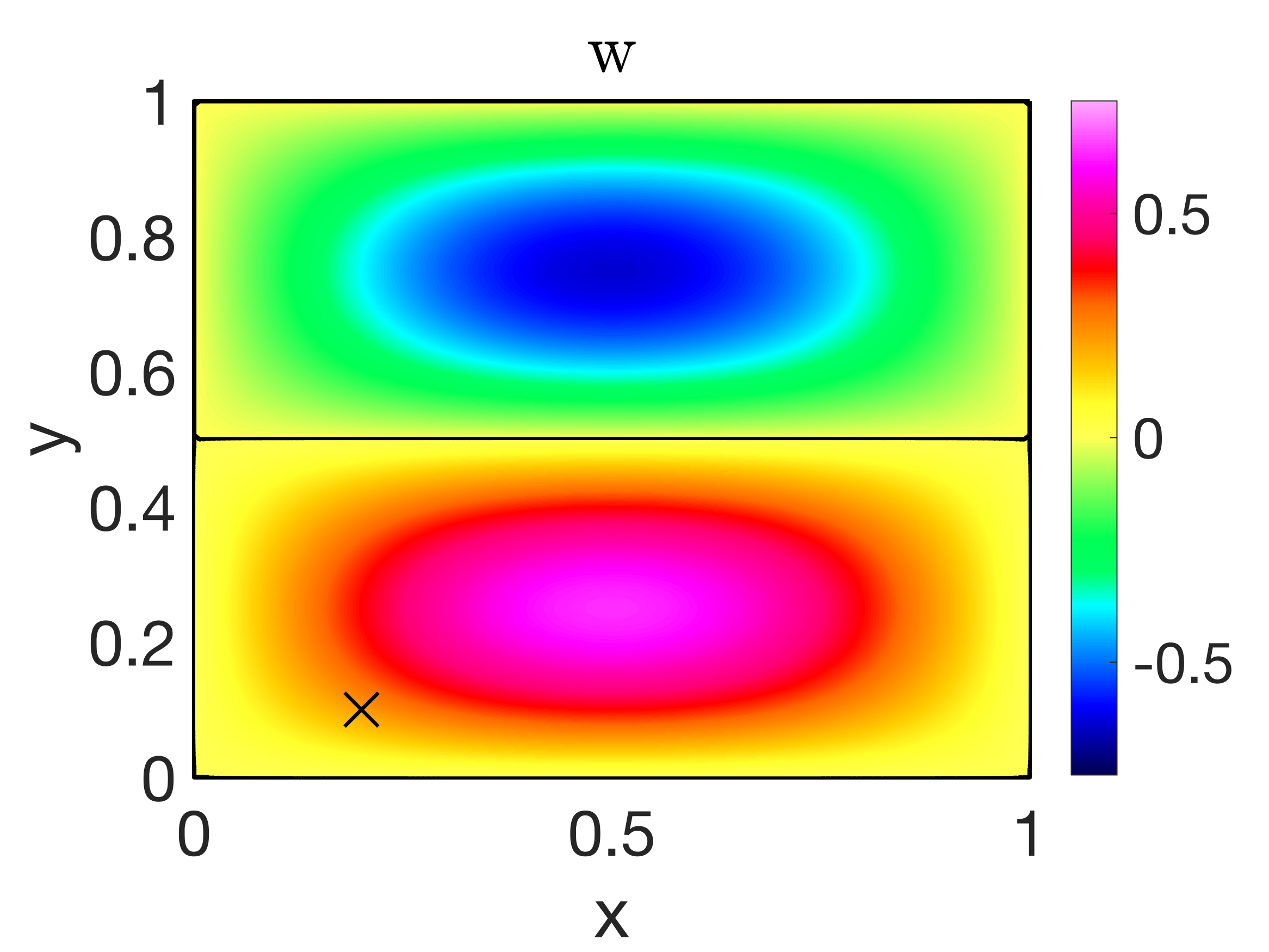}
            \caption{\footnotesize $m=1,n=2$}\label{fig:probedCase}
  \end{subfigure} 
    \begin{subfigure}[b]{0.32\linewidth}
    \centering
    \includegraphics[width=1\linewidth]{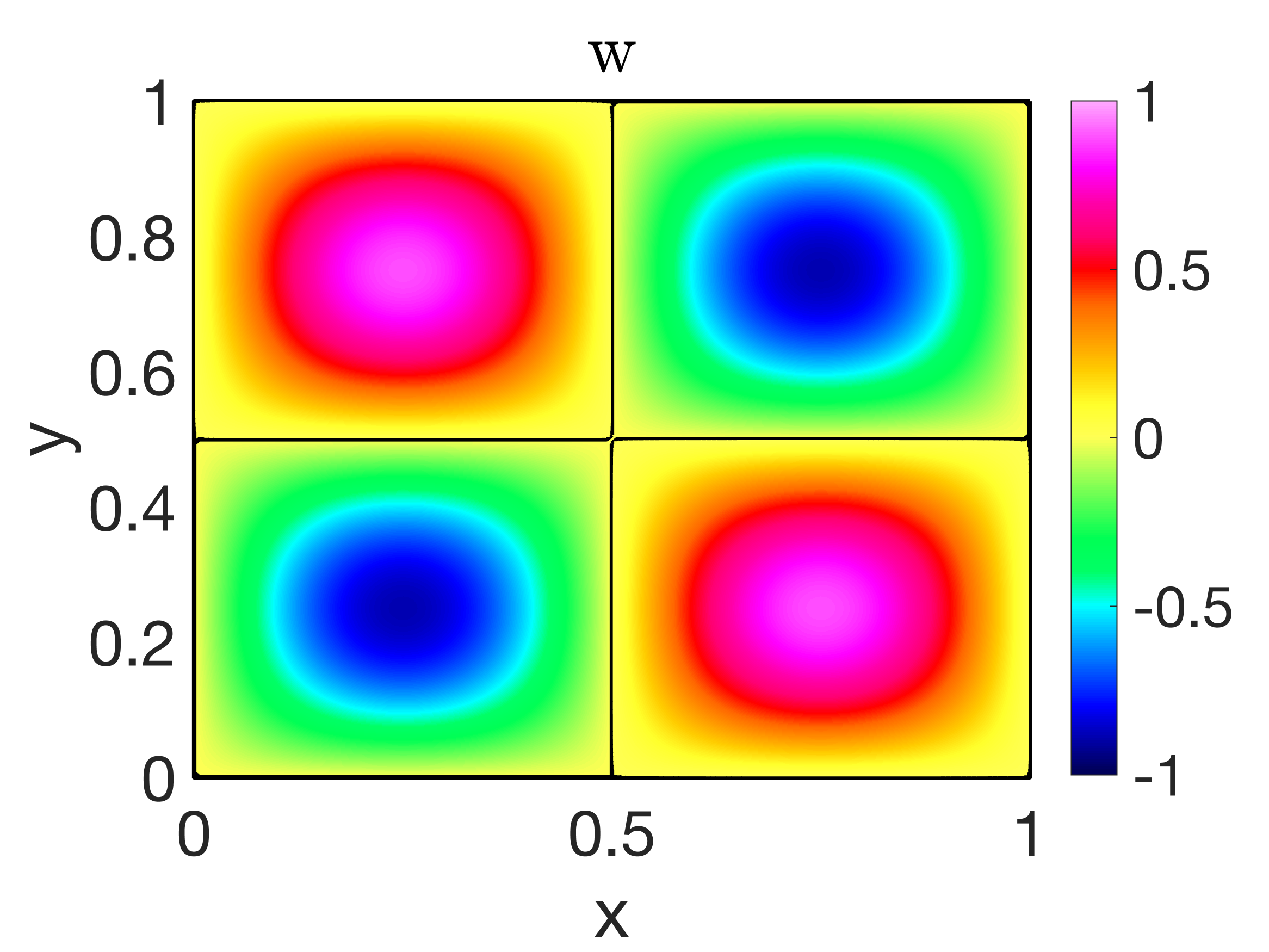}
            \caption{\footnotesize $m=2,n=2$ }
  \end{subfigure}
  \begin{subfigure}[b]{0.32\linewidth}
    \centering
    \includegraphics[width=1\linewidth]{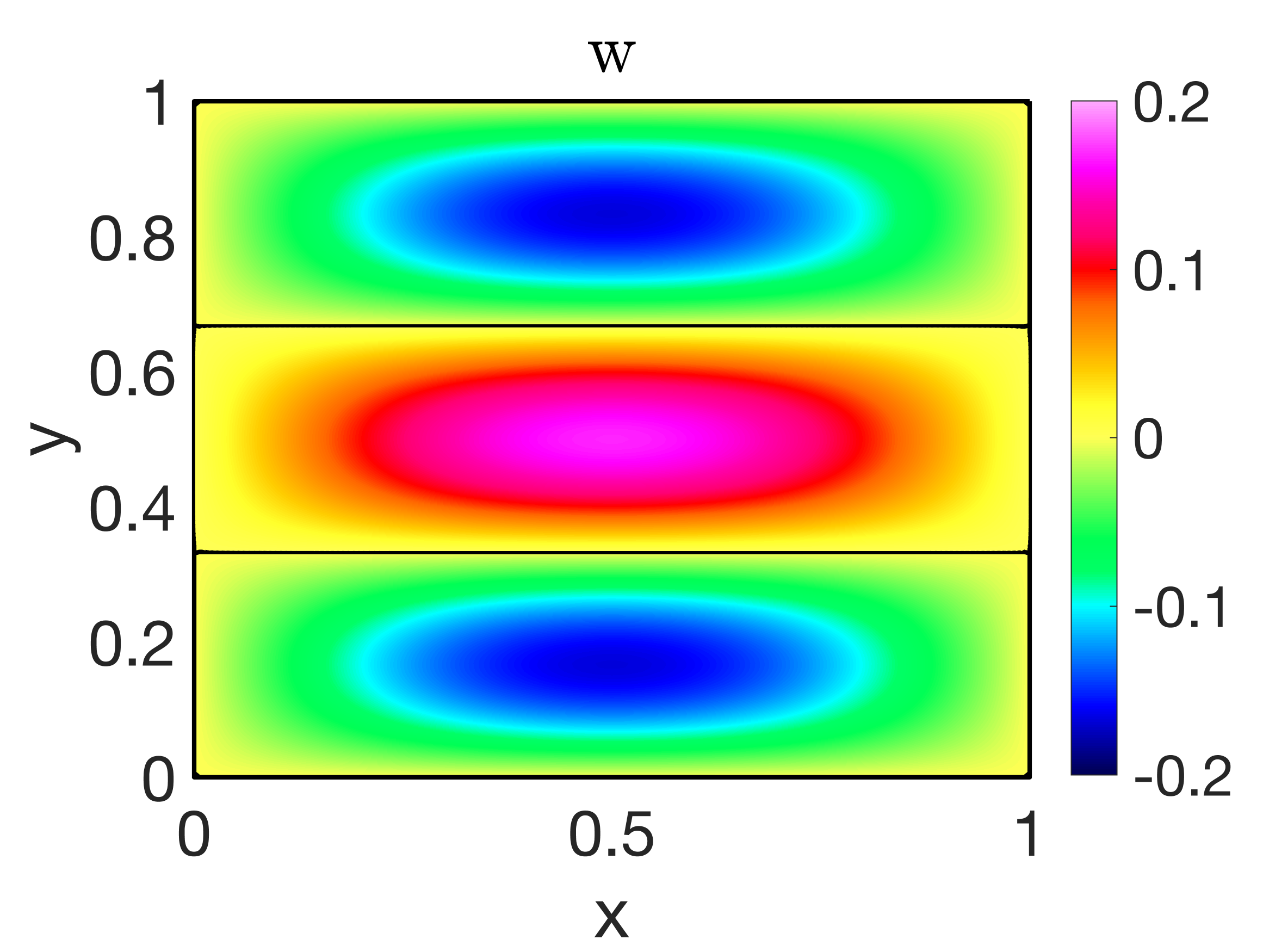}
            \caption{\footnotesize $m=1,n=3$ }
  \end{subfigure} 
    \begin{subfigure}[b]{0.32\linewidth}
    \centering
    \includegraphics[width=1\linewidth]{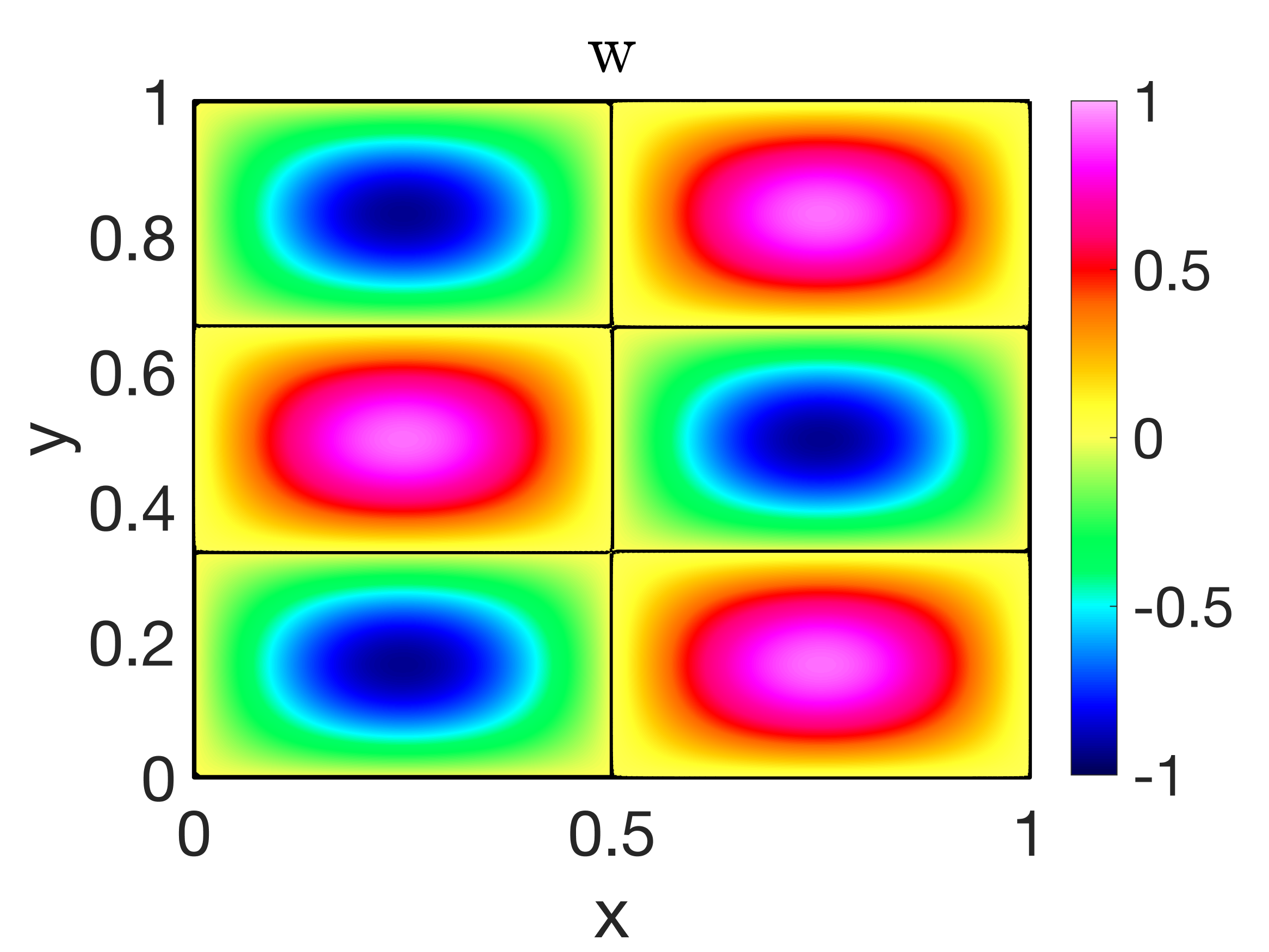}
            \caption{\footnotesize $m=2,n=3$ }
  \end{subfigure} 
    \begin{subfigure}[b]{0.32\linewidth}
    \centering
    \includegraphics[width=1\linewidth]{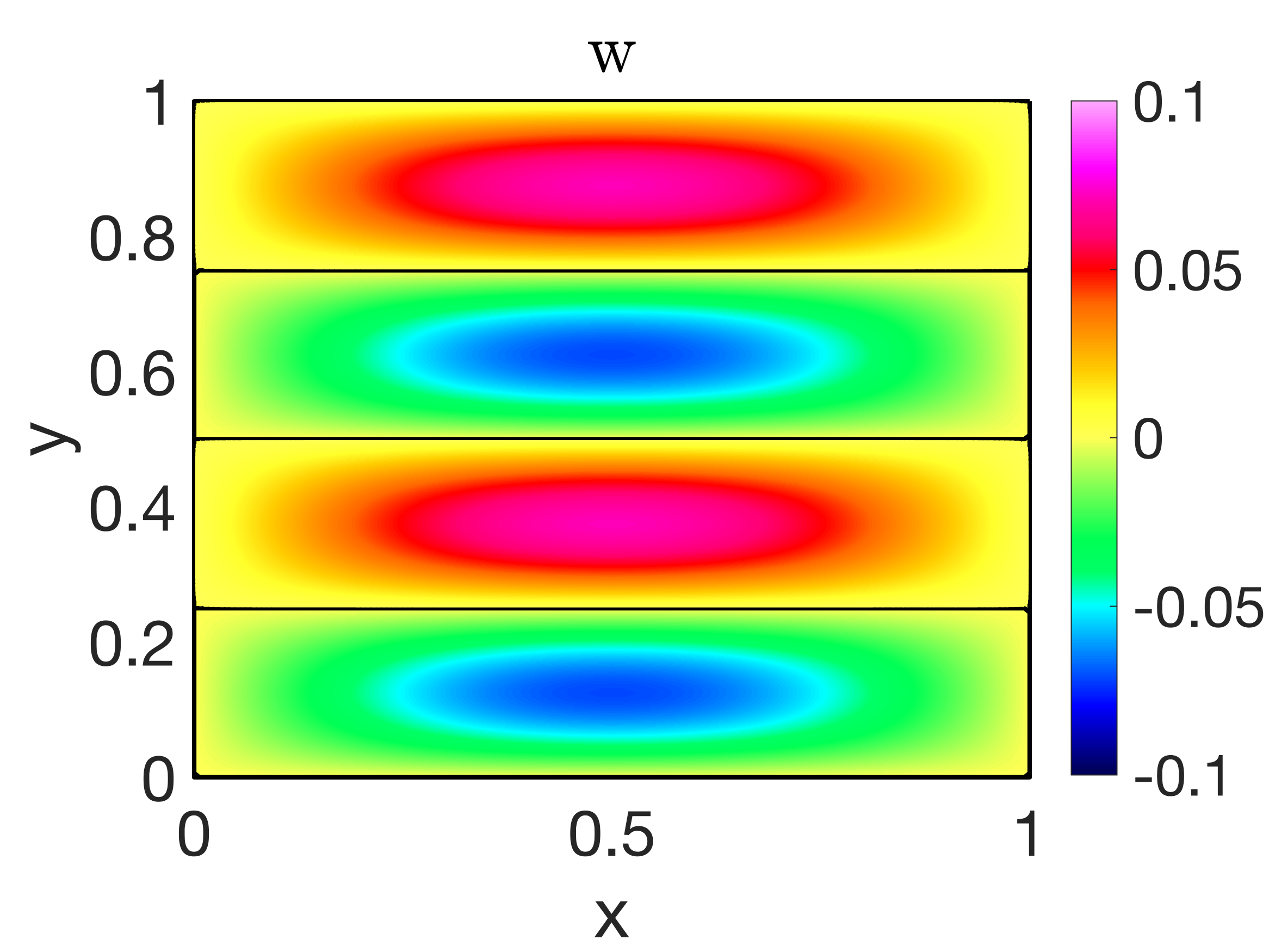}
            \caption{\footnotesize $m=1,n=4$}
      \end{subfigure} 
  \caption{Standing wave solutions  with the nodal lines  at time $t=1$  for some $(m, n)$ values. Simulations are performed using the PC22 scheme on grid $\Gc_{160}$. Results obtained using NB2 scheme are similar.} \label{fig:ModeShapesSupportedFreeVibABAM}
\end{figure}

 We solve the standing wave test problem for a  few modes to validate our numerical schemes. That is to say,  we specify the initial conditions \eqref{eq:freeVibrationIC} with  various values of $(m,n)$,   and compare the numerical approximations with the exact solution \eqref{eq:freeVibrationExact}.
 For simplicity, the rectangular domain is restricted to  a unit square in the computations; namely, we set  $L=H=1$ in  $\Omega=[0,L]\times [0,H]$. The parameters for this test are specified as   $\rho h=2.7$, $K_0=0$, $T=0$, $D=6.4527$, $K_1=0$, $T_1=0$ and $\nu=0.33$.  Note that this is a classical Kirchhoff-Love  model because   only the bending dynamics is accounted for (namely, $D\neq0$).
 Both of the proposed schemes  are  used to conduct numerical simulations, but only the results of the PC22 scheme are presented in Figure~\ref{fig:ModeShapesSupportedFreeVibABAM} since  the other scheme produces  comparable results. In Figure~\ref{fig:ModeShapesSupportedFreeVibABAM}, we show the contour plots of the displacement $w$ at $t=1$ for  a few $(m,n)$-tuples. In the plots, we also show  the zero contours, which represent the nodal lines of the standing wave solutions. It is clear that the patterns of these nodal lines  exhibited in  the numerical solutions resemble those of the corresponding modes of vibration  used in the  initial conditions \eqref{eq:freeVibrationIC}.
 
\begin{figure}[h] 
  \begin{subfigure}[b]{0.32\linewidth}
    \centering
    \includegraphics[width=1\linewidth]{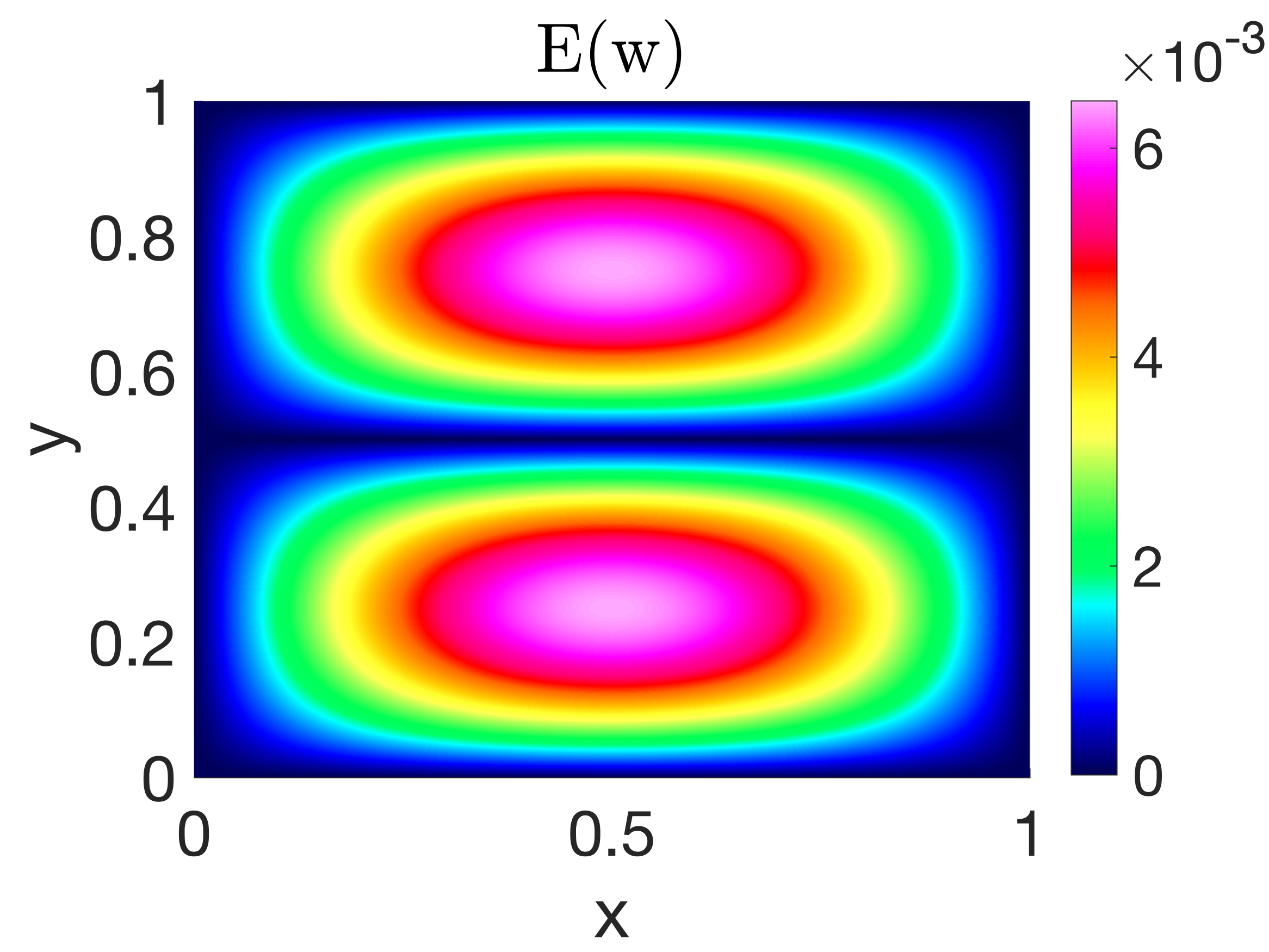}
    \caption*{\footnotesize Error}
  \end{subfigure}
    \bigskip
  \begin{subfigure}[b]{0.32\linewidth}
    \centering
    \includegraphics[width=1\linewidth]{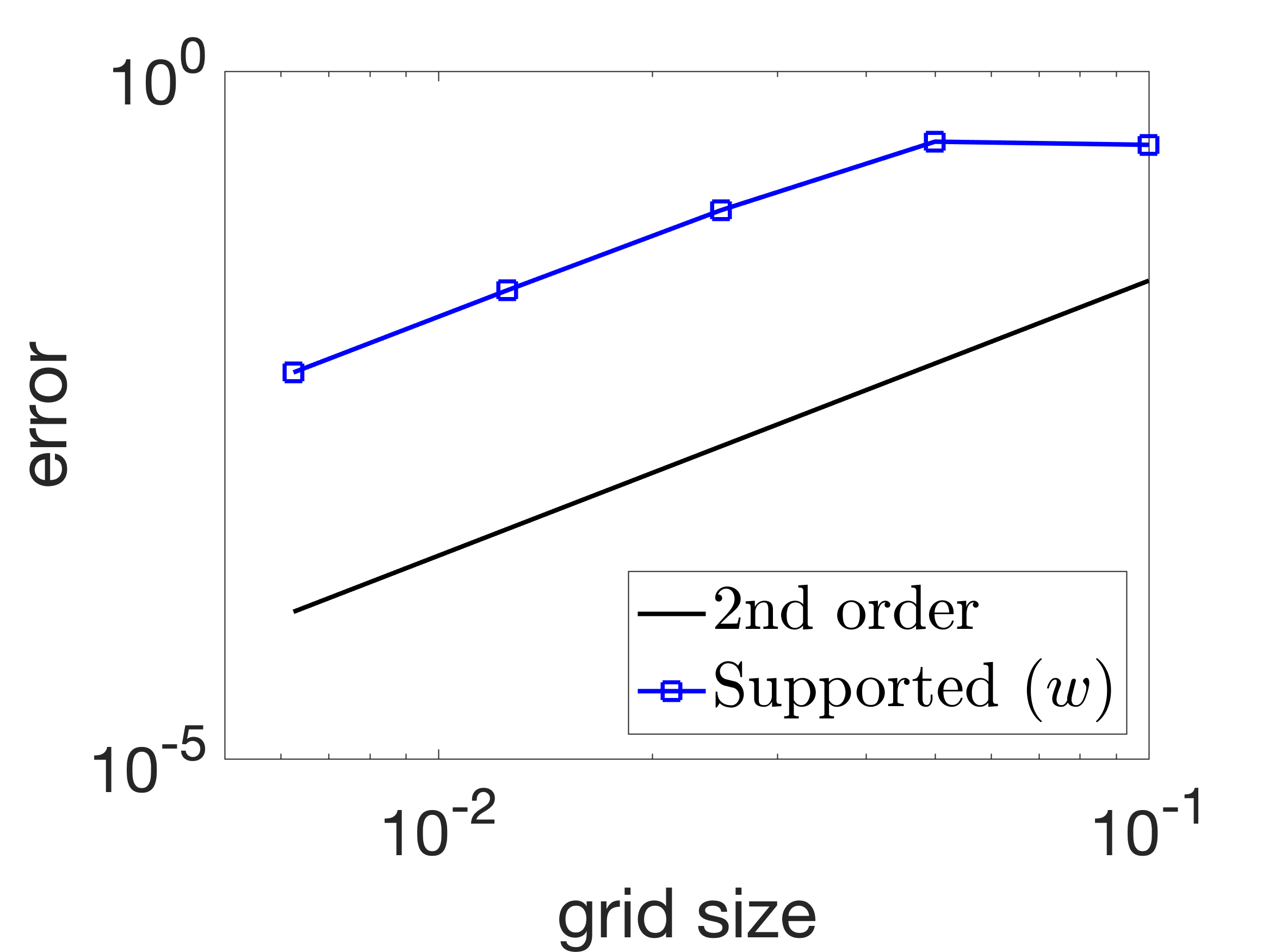}
    \caption*{\footnotesize Convergence rate} 
  \end{subfigure} 
    \bigskip
  \begin{subfigure}[b]{0.32\linewidth}
    \centering
    \includegraphics[width=1\linewidth]{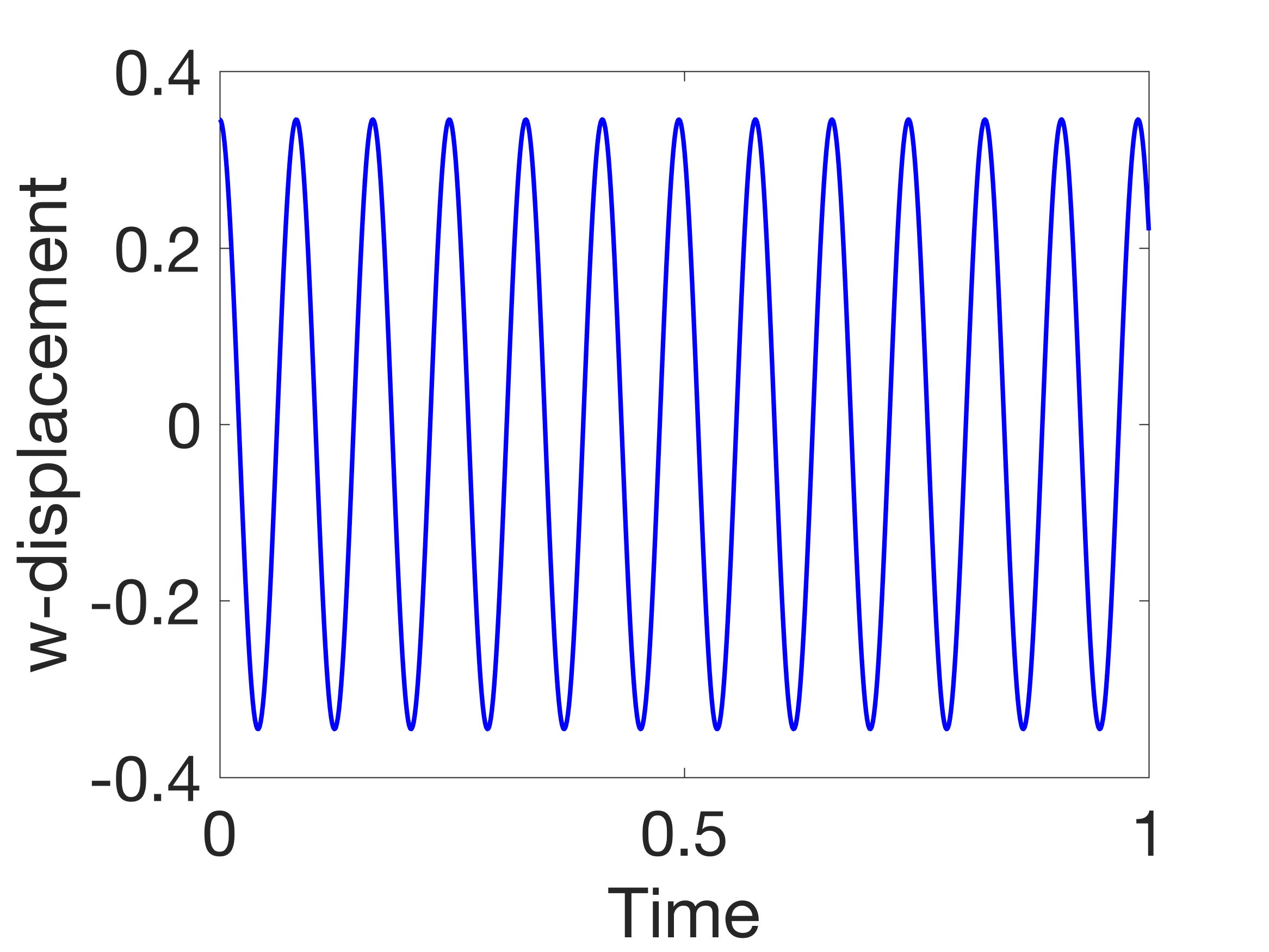}
    \caption*{\footnotesize Displacement at $\xv_p$} 
  \end{subfigure}
  \vspace{-0.4 in}
  \caption{More results for the  case   with  $m=1, n=2$. Left: the error plot of the displacement $w$ at $t=1$. Middle:  the convergence rate of $w$. Right: displacement at point $\xv_p$; the probed location is marked by a cross symbol   in  Figure~\ref{fig:probedCase}.}
  \label{fig:ErrorConvTrackPointSupportedABAMFreeVib}
\end{figure}

Given that  the exact solutions for these standing waves   are available in \eqref{eq:freeVibrationExact},  it is possible for us to   perform  mesh refinement studies to show the accuracy and the convergence of the numerical solutions. In Figure~\ref{fig:ErrorConvTrackPointSupportedABAMFreeVib}, we present  more results for the case with $m=1$ and $n=2$, which include the error  of $w$  at $t=1$ plotted in the left image, the convergence rate of $w$ shown in the middle image, and the evolution of $w$ at a probed location depicted in the right image.  The probed location is $\xv_p=(0.2,0.1)$, which is  marked by a cross symbol   in  Figure~\ref{fig:probedCase}. The error and convergence rate plots confirm that the scheme is accurate and the rate of  convergence is second order.

The  evolution of $w$ at  a point, as   the one shown in the right image of Figure~\ref{fig:ErrorConvTrackPointSupportedABAMFreeVib}, demonstrates how the standing wave solutions oscillate in time. From the evolution curve,   we can estimate the frequency of oscillation  by measuring  the number of cycles  per unit time, which  should match up with the natural frequency of the plate. As another indication of  the accuracy of our proposed schemes, we compare the numerically estimated  frequencies  with the natural frequencies. Please note  that  the natural frequencies defined in  \eqref{eq:NaturalFrequencySupported} are actually angular frequencies that measure the number of oscillations in $2\pi$ units of time.  For comparison,  we use the  ordinary frequency (measured in hertz) that is given by $f_{mn}=\omega_{mn}/(2\pi)$. 
We track the evolution of the displacement   at    $\xv_p=(0.2,0.1)$ for the modes  that correspond to the  first 9 natural frequencies, and estimate the frequencies of oscillation from the evolution curves. The  frequencies  inferred  from the numerical solutions on grid $\Gc_{160}$ are summarized in Table~\ref{Tab:NaturalFrequencySSSSFreeVib}.  We can see that, for both numerical methods, the discrepancies between the estimated frequencies and $f_{mn}$ are small for all the examined cases.

\begin{table}[h]
\begin{center}
\centering
\begin{tabular}{cccccc}
\hline
\multirow{2}{*}{$(m,n)$} & \multirow{2}{*}{ ${f_{mn}}$ } & \multicolumn{2}{c}{PC22 scheme} & \multicolumn{2}{c}{NB2 scheme} \\ \cline{3-6} 
                         &                               &    frequency (est.)         &  error (\%)     &   frequency (est.)      &  error (\%)   \\ \hline
$(1,1)$                  & $4.8567$                         & $4.8565$     & $0.0037\%$   & $4.8541$  & $0.0540\%$ \\
$(1,2)$                  & $12.1417$                       & $12.1403$    & $0.0112\%$   & $12.1359$ & $0.0480\%$ \\
$(2,2)$                  & $19.4267$                      & $19.4242$    & $0.0130\%$   & $19.4203$ & $0.0331\%$ \\
$(1,3)$                  & $24.2834$                       & $24.2769$    & $0.0268\%$   & $24.2691$ & $0.0589\%$ \\
$(2,3)$                  & $31.5684$                       & $31.5608$    & $0.0239\%$   & $31.5544$ & $0.0443\%$ \\
$(1,4)$                  & $41.2817$                      & $41.2617$    & $0.0485\%$   & $41.2344$ & $0.1146\%$ \\
$(3,3)$                  & $43.7100$                      & $43.6975$    & $0.0286\%$   & $43.6542$ & $0.1277\%$ \\
$(2,4)$                  & $48.5667$                      & $48.5455$    & $0.0436\%$   & $48.4883$ & $0.1615\%$ \\
$(3,4)$                  & $60.7084$                     & $60.6821$    & $0.0434\%$   & $60.5840$ & $0.2048\%$ \\ \hline
\end{tabular}
\caption{ Comparison between the  frequencies estimated from the  numerical solutions on grid $\Gc_{160}$  and the  first $9$ natural frequencies. Discrepancies are measured using the percentage of the relative errors.
} \label{Tab:NaturalFrequencySSSSFreeVib}
\end{center}
\end{table}

{\bf $\bullet$ Forced vibration.}
Now, we consider the vibration of the  classical Kirchhoff-Love plate driven by a time-dependent  sinusoidal force
 $F(\xv,t)=F_0\sin(\xi t)$, where $F_0$ and $\xi $ are constants for the magnitude and frequency of the sinusoidal force. The plate  is assumed to be undeformed and at rest initially; that is, $w_0=v_0=0$.  Using method of eigenfunction expansion,  we find the exact solution to the  forced vibration problem,
\begin{equation}\label{eq:SupportedGeneralSol}
 w(x,y,t)=\sum_{m=1}^{\infty}\sum_{n=1}^{\infty}\sin\frac{m\pi x}{L}\sin\frac{n\pi y}{H}T_{mn}(t),
 \end{equation}
where the time-dependent coefficient $T_{mn}(t)$ is given by 
{\small$$
  T_{mn}(t)=\frac{2F_0 \left(1-\cos(m\pi)\right)\left(1-\cos(n\pi)\right)}{\rho h mn\pi^2\omega_{mn}}\left(\frac{\sin(\xi t)+\sin(\omega_{mn}t)}{\xi+\omega_{mn}}-\frac{\sin(\xi t)-\sin(\omega_{mn}t)}{\xi-\omega_{mn}}\right).
  $$
    }

For numerical test, we consider a classical  Kirchhoff-Love plate with the parameters specified as  $\rho h=1$, $K_0=0$, $T=0$, $D=0.1$, $K_1=0$, $T_1=0$ and $\nu=0.3$ on a rectangular domain   $\Omega=[0,0.4]\times [0,0.2]$;  namely $L=0.4$ and $H=0.2$.    The  magnitude and frequency  of the driving force are set  as   $F_0=1000$ and $\lambda=40$.   This test problem is solved using both of the proposed schemes on a uniform Cartesian grid with grid  spacings $h_x=h_y=1/300$.  The numerical results are compared  with an analytical solution truncated from the exact  solution \eqref{eq:SupportedGeneralSol} by keeping 49 modes (i.e., $m=1,\dots,7$ and $n=1,\dots,7$).

{
\newcommand{\figWidth}{6cm}
\def\xa{13.}
\def\ya{3}
\newcommand{\trimfig}[2]{\trimw{#1}{#2}{0.07}{0.}{0.2}{0.18}}
\begin{figure}[h]
\begin{center}
\begin{tikzpicture}[scale=1]
  \useasboundingbox (0.0,0.0) rectangle (\xa,\ya);  

\draw(-0.5,0.0) node[anchor=south west,xshift=0pt,yshift=0pt] {\trimfig{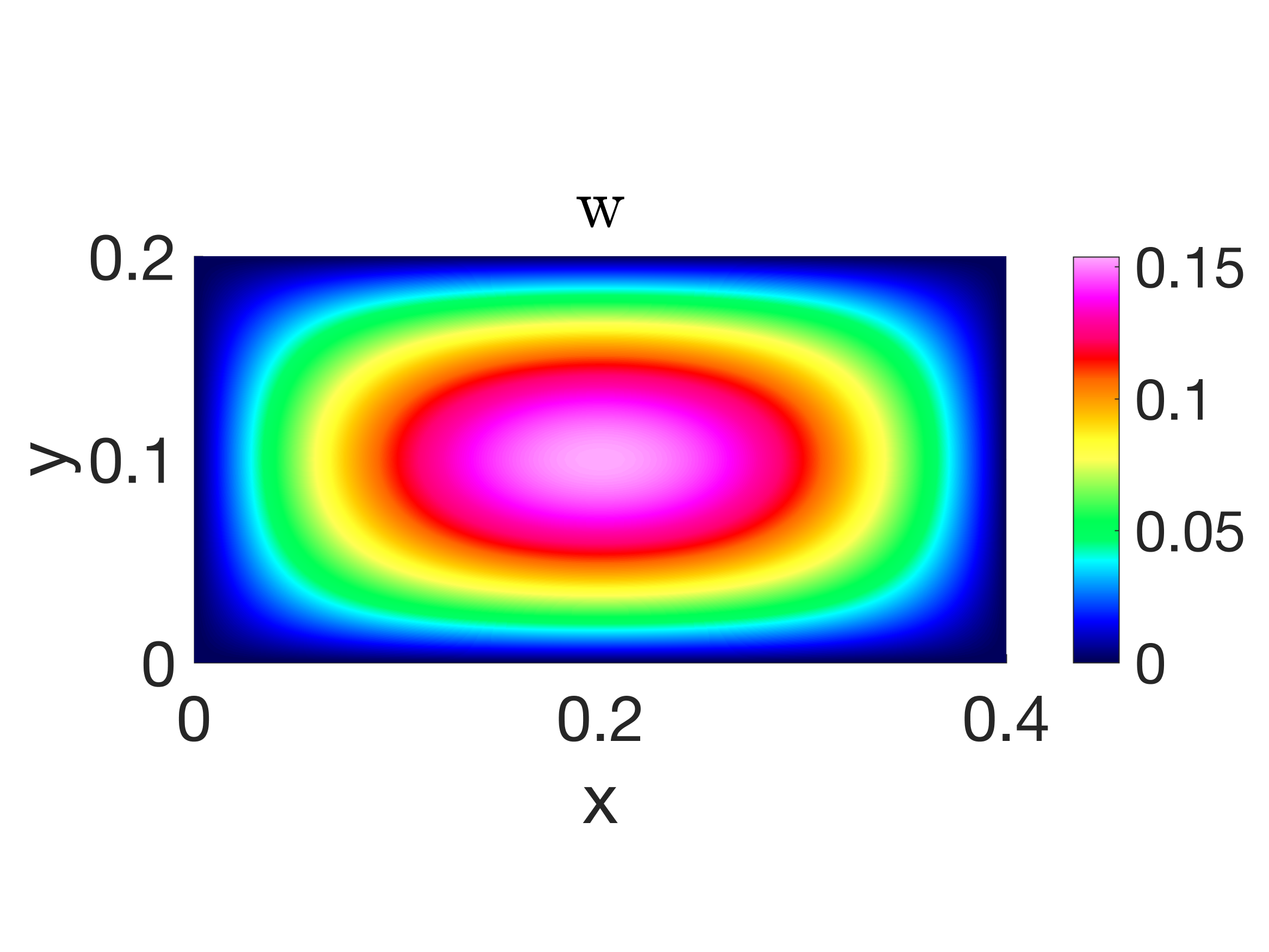}{\figWidth}};
\draw(6.5,0.0) node[anchor=south west,xshift=0pt,yshift=0pt] {\trimfig{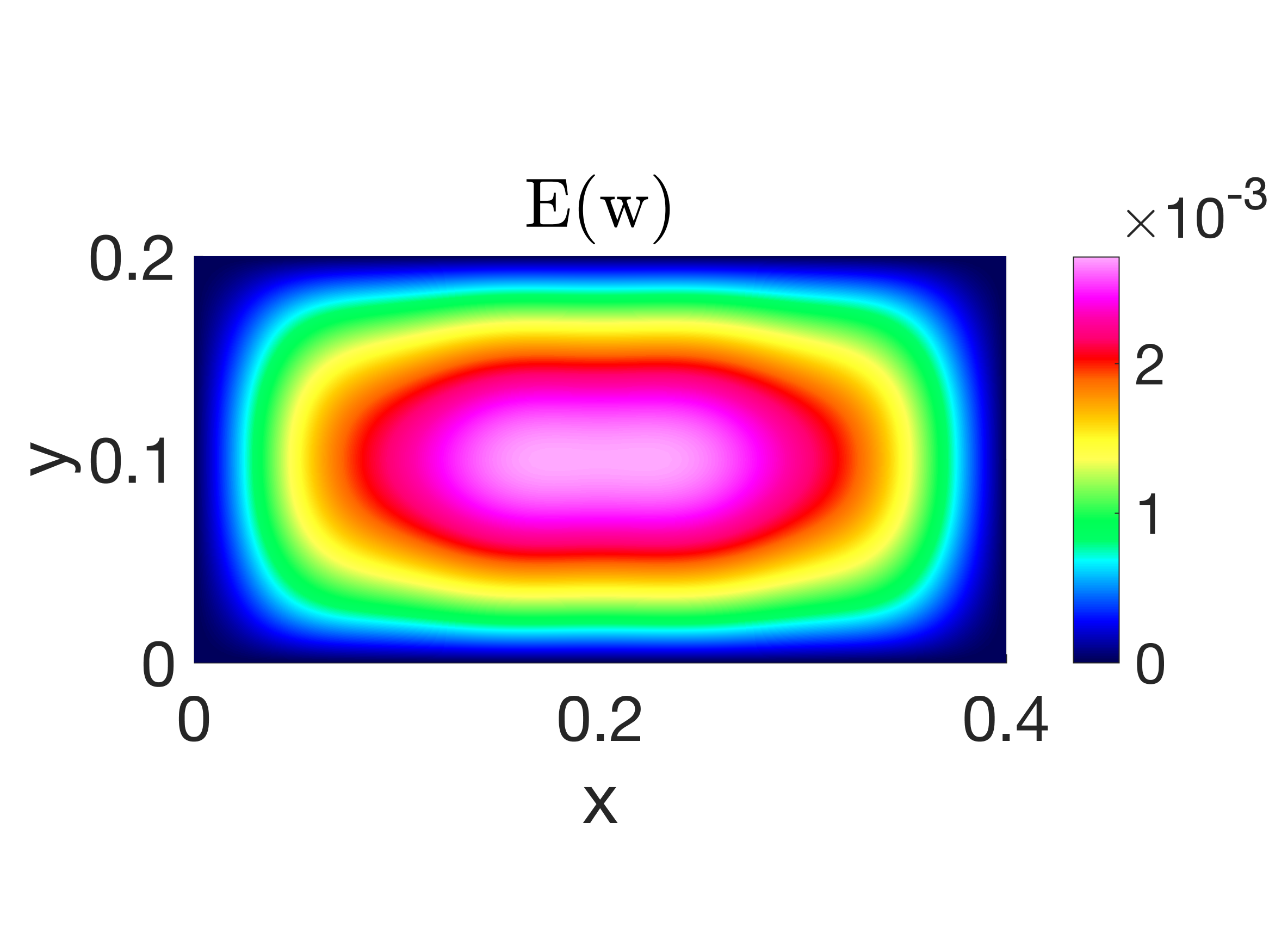}{\figWidth}};


%
\end{tikzpicture}

\end{center}
\caption{ Contour plots of  the numerical solution (left) and  the error (right) of $w$ at $t=1$. Results   are  obtained  using the NB2 scheme on a uniform Cartesian grid with grid spacings $h_x=h_y=1/300$. } \label{fig:SolErrorTimeHarmonicSupported}
\end{figure}
}

{
\newcommand{\figWidth}{6cm}
\def\xa{13.}
\def\ya{5}
\newcommand{\trimfig}[2]{\trimw{#1}{#2}{0.}{0.}{0.}{0.}}
\begin{figure}[h]
\begin{center}
\begin{tikzpicture}[scale=1]
  \useasboundingbox (0.0,0.0) rectangle (\xa,\ya);  

\draw(-0.5,0.0) node[anchor=south west,xshift=0pt,yshift=0pt] {\trimfig{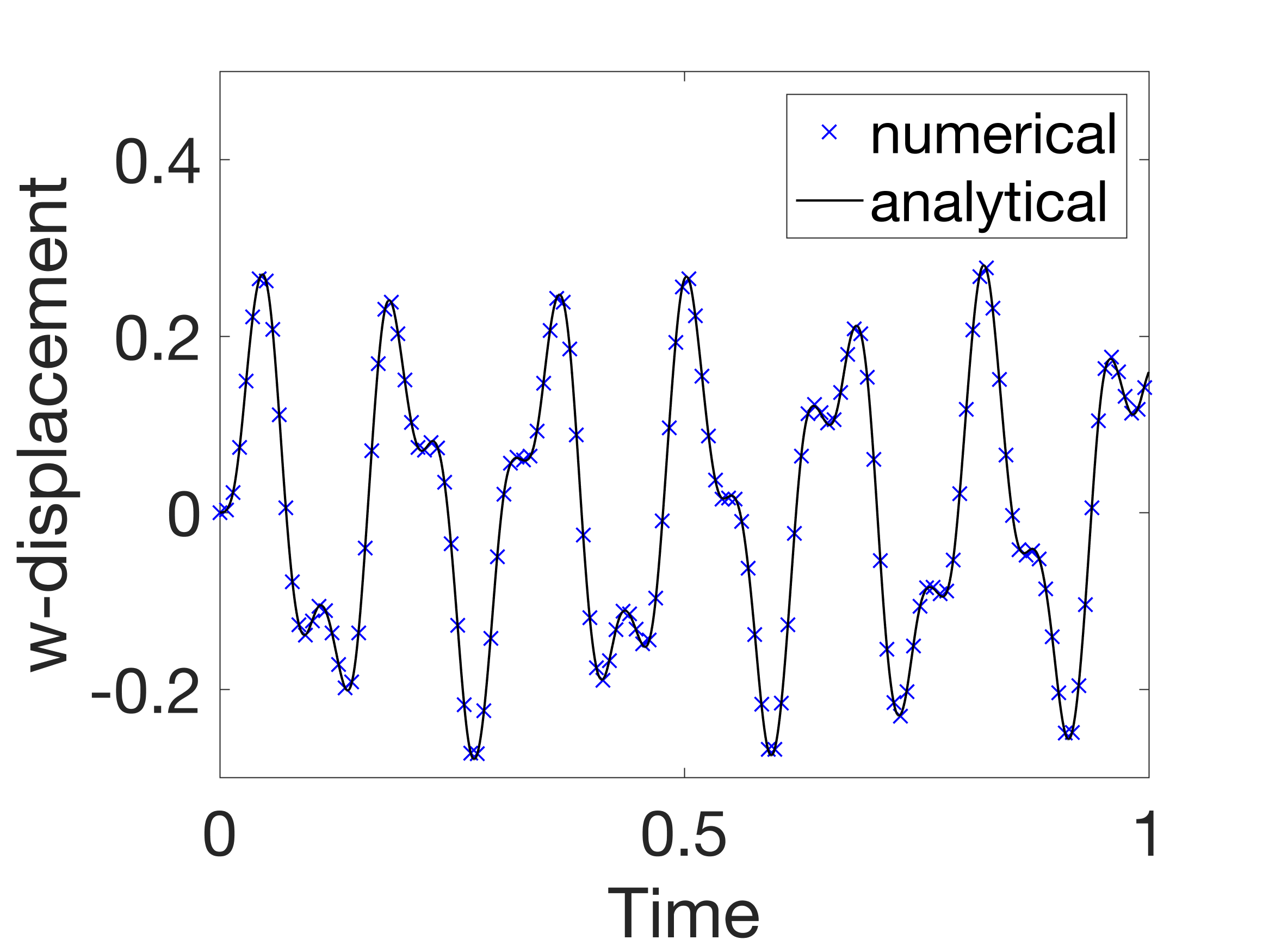}{\figWidth}};
\draw(6.5,0.0) node[anchor=south west,xshift=0pt,yshift=0pt] {\trimfig{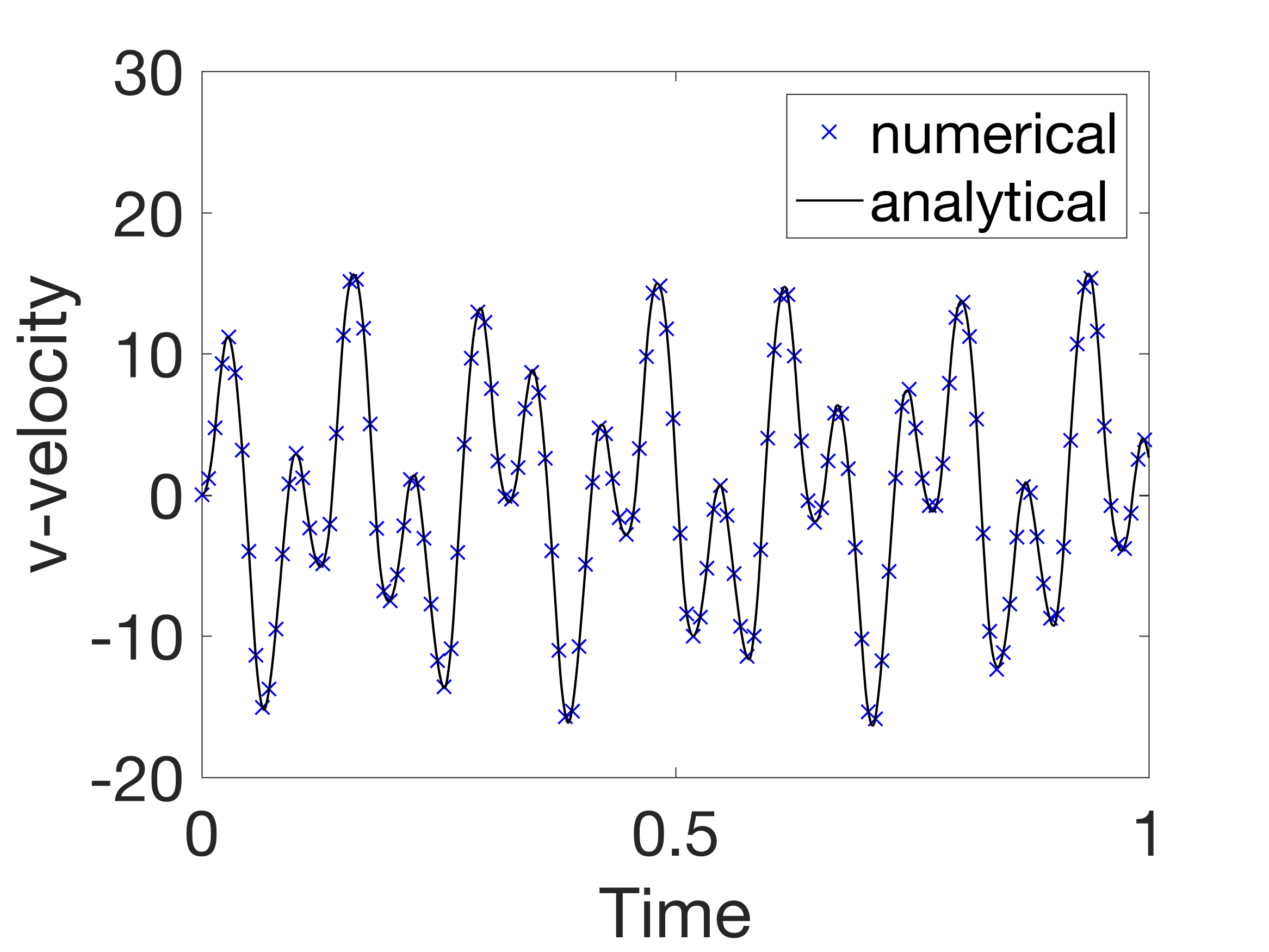}{\figWidth}};

\draw(3.2,4.5)  node[anchor=south,xshift=0pt,yshift=0pt] {Displacement};
\draw(10.3,4.5)  node[anchor=south,xshift=0pt,yshift=0pt] {Velocity};

%
\end{tikzpicture}

\end{center}
\caption{ The displacement and velocity of the plate at point $\xv_p=(0.2,0.1)$. Simulation is performed  using the NB2 scheme on a uniform Cartesian grid with grid spacings $ h_x=h_y = 1/300$.} \label{fig:TrackPointTimeHarmonicSupported}
\end{figure}
}

Both schemes perform comparably  well and produce similar solutions, so we only show the one  obtained with the NB2 schemes here.   In Figure~\ref{fig:SolErrorTimeHarmonicSupported}, the contour plots of    the numerical solution of $w$ at $t=1$  and the error compared with the truncated exact solution  are presented. To show the accuracy of the numerical results over time, we track the displacement $w$, as well the velocity $v$, at the point $\xv_p=(0.2,0.1)$. The time evolution of the numerical displacement and velocity at $\xv_p$ are plotted on top of the referenced analytical solutions 
in Figure~\ref{fig:TrackPointTimeHarmonicSupported}; it is easily seen that  our numerical results agree well with the analytical solution over time.

\subsubsection{Vibration of the generalized Kirchhoff-Love  plate}
For the generalized Kirchhoff-Love model \eqref{eq:generalizedKLPlate} that cannot be solved analytically, we  numerically solve  the frequency domain eigenvalue problem to identify the  natural frequencies and modes of vibration for  the plate, and then utilize the computed eigenvalues and eigenvectors to construct standing wave solutions.  The nodal line patterns (i.e., Chladni figures) of the standing wave solutions obtained from our numerical simulations are compared against those solved from the  eigenvalue problem  for the validation of the numerical schemes. 

Specifically, a square plate ($\Omega_S=[0,0.25]\times[0,0.25]$) and  an annulus  plate ($\Omega_A=\{\xv: 0.1\leq |\xv|\leq 0.5\}$) are considered. For both plates, we  assume the same  parameters,   $\rho h=1,K_0=2,T=1,D=2,K_1=0,T_1=0$, and $\nu=0.1$, noting that these parameters specify an undamped plate.  
On the edges of the plates, we impose the clamped  boundary conditions \eqref{eq:clampedBC} for square plate, and the simply supported boundary conditions \eqref{eq:supportedBC} for the annular plate.

First, let's consider the  eigenvalue problem for the undamped plate on mesh $\Omega_h$,
\begin{equation}\label{eq:eigenProblem}
 \Kc_h \phi_\iv=\lambda \phi_\iv, ~\xv_\iv\in\Omega_h,
\end{equation}
where  $\phi_\iv=\phi(\xv_\iv)$ is the mode function (or eigenfunction) for the eigenvalue $\lambda$. The definition of the difference operator   $\Kc_h$ is given in \eqref{eq:KBOperators}.  To get the eigenfunction-eigenvalue pairs $(\phi_n(\xv),\lambda_n)$, we numerically solve the eigenvalue problem \eqref{eq:eigenProblem} subject to  the appropriate numerical  boundary conditions; that is,  clamped  \eqref{eq:discreteClampedBC} for the square plate and supported \eqref{eq:discreteSupportedBC} for the annular plate.   The   \texttt{eigs} function   in MATLAB is used here.
To save space, we put the results in Appendix~\ref{sec:appendix}, where   the nodal lines of the first 25 eigenmodes (with multiplicity) for the plates (square and annular) are presented  in Figures~\ref{fig:ClampedModeShapeSquareMATLAB} \& \ref{fig:SupportedModeShapeAnnMATLAB}. Note that, following the tradition in structural engineering, the values of  natural frequencies, rather than the eigenvalues,   are reported in the plots.  The natural frequency corresponding to $\lambda_n$ is given by 
$$
f_n=\frac{1}{2\pi} \sqrt{\frac{\lambda_n}{\rho h}}.
$$

\begin{figure}[h]
  \begin{subfigure}[b]{0.325\linewidth}
    \centering
    $f_{11}=395.1906$
    \includegraphics[width=1\linewidth]{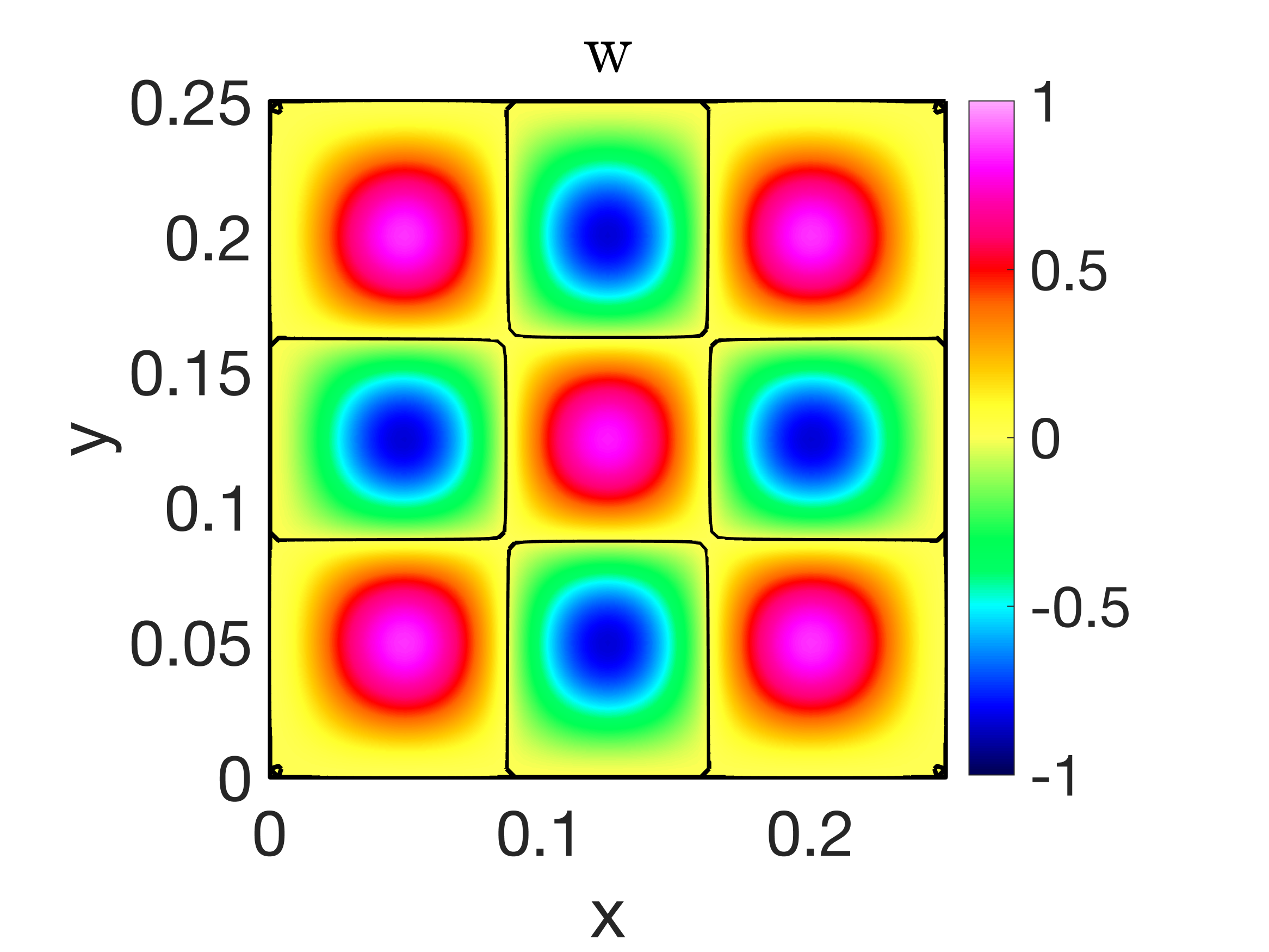}
  \end{subfigure} 
    \begin{subfigure}[b]{0.325\linewidth}
      \centering
      $f_{17}=553.9450$
    \includegraphics[width=1\linewidth]{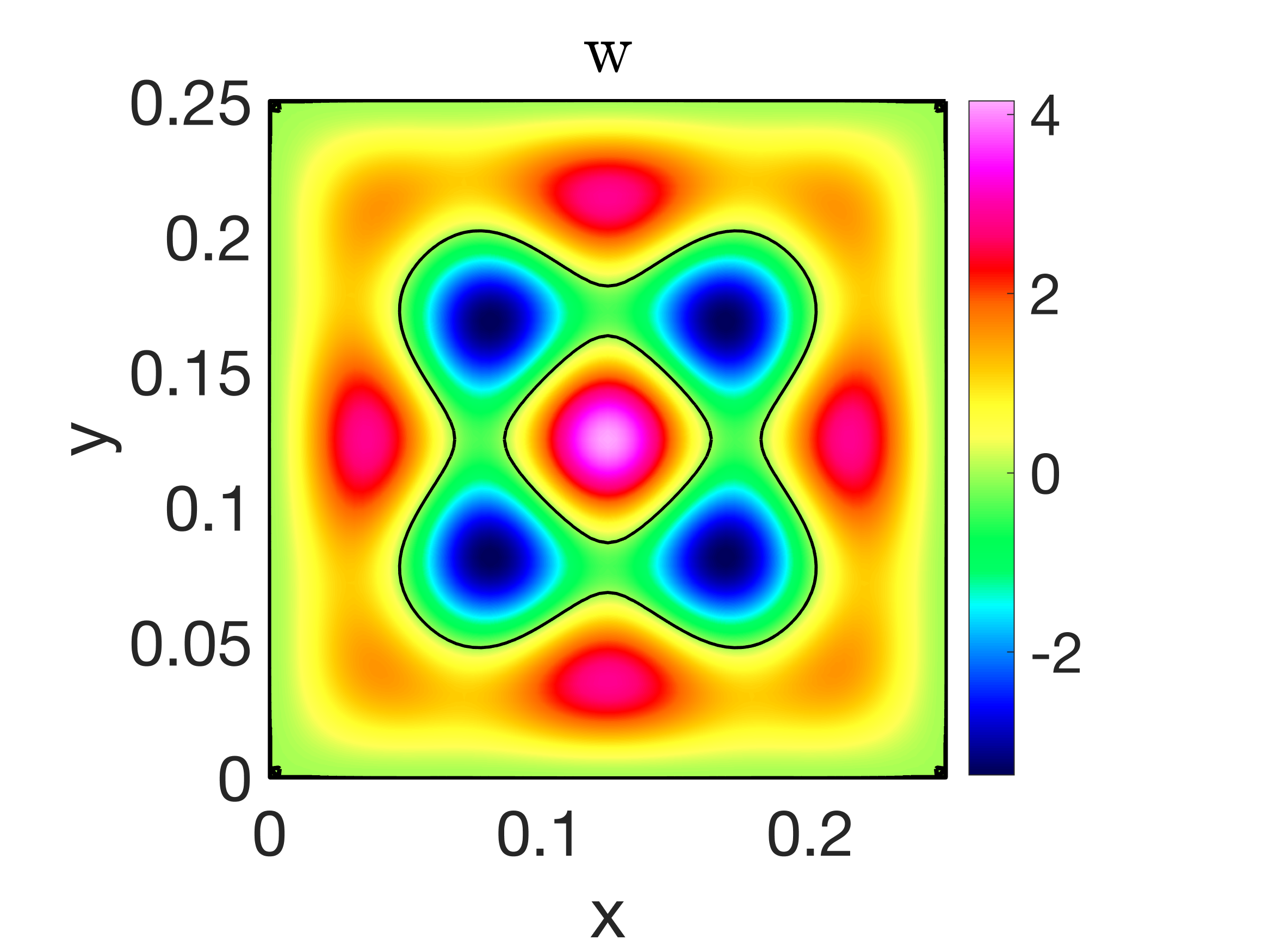}
  \end{subfigure}
  \begin{subfigure}[b]{0.325\linewidth}
    \centering
     $f_{22}=705.9267$
    \includegraphics[width=1\linewidth]{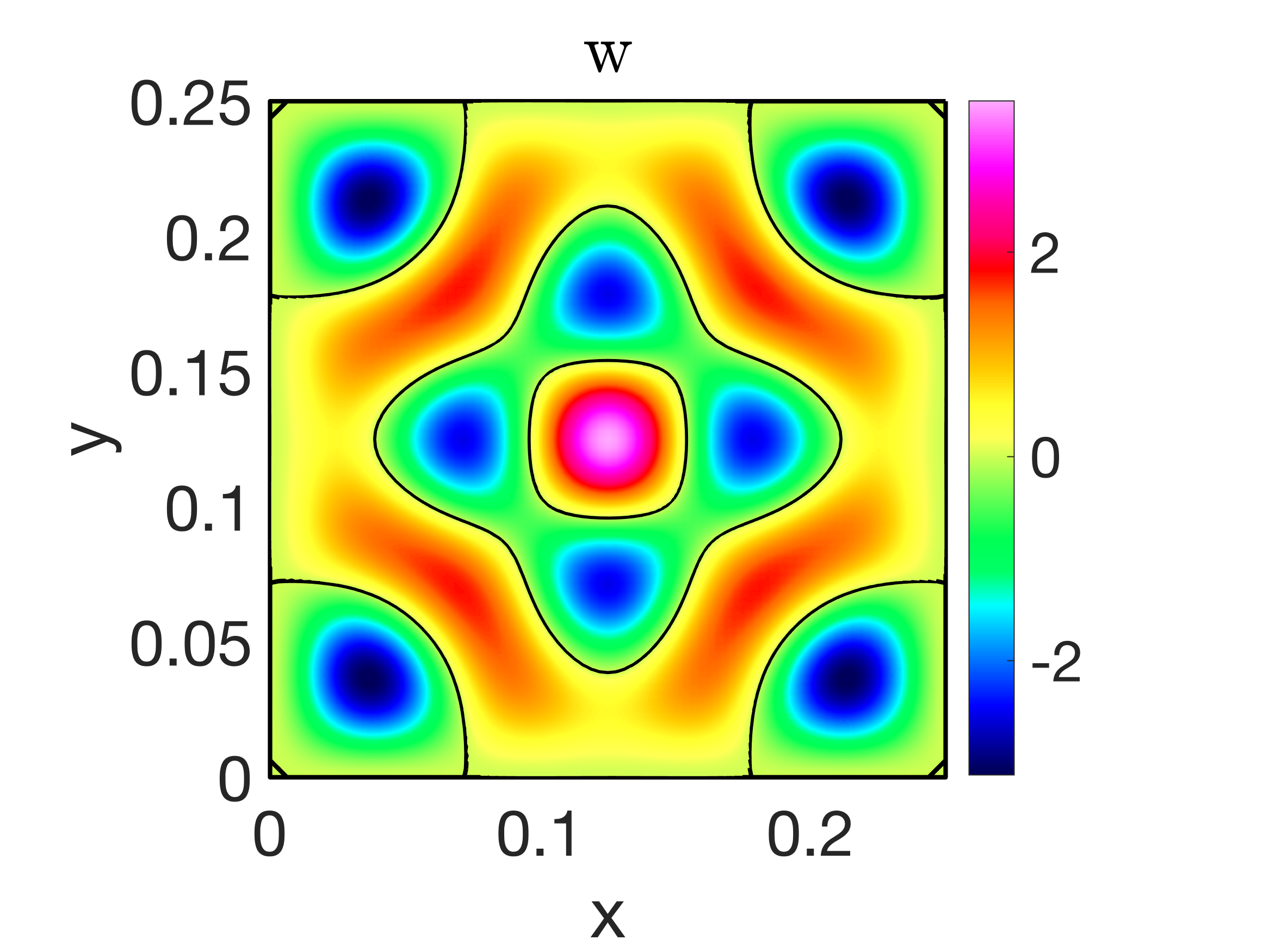}
  \end{subfigure} 
  \caption{ Standing waves in  the square plate with clamped edges  at time $t = 1$ for three eigenvalue cases. Zero contour lines of the solutions  are also plotted to indicate the nodal line patterns. Simulations are performed using the NB2 scheme  on grid $\Gc_{80}$.}
  \label{fig:ClampedModeShapeNB}
\end{figure}

\begin{figure}[h]
   \begin{subfigure}[b]{0.325\linewidth}
     \centering
     $f_{2}=8.8174$
    \includegraphics[width=1\linewidth]{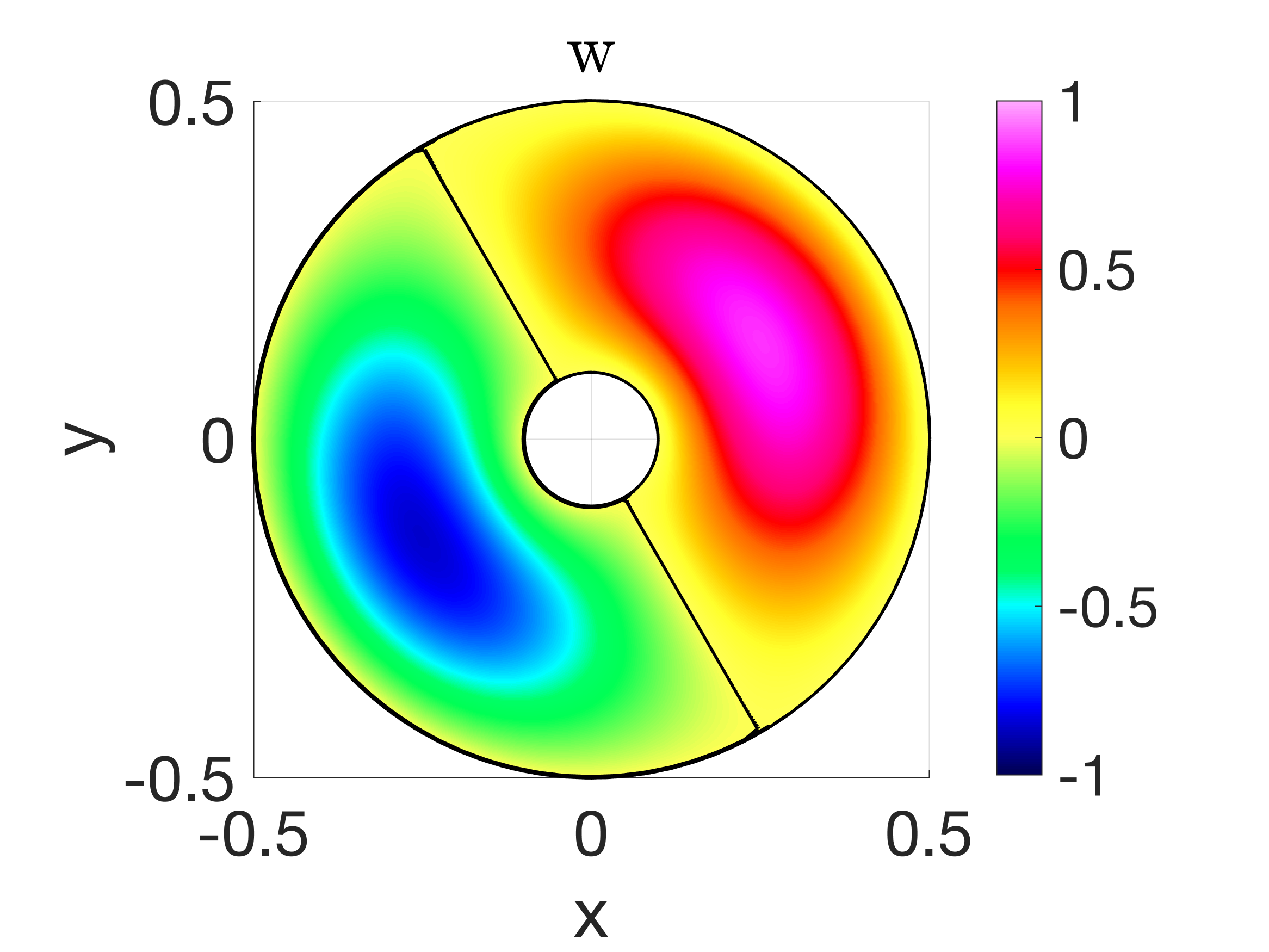}
  \end{subfigure} 
    \begin{subfigure}[b]{0.325\linewidth}
      \centering
      $f_{10}=28.9529$
    \includegraphics[width=1\linewidth]{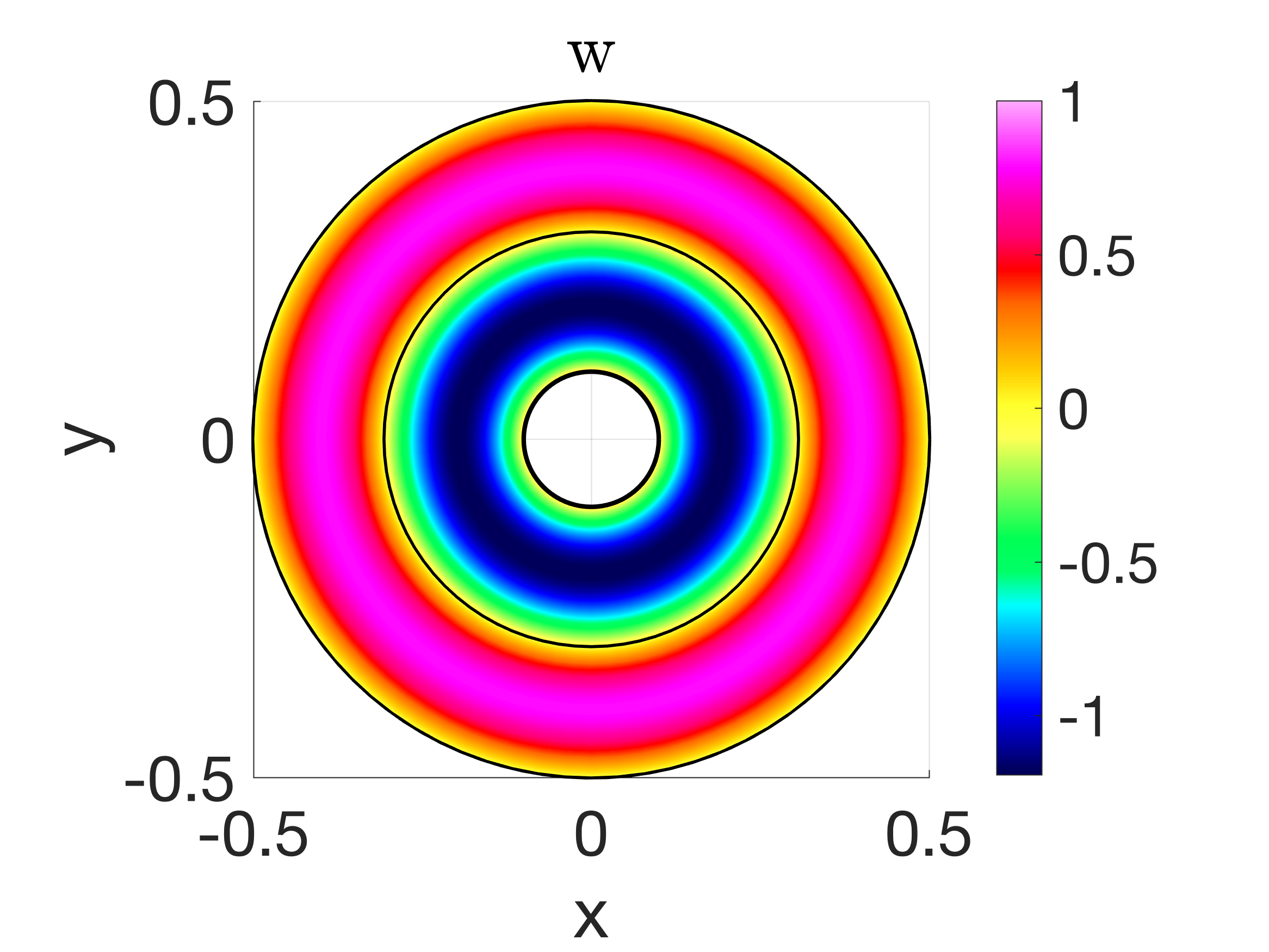}
  \end{subfigure}
  \begin{subfigure}[b]{0.325\linewidth}
    \centering
    $f_{20}=43.8879$
    \includegraphics[width=1\linewidth]{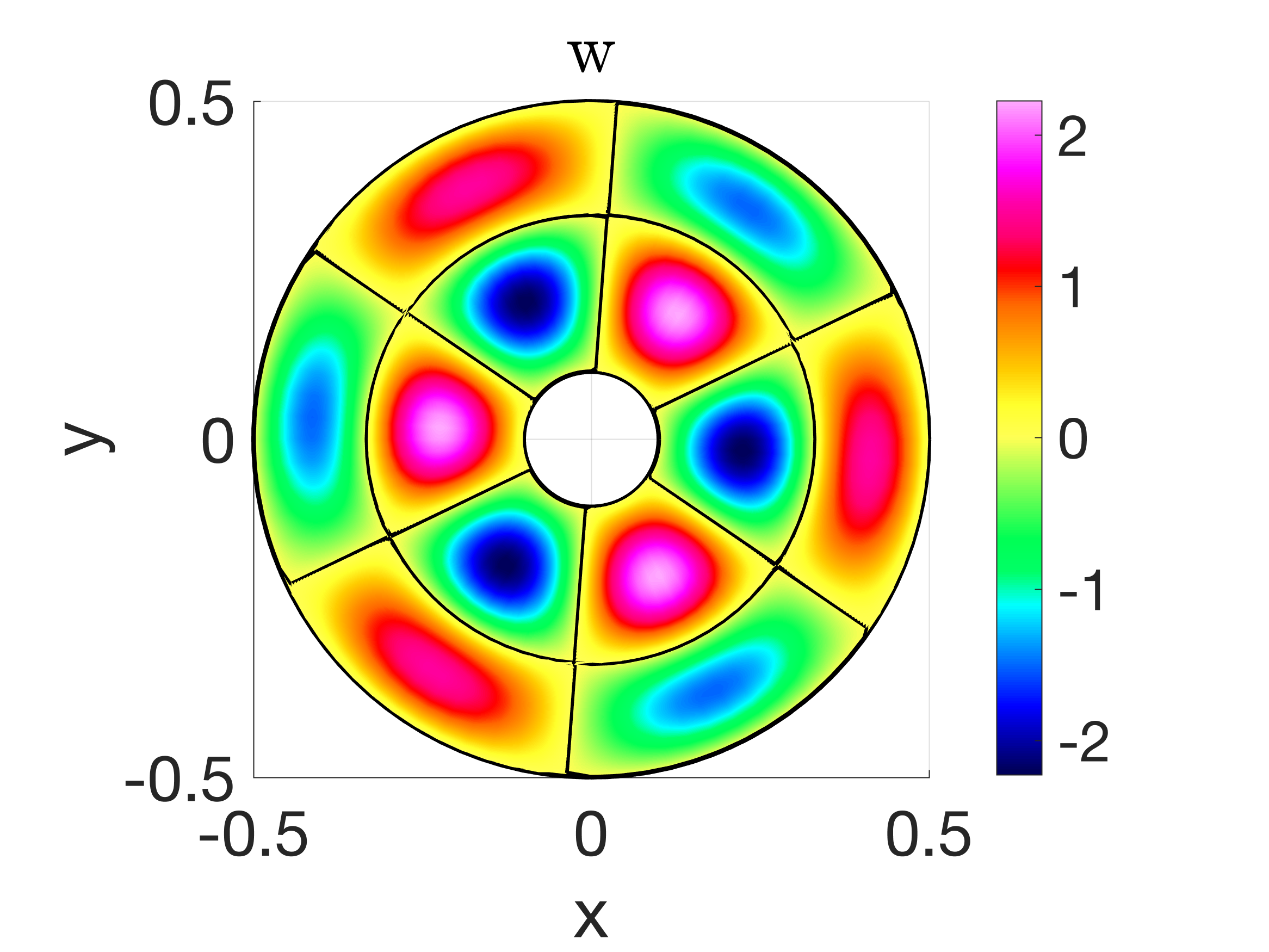}
  \end{subfigure} 
  \caption{ Standing waves in  the annular plate with simply supported edges  at time $t = 1$ for three eigenvalue cases. Zero contour lines of the solutions  are also plotted to indicate the nodal line patterns. Simulations are performed using the NB2 scheme  on grid $\Gc_{80}$.}
  \label{fig:AnnClampedModeShapeNB}
\end{figure}

Next, we solve for   standing waves in the  generalized Kirchhoff-Love model \eqref{eq:generalizedKLPlate} with the aforementioned parameters and   boundary conditions  numerically using the NB2 scheme.  Same  as before,  standing wave test problems are generated by assigning the initial conditions with the eigenmodes, 
$
w_0(\xv)=\phi_n(\xv) \quad\text{and}\quad v_0(\xv)=0,
$
and then  let the plate  vibrate freely (i.e., zero external forcing).
Results  for the square and annular plates are respectively shown  in Figure~\ref{fig:ClampedModeShapeNB} and Figure~\ref{fig:AnnClampedModeShapeNB};  three modes for each plate is selected  for presentation. Nodal lines of the numerical solutions are also plotted on top of the contour images. 
The fact that the nodal line patterns obtained from snapshots (at $t=1$) of  the numerical solutions to the dynamical PDE \eqref{eq:generalizedKLPlate}  clearly match those  solved from  the eigenvalue problem \eqref{eq:eigenProblem} (shown in  Figures~\ref{fig:ClampedModeShapeSquareMATLAB} \& \ref{fig:SupportedModeShapeAnnMATLAB} in  Appendix~\ref{sec:appendix})  is a strong evidence indicating the  accuracy of our numerical methods.

\subsubsection{Cross-validation with experiments}

As a final test, we compare our numerical results with  existing experimental results. In \cite{Tuan-2015}, Tuan {\em et al.} experimentally measured  the Chladni nodal line patterns and resonant frequencies for a thin plate excited by an electronically controlled mechanical oscillator. For one of their reported experiments, a thin  square plate with a length of $L=0.24$~m and  a  thickness of $h=0.001$~m was used.  The plate was made of   aluminum sheet that has the following  material parameters:  $E=69$~GPa, $\rho=2700$~kg/m$^3$ and $\nu=0.33$. The center of the plate  was fixed with a screw supporter that can be driven  with an electronically controlled mechanical oscillator. Silica sands  with grain size of $0.3$~mm were placed on the top surface of the plate. When the  oscillator drove the plate to vibrate at a resonant (natural) frequency,  the sand particles  stopped at the nodes of the resonant modes and therefore  manifested the nodal line patterns for the vibrating plate. 

For this test, we attempt to simulate  the experiment  and reconstruct comparable nodal line patterns numerically. To mimic the experiment, we consider the  Kirchhoff-Love plate  \eqref{eq:generalizedKLPlate}   on the square  domain, $\Omega=[0,0.24]\times [0,0.24]$. The  edges  of the plate are assumed to move freely, so the free boundary conditions \eqref{eq:freeBC} are applied; and  the center  of the plate is fixed, i.e., $w(\xv_c,t)=0$.  The parameters of the governing equation are chosen to represent  the material properties of the aluminum sheet; specifically, we set  $\rho h=2.7$, $K_0=0$, $T=0$, $D = 6.4527$, $K_1=0$, $T_1=0$, and $\nu=0.33$. The plate is assumed to be at rest and undeformed at $t=0$; that is, we have $w_0=v_0=0$ for the initial conditions \eqref{eq:IC}.

To account for the driving  force exerted  by the mechanical oscillator used in the experiment, we specify the external forcing of the model as   a time-dependent sinusoidal function that is none-zero on a small square area $\Ac$ at the center $(x_c,y_c)$ of the plate,
\begin{equation}\label{eq:Chladni1}
  F(\xv,t)=\begin{cases}
  F_0\cos\xi t,&  \xv\in \Ac \\
  0, & \xv\not\in \Ac \end{cases},
\end{equation}
where $F_0$ and  $\xi$ are the magnitude and the angular frequency of the driving force, and the square area is  $\Ac=\left[x_c-0.01,x_c+0.01\right]\times \left[y_c-0.01,y_c+0.01\right]$.


To reconstruct  the nodal line patterns numerically, we need  the driving force \eqref{eq:Chladni1} to oscillate   at a resonant (natural) frequency. So we  first  solve the  eigenvalue problem \eqref{eq:eigenProblem} using the \texttt{eigs} function in MATLAB to  find the natural frequencies for  the model plate. The relation between  the $k$th  resonance angular frequency   and the corresponding  eigenvalue is $\xi_k=\sqrt{{\lambda_k}/({\rho h})}$. We  choose the magnitude of the force to be $F_0=10^{10}$, and perform the simulations using the NB2 scheme on grid $\Gc_{160}$. The reason why  such a large magnitude is used for the driving force is that we hope to quickly  force the plate to vibrate in the resonant mode.

\begin{figure}[h!] 
  \begin{subfigure}[b]{0.32\linewidth}
    \centering
    {$f_1=609.7$Hz}
    \includegraphics[width=1\linewidth]{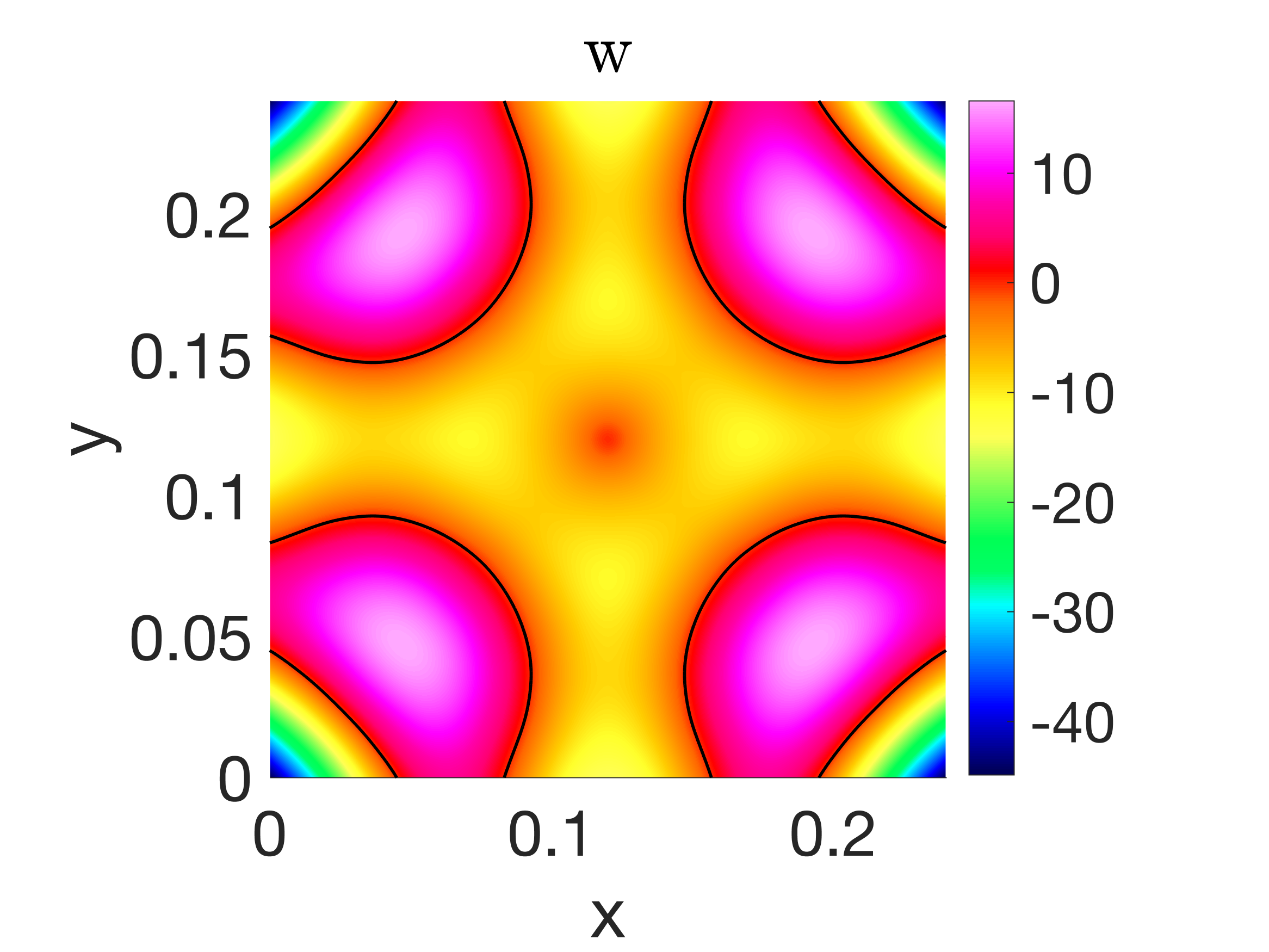}
  \end{subfigure}
    \bigskip
  \begin{subfigure}[b]{0.32\linewidth}
    \centering
    {$f_2=995.0$Hz}
    \includegraphics[width=1\linewidth]{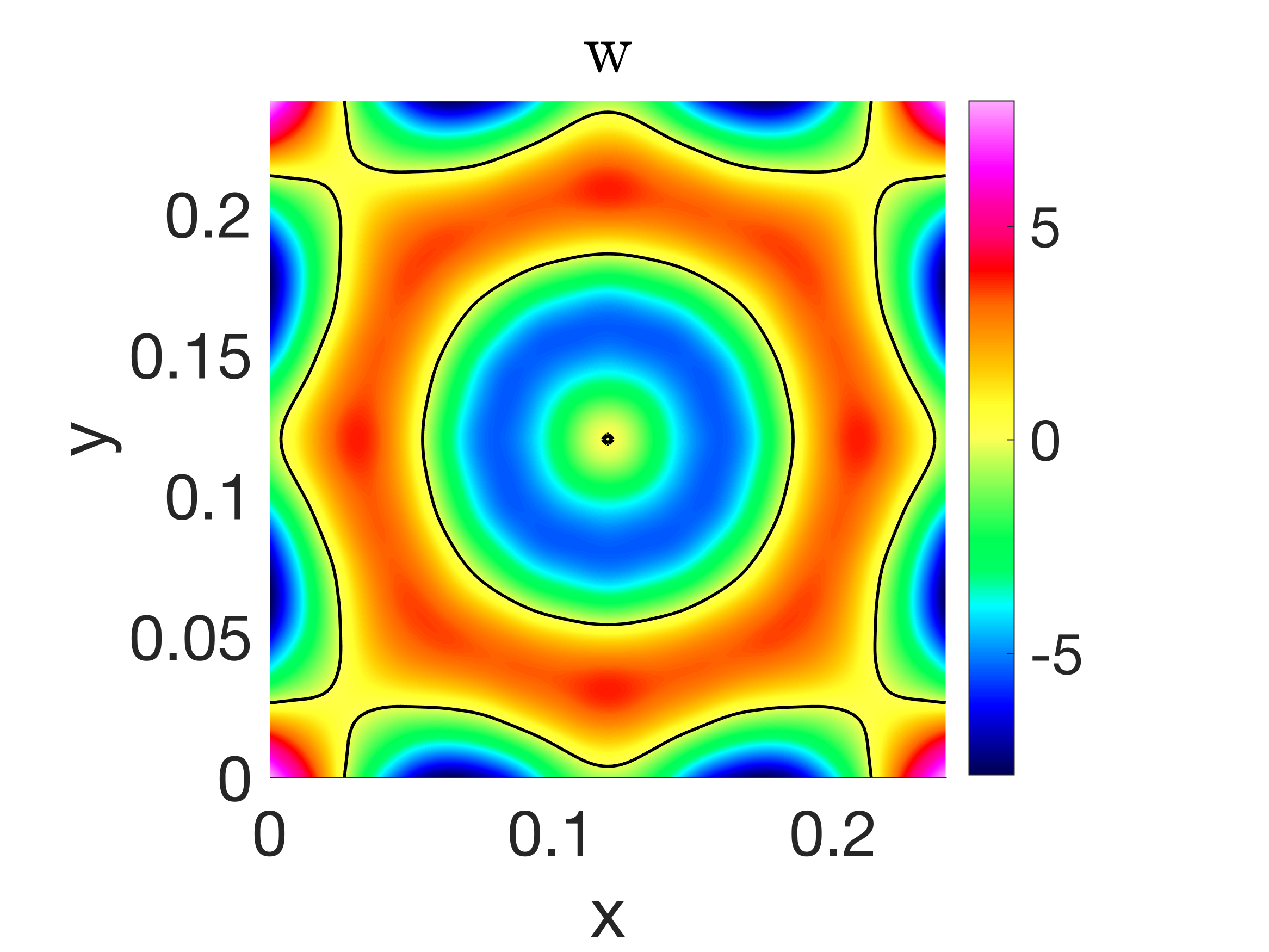}
  \end{subfigure} 
  \begin{subfigure}[b]{0.32\linewidth}
    \centering
     {$f_{10}=3443.9$Hz}
    \includegraphics[width=1\linewidth]{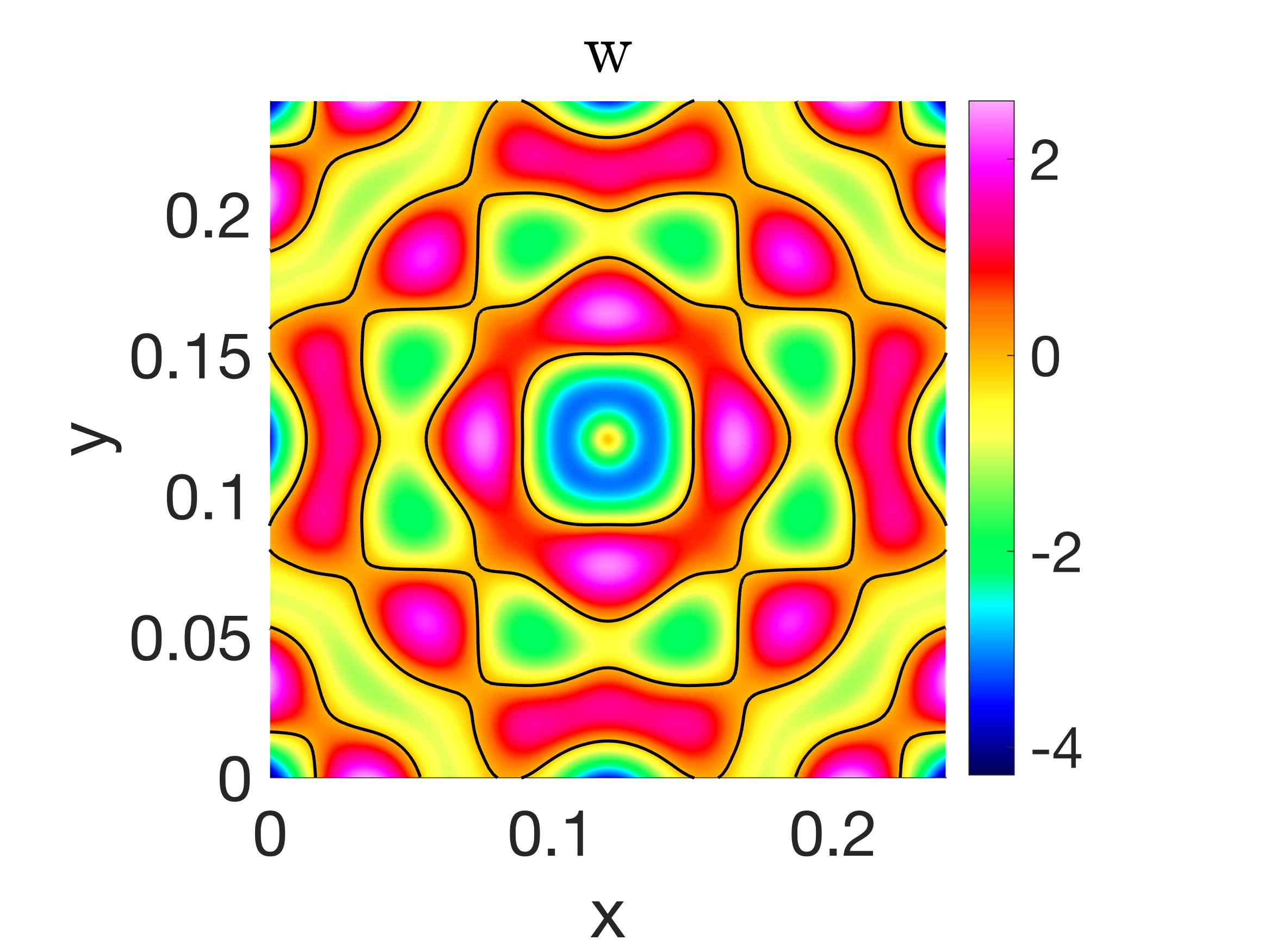}
  \end{subfigure} 
  \caption{ Contour plots of the displacement and the nodal lines (i.e., zero contours)  at $t=1$ for three typical  resonant frequencies. The  results shown here are obtained  from  the simulation of the NB2 scheme on grid $\Gc_{160}$.}
  \label{fig:ChladniFreeModeShapeNB}
\end{figure}

The  results obtained  from  the simulations using the NB2 scheme on grid $\Gc_{160}$ are presented in Figure~\ref{fig:ChladniFreeModeShapeNB}, which includes contour plots of the displacement and the nodal lines at $t = 1$ for three typical resonant frequencies. 
In order to directly compare  with the regular  frequencies  reported in the experiment,  we also report in Figure~\ref{fig:ChladniFreeModeShapeNB} the regular frequencies that are converted from  the angular frequencies by $f_k={\xi_k}/{2\pi}$.
The nodal line patterns manifested in our numerical results are in excellent agreement to the experimental results for all the frequencies (except for the degenerate eigenvalues); note that the specific experimental results we are comparing with are reported in Figure~3b in \cite{Tuan-2015}.

It is worth pointing out that there are noticeable discrepancies between the  values of the  numerical and experimental  resonant frequencies. This is because the plate model we used for this test is a simple classical  Kirchhoff-Love plate that does not consider the  influences of the ambient air and   the extra mass  of the sand particles on  the plate. The model could be improved by tuning the various parameters in \eqref{eq:generalizedKLPlate}, which is beyond the scope of this study.  The purpose  of this test is to validate our numerical methods; and  the fact that  the  numerically reconstructed  resonant nodal line patterns agree well with the experimental ones and the values of resonant frequencies are in qualitative agreement serves that purpose.

\subsection{Application}  As an application of our schemes, we   numerically explore the interesting physical phenomena known as   resonance  and beat that occur when the driving frequency is  right at  or close to  a natural frequency.  Here we  demonstrate the application by considering an annular plate ($\Omega=\{\xv: 0.1\leq|\xv|\leq 0.5\}$) with no external forcing as an example. The plate satisfies the generalized  Kirchhoff-Love equation \eqref{eq:generalizedKLPlate} and is driven to vibrate by  the following time-dependent   clamped boundary conditions that prescribe the displacement of  the inner edge; i.e.,
\begin{equation} \label{eq:ModifiedClampedBC}
w(\xv,t)=W_\text{in}\cos\xi t, \quad \pd{w}{\nv}(\xv,t)=0 ~\text{for}~|\xv|=0.1,
\end{equation}
where $W_\text{in}$ and $\xi$   are the  maximum value (amplitude)  and  angular frequency of the prescribed boundary displacement, respectively.  For this example, we set $W_{in}=1$ and vary   $\xi$ to investigate its effects.
The outer edge of the plate is allowed to   move freely; namely,  the free boundary conditions \eqref{eq:freeBC} are applied  at $|\xv|=0.5$.  Initially, we assume $w_0(\xv,0)=v_0(\xv,0)=0$. The setup of this problem can easily be replicated experimentally  by clamping the inner edge of an annular plate to a  mechanical oscillator undergoing a sinusoidal motion.
The material parameters of the plate are  $\rho h =1,D=0.01,T=0,K_0=0$, and $\nu=0.3$. With the  intention to study the effects of the damping terms in the equation, various values for $K_1$ and $T_1$ are considered below.

To begin with, we consider the undamped case (i.e.,  $K_1=T_1=0$), and specify a value to  $\xi$ in \eqref{eq:ModifiedClampedBC} that is either close to or at a natural frequency of the plate.   Due to the complexity of the generalized plate equation and the time-dependent boundary conditions, it is non-trivial to analytically find the  frequency domain  eigenvalue problem  and then solve it  for the natural frequencies and modes  as were done in Section~\ref{sec:vibrationOfPlates}. Following a procedure  proposed in  \cite{Chugh2007},   we illustrate a more general  strategy  to  identify the natural frequencies of  a plate using our numerical methods in conjunction with a fast Fourier transformation (FFT) power spectrum analysis of the numerical data.

\begin{figure}[h] 
  \centering
  \begin{subfigure}[b]{0.49\linewidth}
    \centering
    $w(\xv_p,t)$
    \includegraphics[width=1\linewidth]{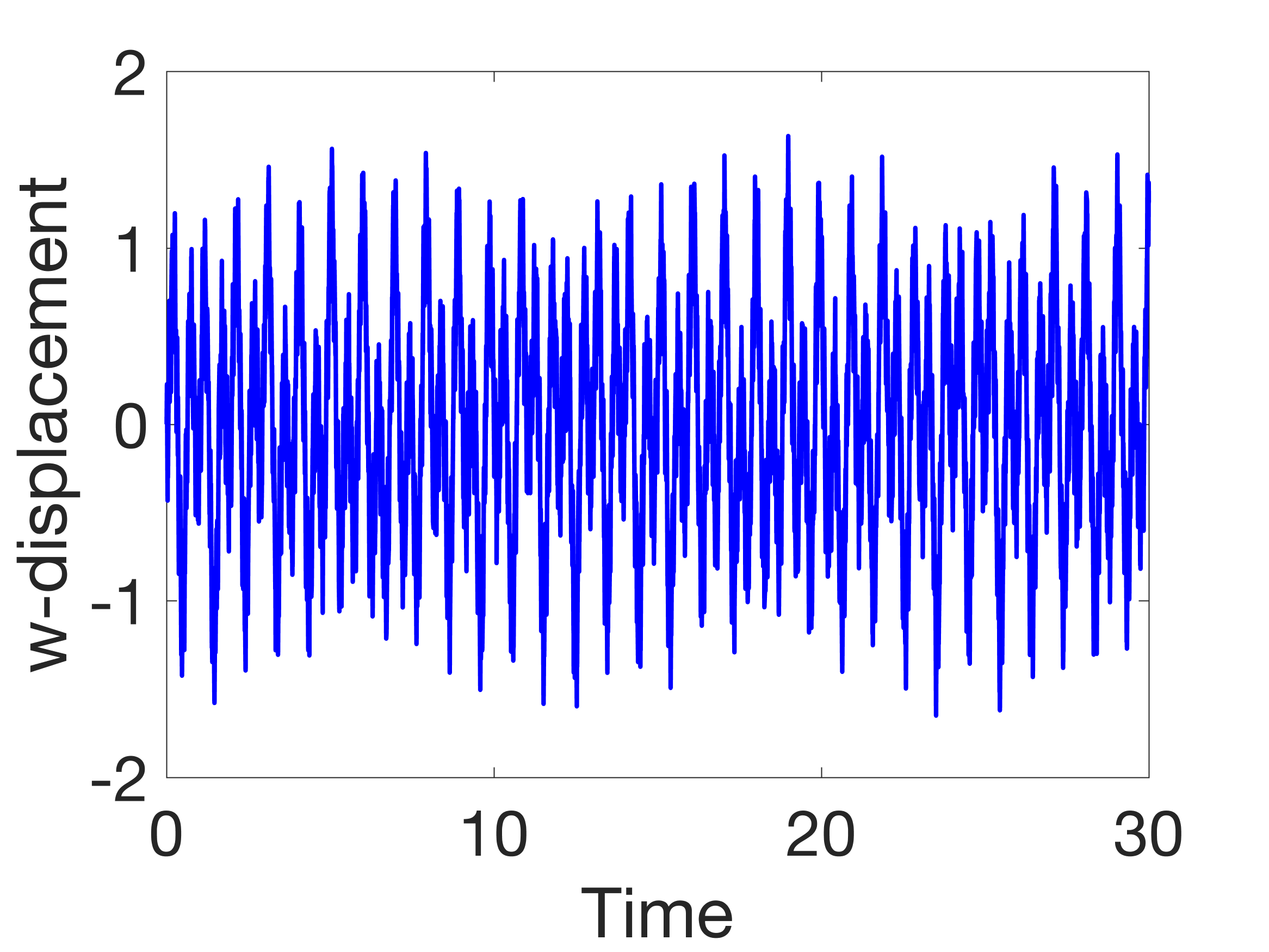}
  \end{subfigure}
    \bigskip
  \begin{subfigure}[b]{0.49\linewidth}
    \centering
    FFT of $w(\xv_p,t)$
    \includegraphics[width=1\linewidth]{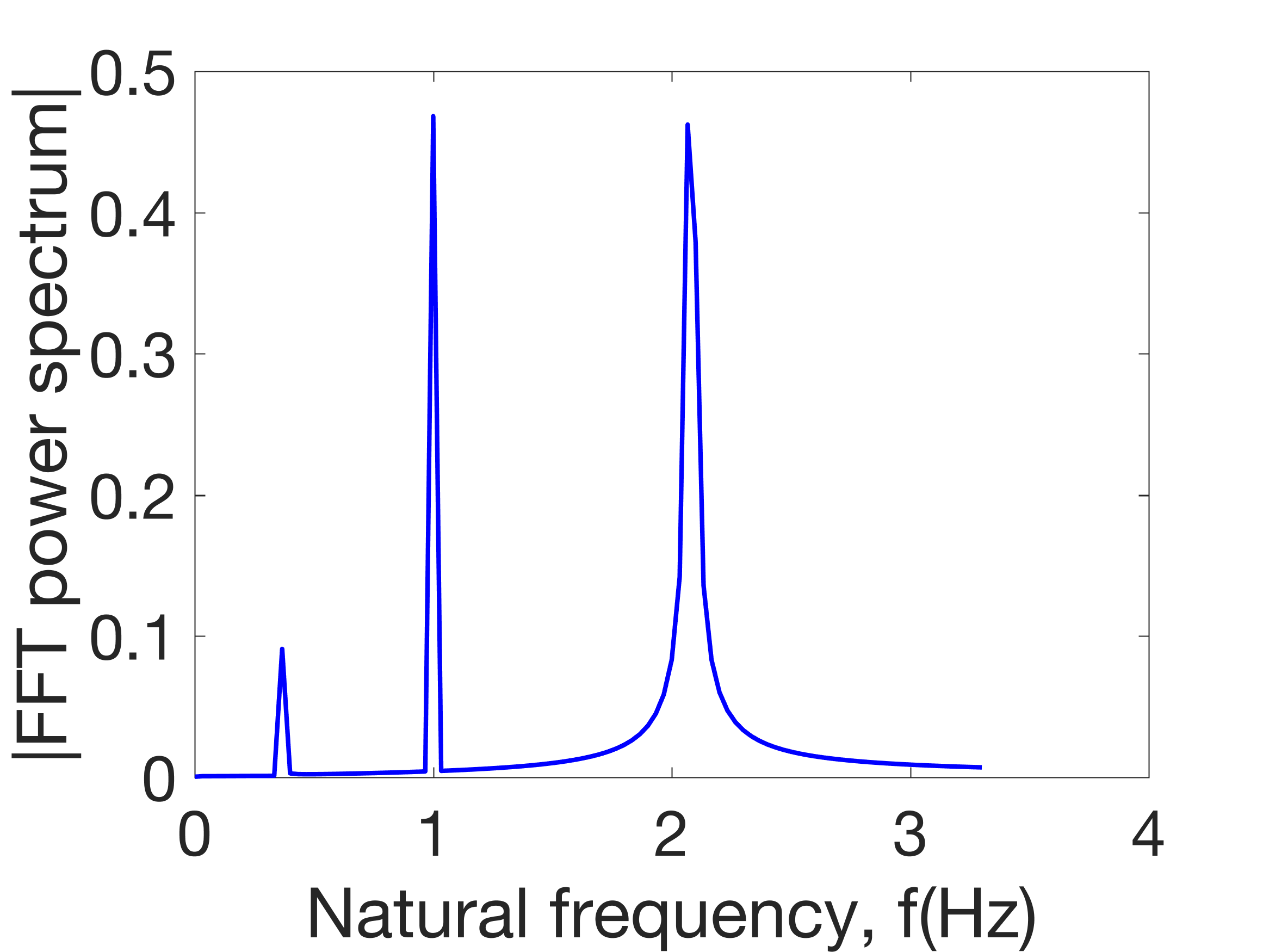}
  \end{subfigure} 
  \caption{  Vibration  of the plate when the driving frequency is $\xi=2\pi$ (or $f_d=1$~Hz). The simulation  runs until $t=30$ using the NB2 scheme on grid $\Gc_{80}$. 
    Left: the displacement response at the selected location $\xv_p=(-0.2,0)$ over time. Right: the FFT power spectrum of the  displacement response.}
       \label{fig:MovingClampedOmega1}%
\end{figure}

The strategy for finding a natural frequency goes as following. We first  simulate  the problem for an arbitrary driving frequency; say, $\xi=2\pi$ (or $f_d=1$~Hz), and   trace  the response of the plate at $\xv_p=(-0.2,0)$. The simulation  runs until $t=30$ using the NB2 scheme  on grid $\Gc_{80}$.    The left image of Figure~\ref{fig:MovingClampedOmega1} shows the displacement response at the selected location over time.
We then perform  FFT to the   displacement data using the \texttt{fft} function in MATLAB, and  present its power spectrum  in the right image of  Figure~\ref{fig:MovingClampedOmega1}. 
From this graph, we are able to identify two natural  frequencies ($f_1=0.367$,  $f_2=2.067$) and the driving frequency ($f_d=1$). More natural frequencies can be identified this way by sampling   different  values  for the driving frequency $\xi$. 

 \begin{figure}[h]    
    \begin{subfigure}[b]{0.32\linewidth}
    \centering
    \includegraphics[width=1\linewidth]{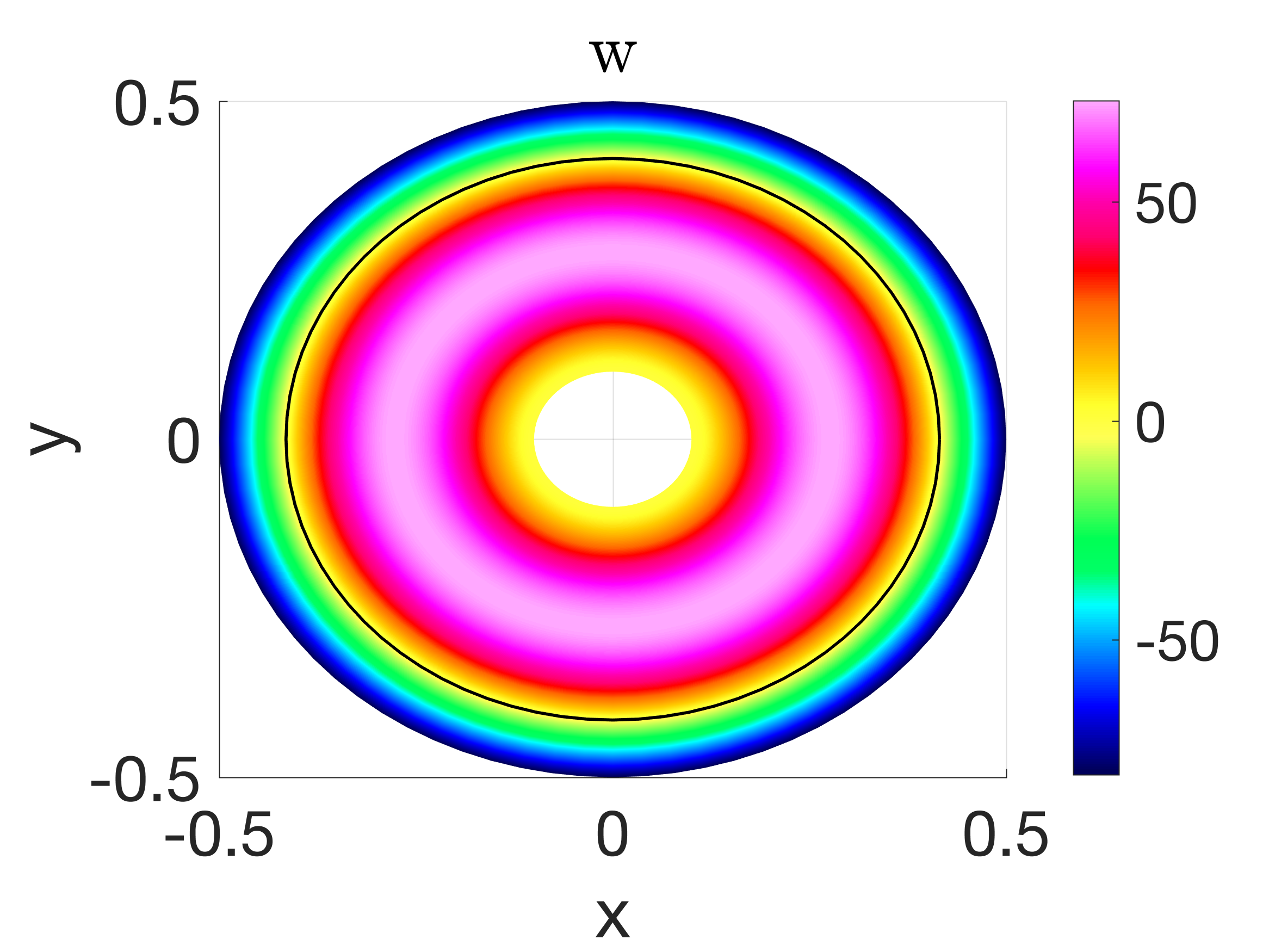}
    \caption{\footnotesize   contour of  $w(\xv,30)$} \label{fig:MovingClampedResonanceA}
  \end{subfigure} 
    \begin{subfigure}[b]{0.32\linewidth}
    \centering
    \includegraphics[width=1\linewidth]{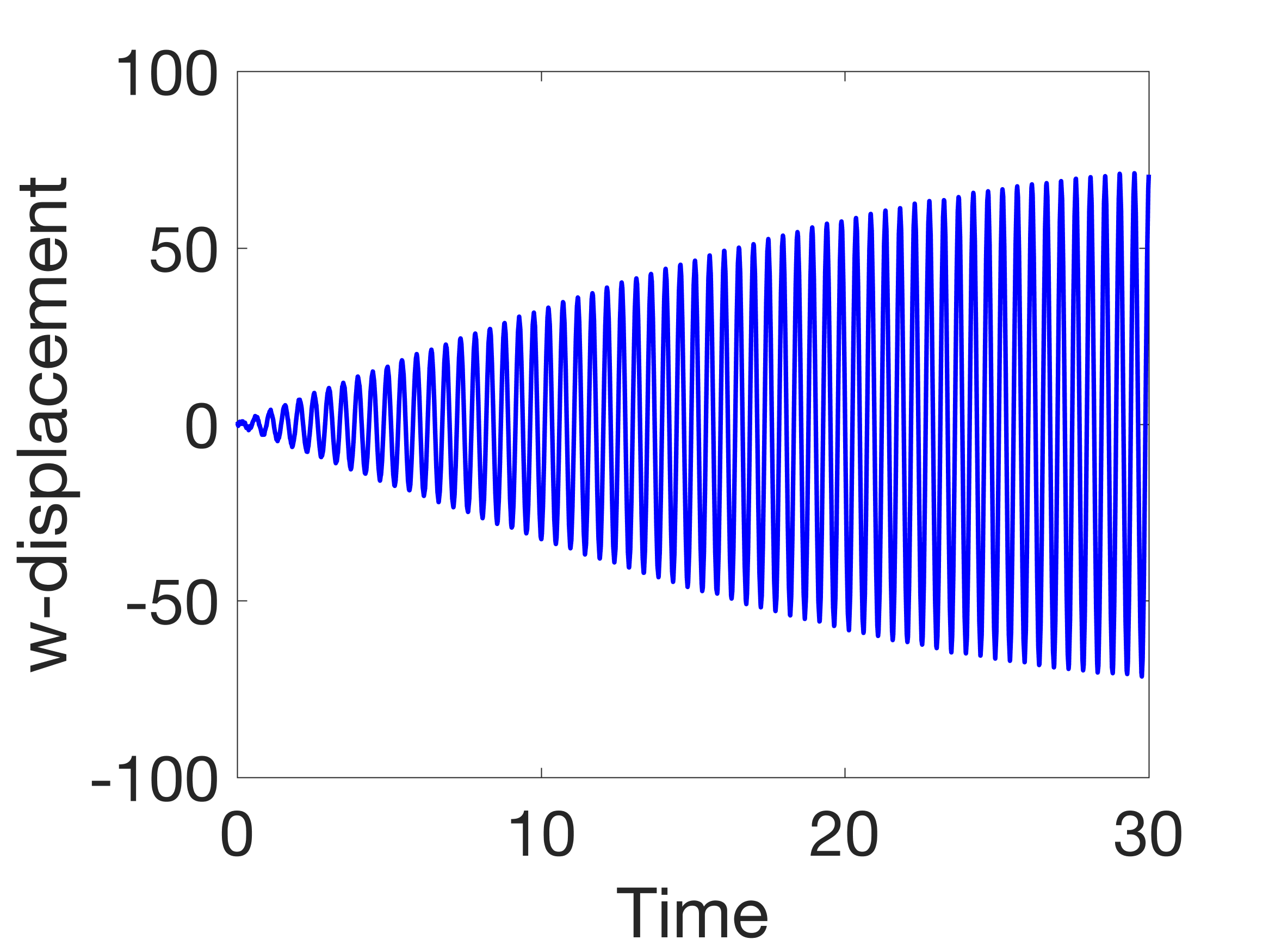}
    \caption{\footnotesize  $w(\xv_p,t)$}  \label{fig:MovingClampedResonanceB}
  \end{subfigure} 
      \begin{subfigure}[b]{0.32\linewidth}
    \centering
    \includegraphics[width=1\linewidth]{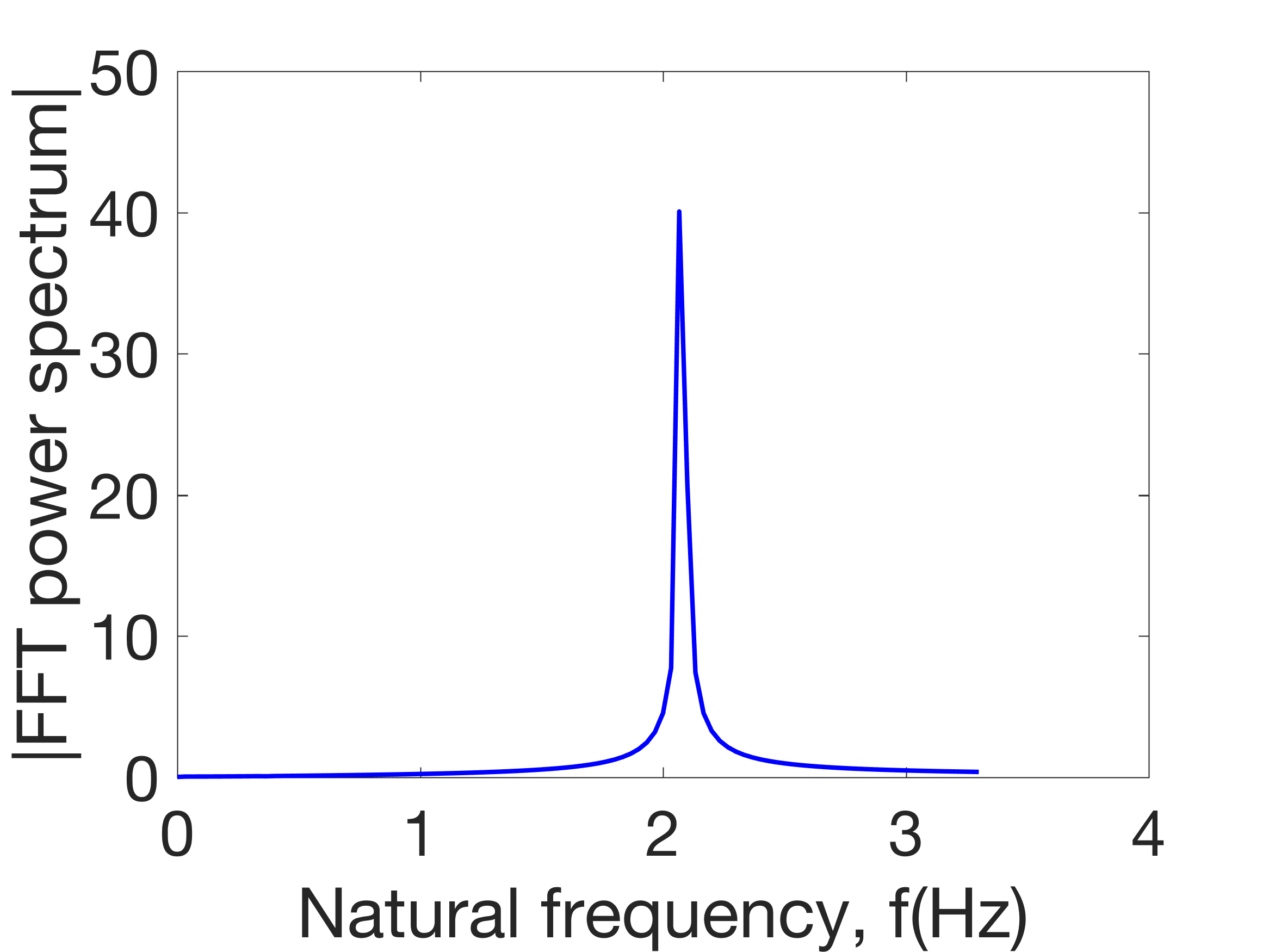}
    \caption{\footnotesize  FFT of  $w(\xv_p,t)$}  \label{fig:MovingClampedResonanceC}
  \end{subfigure} 
  \caption{ Numerical simulation of the resonance phenomenon  using the NB2 method on grid $\Gc_{80}$. Here the driving frequency is $\xi_2=2\pi f_2$ and the probed location is $\xv_p=(-0.2,0)$.}      
  \label{fig:MovingClampedResonance}%
 \end{figure}
 
 Resonance occurs when the driving frequency of the plate is at a natural frequency.  As an example, we simulate  this phenomenon at the natural frequency $f_2$  by setting  the driving frequency as  $\xi_2=2\pi f_2$.    The simulation is carried out using the NB2 scheme until $t=30$, and  the  results are collected in  Figure~\ref{fig:MovingClampedResonance}. In particular, we show in Figure~\ref{fig:MovingClampedResonanceA} the contour plot of $w$  as well as its nodal lines at $t=30$. The nodal line pattern sheds light on the mode shape (eigenfunction) associated with  the natural frequency $f_2$. We also trace the displacement  at the point $\xv_p=(-0.2,0)$ and depict its time history in Figure~\ref{fig:MovingClampedResonanceB}. The resonance phenomenon is clearly observed  as  the amplitude of the vibration increases  over time.  The FFT power spectrum of the displacement data at this point, as is  shown  in Figure~\ref{fig:MovingClampedResonanceC}, also confirms that the plate vibrates at a frequency  consistent with   the natural  frequency  $f_2$.

   \begin{figure}[h]
  \begin{subfigure}[b]{0.32\linewidth}
    \centering
    \includegraphics[width=1\linewidth]{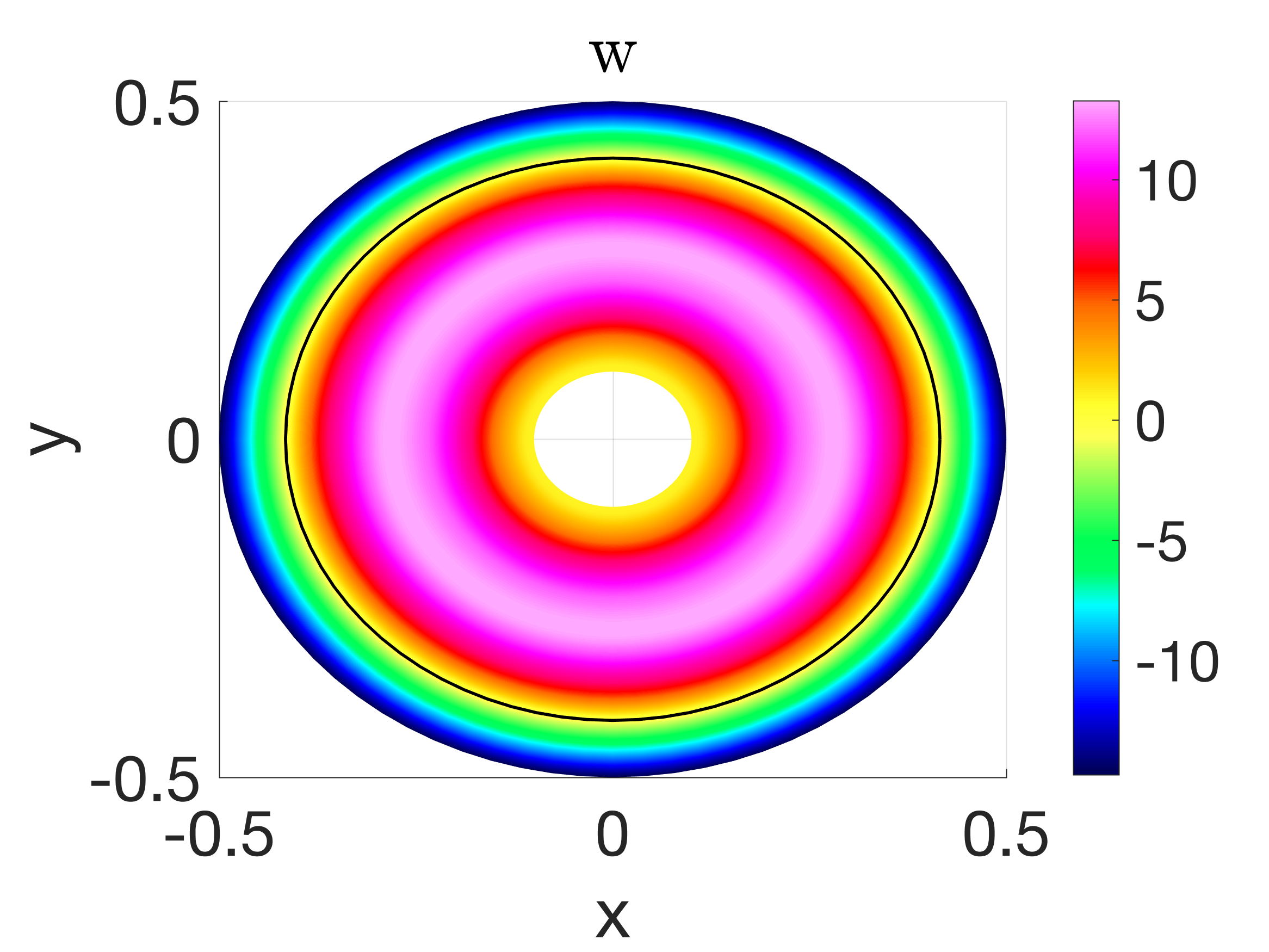}
    \caption{\footnotesize  contour of $w(\xv,30)$}  \label{fig:MovingClampedBeatA}%
  \end{subfigure} 
    \begin{subfigure}[b]{0.32\linewidth}
    \centering
    \includegraphics[width=1\linewidth]{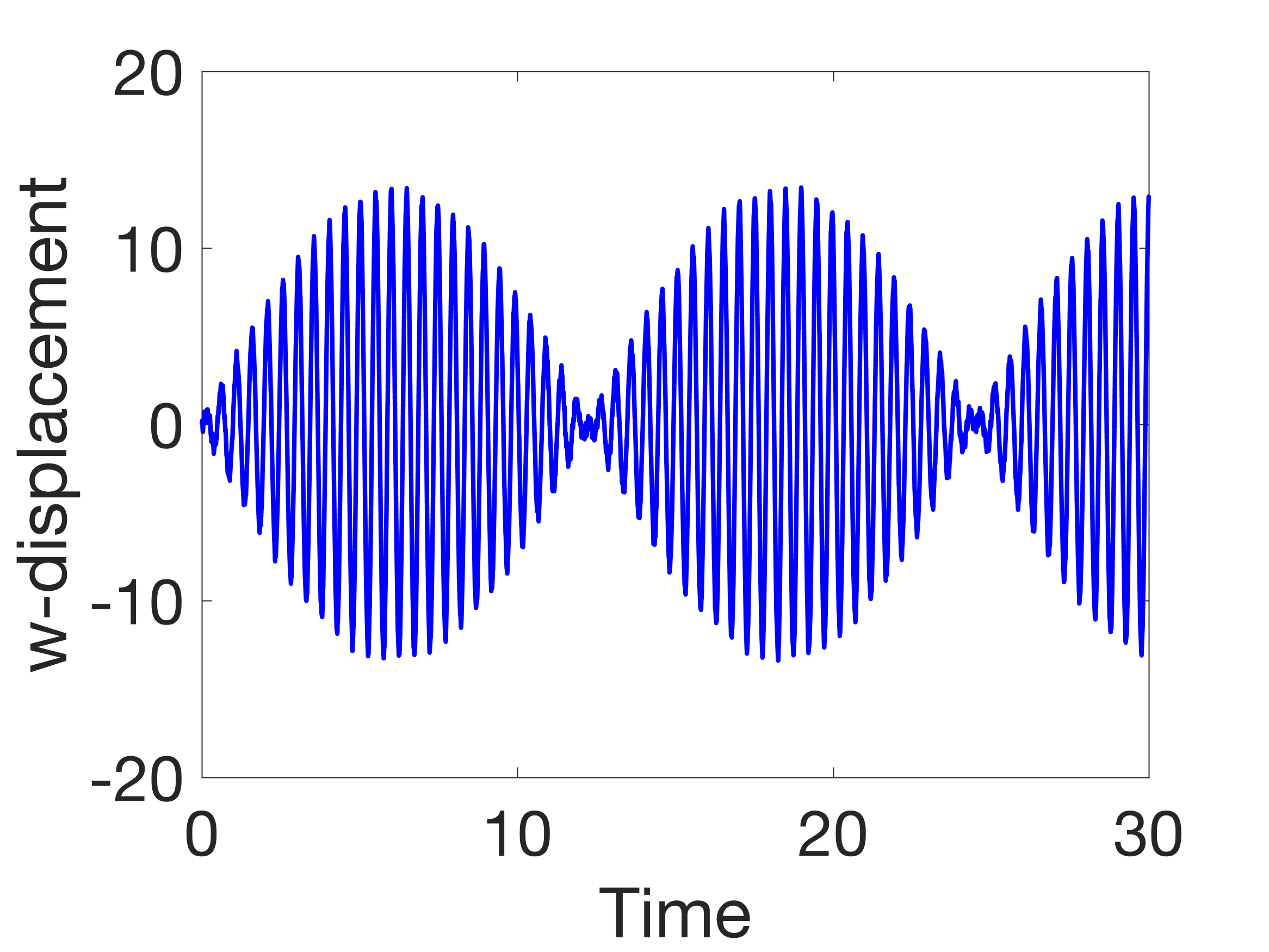}
    \caption{\footnotesize $w(\xv_p,t)$}  \label{fig:MovingClampedBeatB}%
  \end{subfigure} 
    \begin{subfigure}[b]{0.32\linewidth}
    \centering
    \includegraphics[width=1\linewidth]{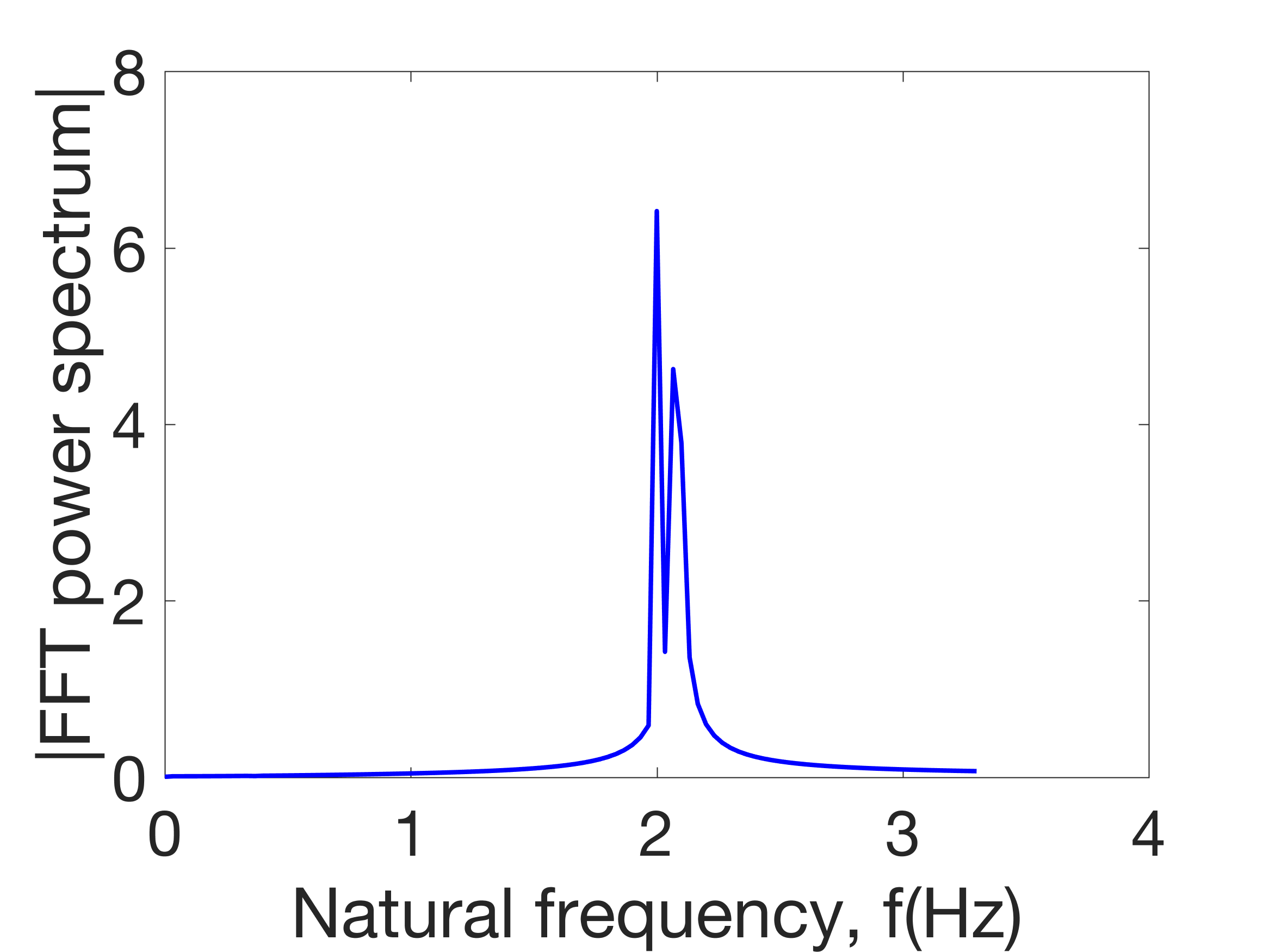}
    \caption{\footnotesize  FFT  of  $w(\xv_p,t)$}   \label{fig:MovingClampedBeatC}%
  \end{subfigure} 
    \caption{
Numerical simulation of the beat phenomenon  using the NB2 method on grid $\Gc_{80}$. Here the driving frequency is $\xi_b=4\pi$ (or $f_b=2$~Hz) that is close to the natural frequency $2\pi f_2$ and the probed location is $\xv_p=(-0.2,0)$.
    }    
  \label{fig:MovingClampedBeat}%
 \end{figure}

   To simulate  beat, we drive the plate  at a frequency that is very close to the previously found  natural frequency $f_2$.  Therefore,  we set the driving frequency in  \eqref{eq:ModifiedClampedBC}  as the so-call beat frequency $\xi_b =4\pi$ (or $f_b=2$~Hz),   noting that the difference between the  beat frequency  and the natural frequency $f_2$ is small for $|f_b- f_2|= 0.067$. Again, the simulation is performed using the NB2 scheme until $t=30$, and a similar collection of results are presented in Figure~\ref{fig:MovingClampedBeat}. The nodal line pattern  for this case (Figure~\ref{fig:MovingClampedBeatA}) is in accordance  with  that shown in Figure~\ref{fig:MovingClampedResonanceA}. Furthermore, Figure~\ref{fig:MovingClampedBeatB} shows the expected oscillation pattern  that resembles the   beat phenomenon.  The FFT power spectrum of the displacement data at $\xv_p$, as is  shown in Figure~\ref{fig:MovingClampedBeatC},  clearly shows the two adjacent   frequencies that correspond respectively to  the beat frequency $f_b$ and the natural frequency  $f_2$.

{
\newcommand{\figWidth}{5.5cm}
\def\xa{12.5}
\def\ya{10}
\newcommand{\trimfig}[2]{\trimw{#1}{#2}{0.08}{0.08}{0.08}{0.06}}
\begin{figure}[h!]
\begin{center}
\begin{tikzpicture}[scale=1]
  \useasboundingbox (0.0,0.0) rectangle (\xa,\ya);  

\draw(3.,5.5) node[anchor=south west,xshift=0pt,yshift=0pt] {\trimfig{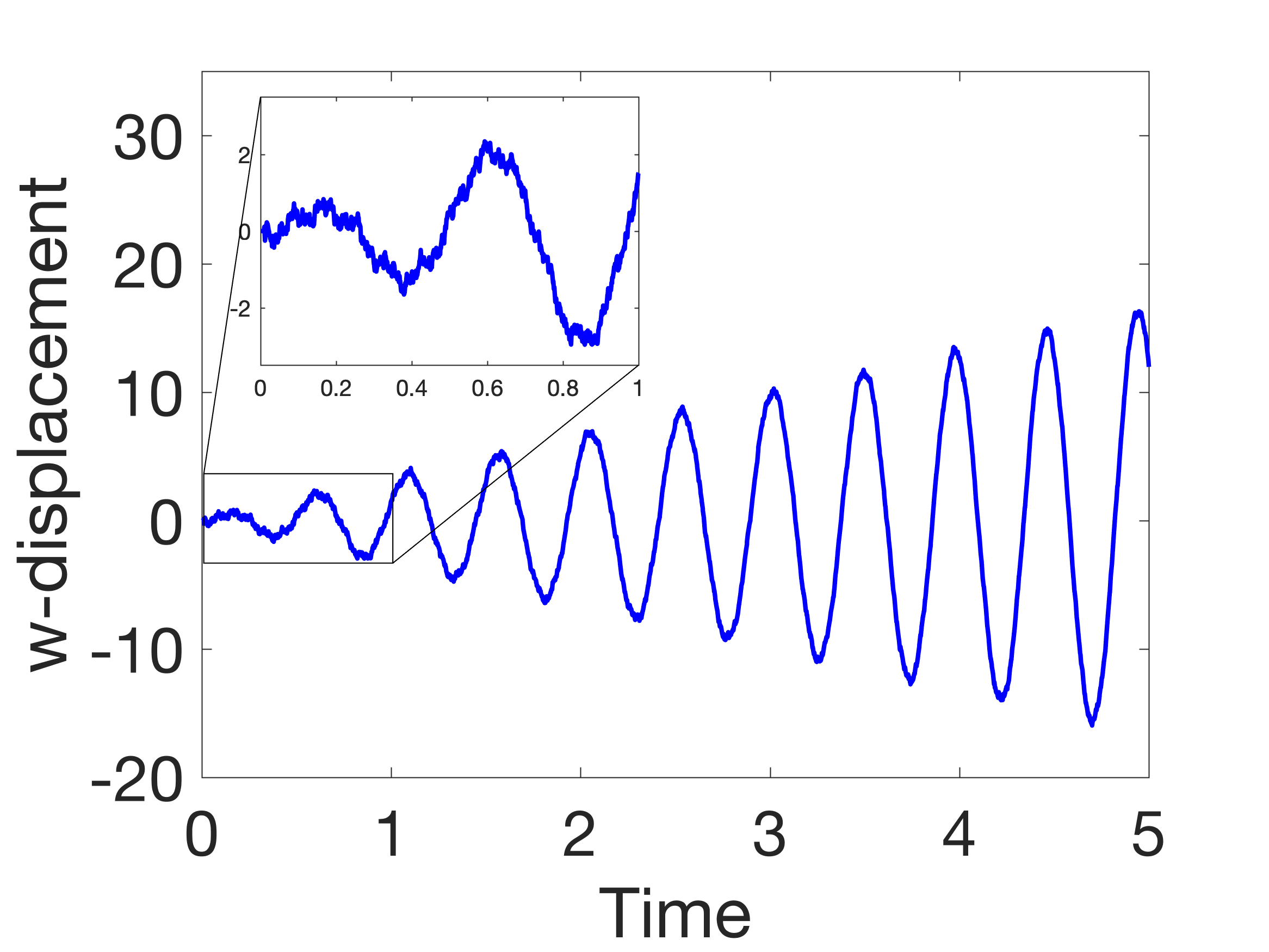}{\figWidth}};
\draw(-.5,0.5) node[anchor=south west,xshift=0pt,yshift=0pt] {\trimfig{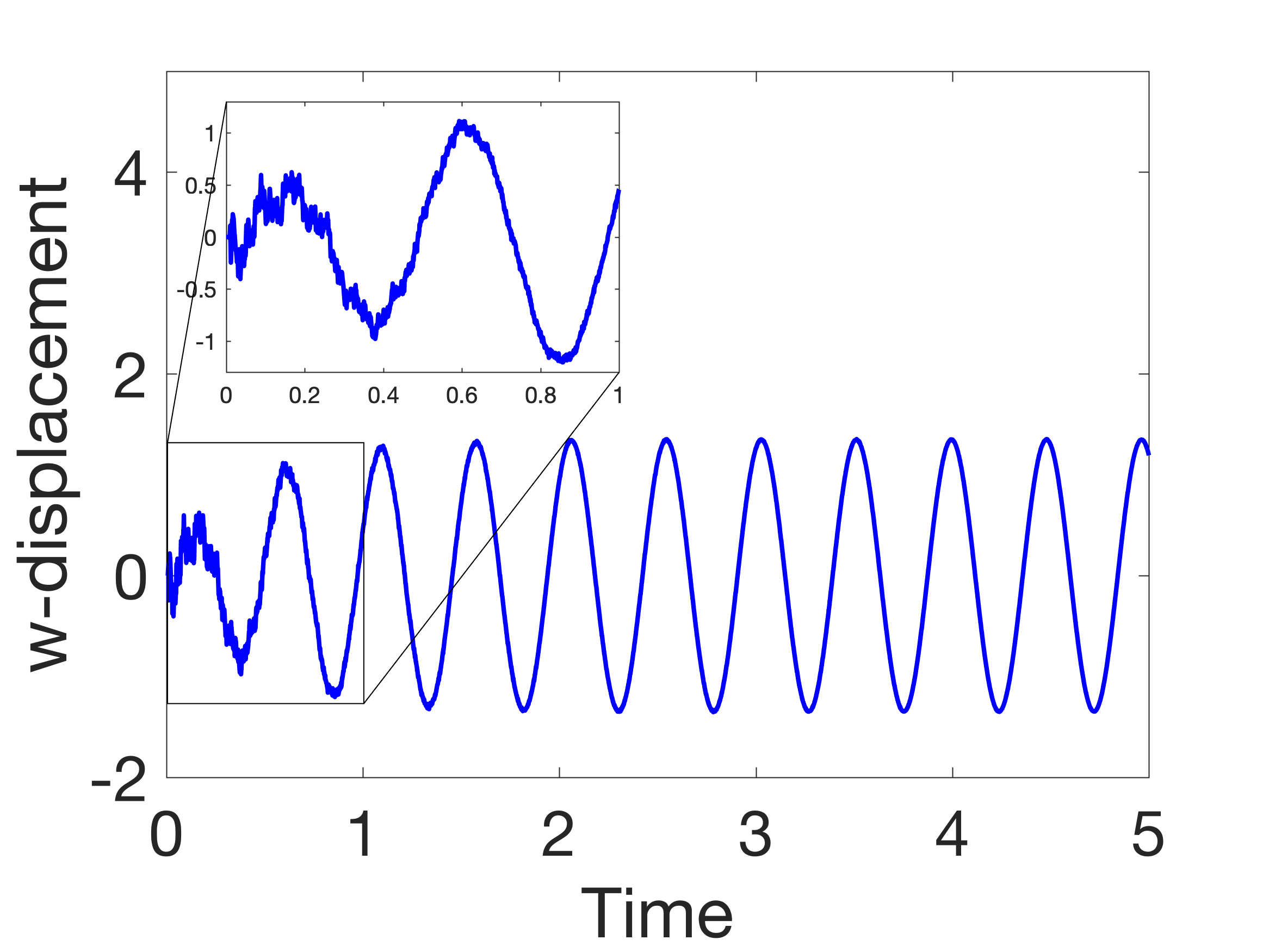}{\figWidth}};
\draw(6.5,0.5) node[anchor=south west,xshift=0pt,yshift=0pt] {\trimfig{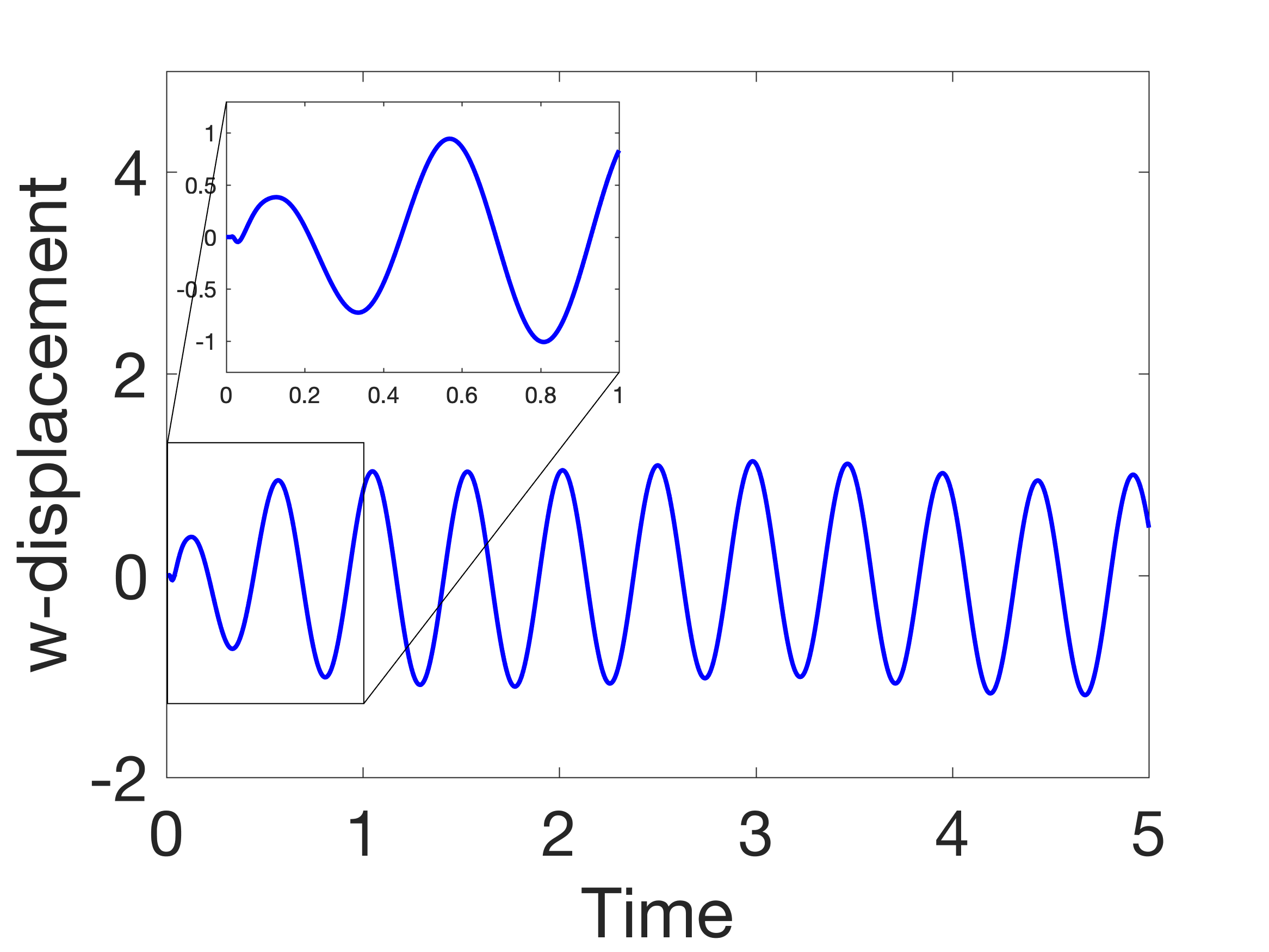}{\figWidth}};

\draw(3,0)  node[anchor=south,xshift=0pt,yshift=0pt] {Damped: $K_1=5, {T_1=0} $};
\draw(10,0)  node[anchor=south,xshift=0pt,yshift=0pt] {Damped : ${K_1=0}, T_1=0.1$};
\draw(6.5,5.)  node[anchor=south,xshift=0pt,yshift=0pt] {Undamped: $K_1=0, T_1=0 $};

%
\end{tikzpicture}

\end{center}
\caption{Plot of  $w(\xv_p,t)$ {\em vs.} $t$ for three different damping scenarios. 
  Simulations are  carried out using the NB2 scheme on grid $\Gc_{80}$.
Here the driving frequency is the same as the resonant frequency $\xi_2=2\pi f_2$ and $\xv_p=(-0.2,0)$.}
  \label{fig:MovingClampedDamped}%
\end{figure}
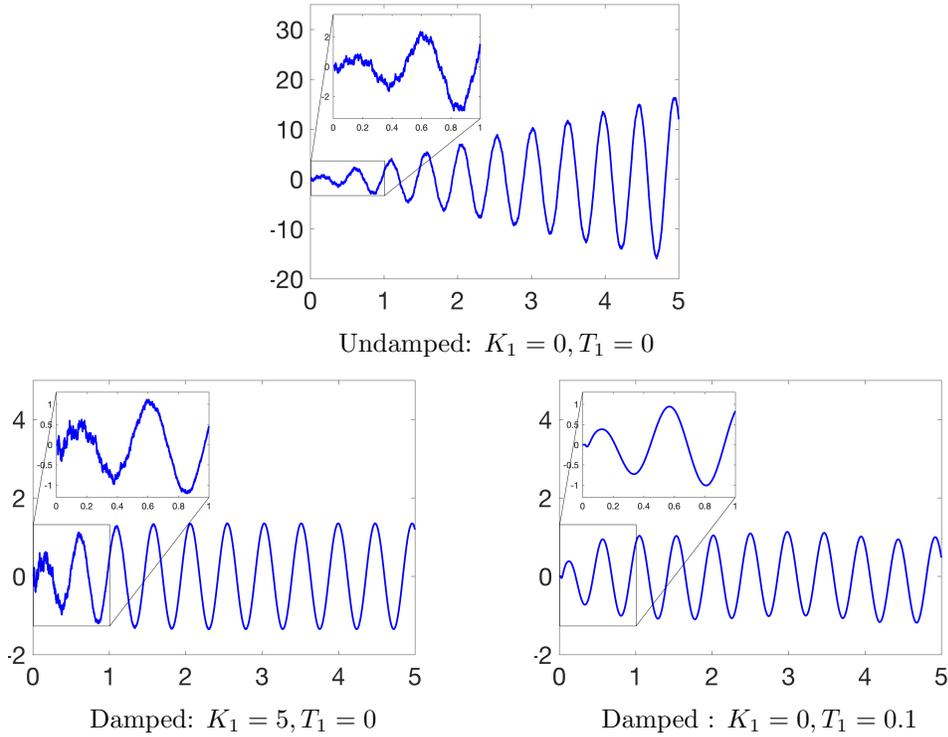
}

Now, we consider the effects of the damping terms in  the generalized  Kirchhoff-Love plate \eqref{eq:generalizedKLPlate}, which are  the  linear damping term with coefficient ${K}_1$  and the  visco-elastic damping term  with  coefficient ${T}_1$. The visco-elastic damping tends to smooth high-frequency oscillations in space and is often added to model vascular structures in haemodynamics. For the simulation, we use the same numerical setup as the resonance example, and  demonstrate the damping effects  by tracking the displacement over time at the point  $\xv_p=(-0.2,0)$   for  three combinations of $K_1$ and $T_1$ values. The time history of the displacements are shown in Figure~\ref{fig:MovingClampedDamped}. For the case when only the linear damping term  $(K_1=5, {T_1=0})$ is added, we can see from the plot that  the amplitude  of the oscillation is damped down to around $2$ when compared with the undamped resonance case. If we zoom in the plot at early times, we do see some  high-frequency oscillations.  However, when  the visco-elastic damping term is  included (${K_1=0}, T_1=0.1$),  we can no longer see  the high-frequency oscillations in the zoomed in plot;  therefore, this term serves to smooth the wave as expected.

Similar numerical results can also be obtained using the PC22 method, although all the simulations are conducted with the NB2 scheme in this section. It is important to  remark  that  the examples considered here   also showcase the accuracy and efficiency of our numerical methods. The fact that we are  able to simulate   the resonance and beat phenomena using a numerically found   natural frequency and that the numerically observed resonance frequency further  corroborates that value strongly suggests the accuracy of our computations.

\section{Conclusions}\label{sec:conclusions}
In this paper, 
 we propose two numerical schemes, referred to as the  PC22 and the NB2 schemes, for the approximation of  a generalized  Kirchhoff-Love plate model. Both schemes are based on centered finite difference methods of second-order accuracy  for  spatial discretization; and the  resulted spatially  discretized 
equations are then advanced  in time using an appropriate time-stepping scheme.   The PC22 scheme uses    an  explicit predictor-corrector scheme that consists of a second-order Adams-Bashforth (AB2) predictor and a second-order
Adams-Moulton (AM2) corrector, and  the NB2 scheme  utilizes an implicit Newmark-Beta  scheme of second-order accuracy.  Stable and accurate numerical boundary conditions are also derived for three  common plate boundary conditions (clamped, simply supported and free).    Stability analysis is performed for the time-stepping schemes to find the regions of absolute stability, which are utilized to determine stable time steps for both proposed schemes.

Carefully designed test problems are solved  to demonstrate the properties and applications of our numerical approaches. The stability and accuracy of the schemes are verified  by mesh  refinement studies using problems with known exact solutions, and  by cross-validation with  experimental results. An interesting application  concerning the exploration of the  resonance and beat phenomena  of an annular plate with general configurations is considered  to further display the accuracy and efficiency of the numerical methods.

The domains of all the examples considered in this paper are restricted to simple ones that can be discretized with a single Cartesian or curvilinear mesh. We would like to extend the schemes for more general geometries using  composite overlapping grids \cite{CGNS}. According to previous studies for wave-like equations on overlapping grids  \cite{max2006b}, we expect weak instabilities to occur  near the interpolation points  of the overlapping grids. Therefore, the investigation  of novel methods,  such as adding high-order spatial  dissipation and   upwind schemes,  to   suppress possible instabilities that can be generated from the overlapping grid interpolation 
would be interesting topics for  future research.

\section*{Acknowledgement}
L. Li is grateful to  Professor W.D. Henshaw of  Rensselaer Polytechnic Institute (RPI) for helpful conversations. Portions of this research were conducted with high performance computational resources provided by the Louisiana Optical Network Infrastructure (http://www.loni.org).

\bibliographystyle{elsart-num}
\bibliography{journal-ISI,LongfeiPapers,henshawPapers,plate}



\clearpage
\begin{appendix}
\section{Nodal line patterns for the eigenvalue problem}\label{sec:appendix}

We show the results of the eigenvalue problem \eqref{eq:eigenProblem} here. Nodal lines of the first $25$ eigenmodes (with multiplicity) for the   square  plate with clamped edges  and the annular plate with simply supported boundaries are shown in  Figures~\ref{fig:ClampedModeShapeSquareMATLAB} \& \ref{fig:SupportedModeShapeAnnMATLAB},  respectively. The eigenmodes plotted for each degenerated pair are arbitrary so they can be asymmetric.

{
\newcommand{\figWidth}{13cm}
\def\xa{13.}
\def\ya{13.}
\newcommand{\trimfig}[2]{\trimw{#1}{#2}{0.02}{0.22}{0.02}{0.0}}
\begin{figure}[h]
\begin{center}
\begin{tikzpicture}[scale=1]
  \useasboundingbox (0.0,0.0) rectangle (\xa,\ya);  

\draw(6.3,6.3) node[anchor=center,xshift=0pt,yshift=0pt] {\trimfig{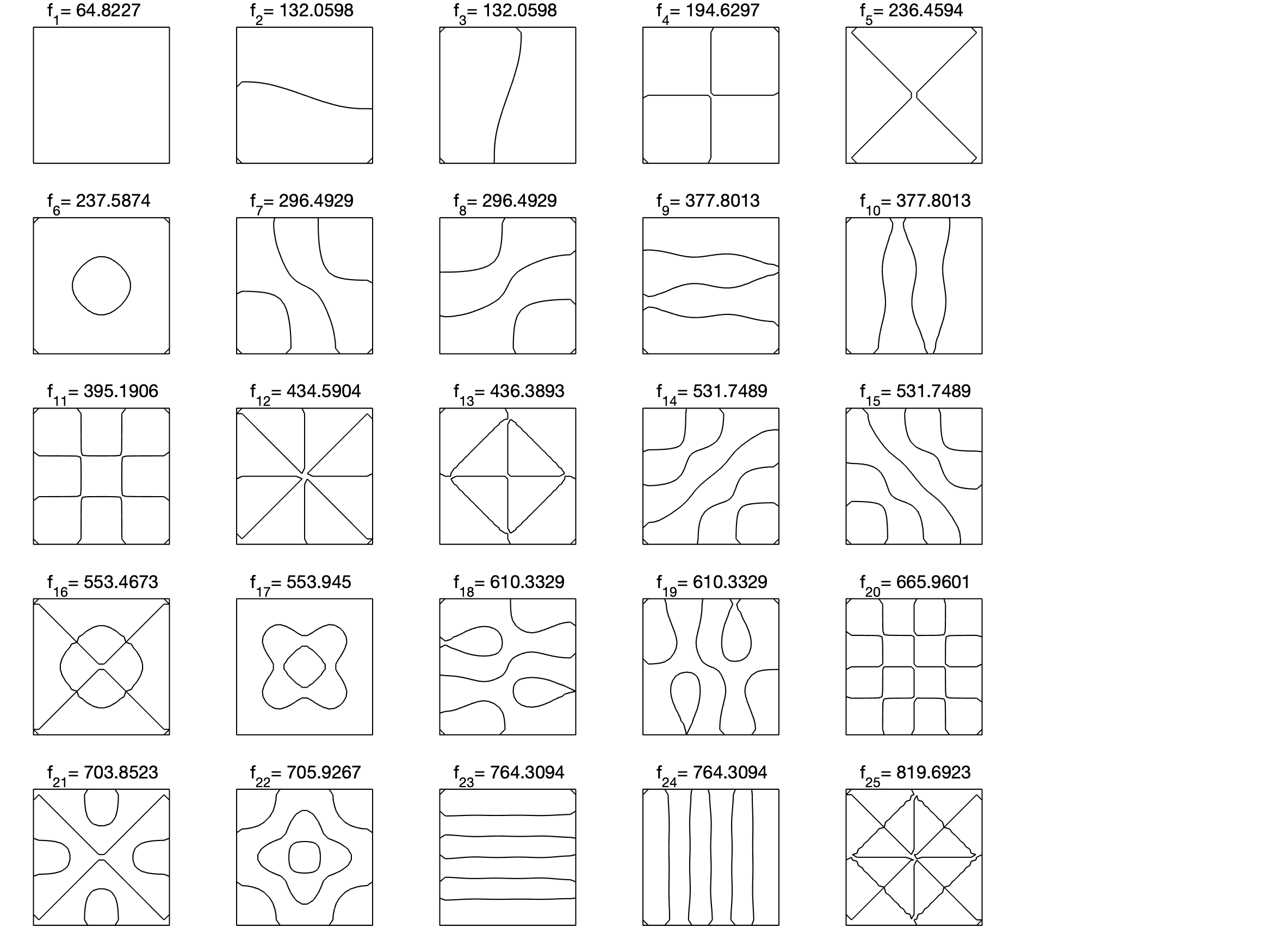}{\figWidth}};

%
\end{tikzpicture}

\end{center}
\caption{Nodal lines of the first $25$ eigenmodes (with multiplicity) of the clamped square plate. There are $6$ degenerate pairs of eigenmodes, so that only $19$ distinct normalized eigenvalues are represented. Values reported in the plots are the natural frequencies in increasing order.
  }
  \label{fig:ClampedModeShapeSquareMATLAB}
\end{figure}
}

{
\newcommand{\figWidth}{13cm}
\def\xa{13.}
\def\ya{13.}
\newcommand{\trimfig}[2]{\trimw{#1}{#2}{0.02}{0.22}{0.02}{0.0}}
\begin{figure}[h]
\begin{center}
\begin{tikzpicture}[scale=1]
  \useasboundingbox (0.0,0.0) rectangle (\xa,\ya);  

\draw(6.3,6.3) node[anchor=center,xshift=0pt,yshift=0pt] {\trimfig{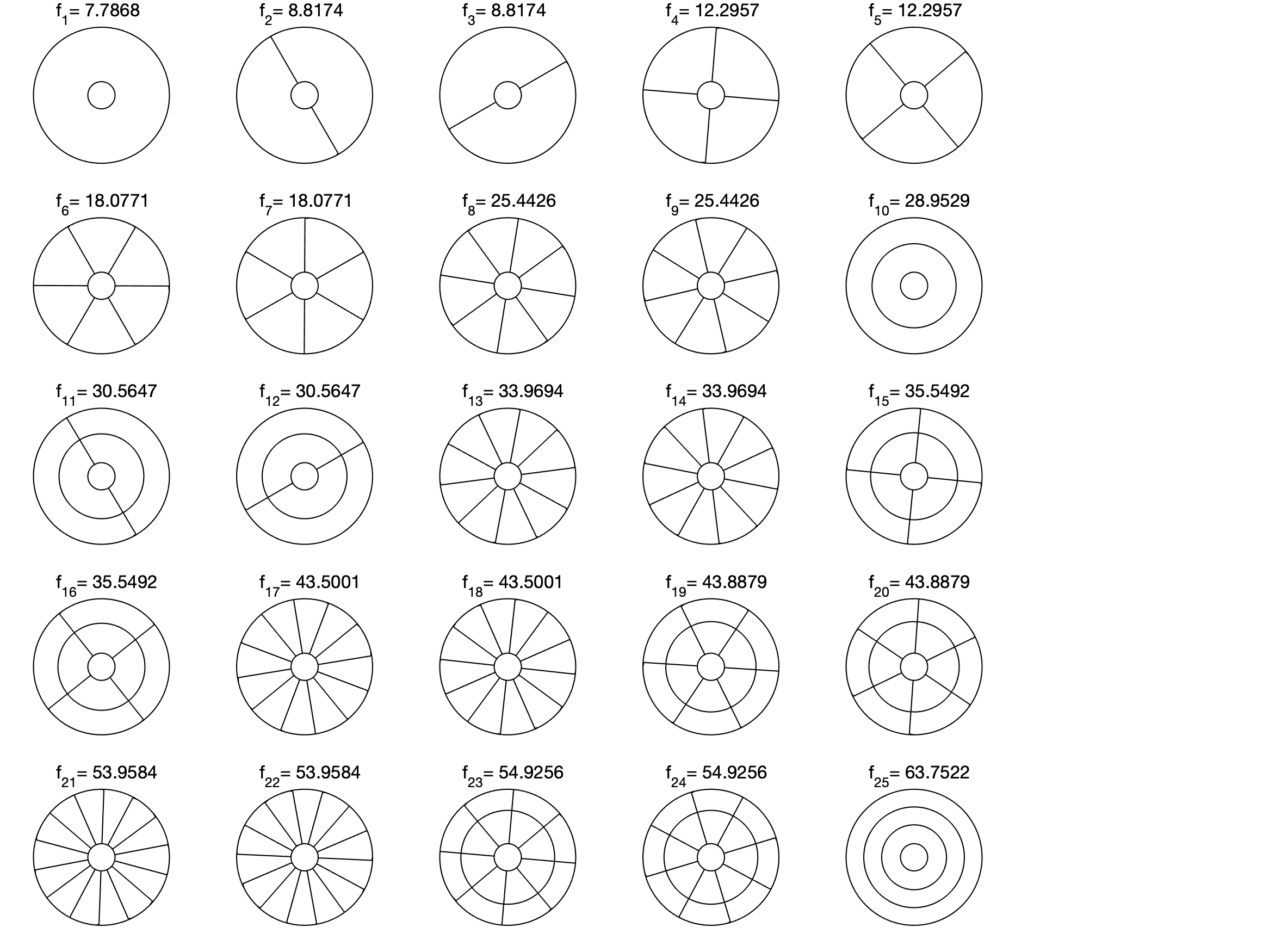}{\figWidth}};

%
\end{tikzpicture}

\end{center}
\caption{ Nodal lines of the first $25$ eigenmodes (with multiplicity) of the simply supported annular plate. There are $11$ degenerate pairs of eigenmodes, so that only $14$ distinct eigenvalues are represented. Values reported in the plots are the natural frequencies in increasing order.
  }
  \label{fig:SupportedModeShapeAnnMATLAB}
\end{figure}
}

  

\end{appendix}


\end{document}